\documentclass[a4paper,10pt,reqno]{amsart}

\usepackage[top=2.5cm, bottom=2.5cm, left=2.5cm, right=2.5cm]{geometry}
\usepackage{amssymb,amsmath,amsthm,ascmac,array,bm,multirow}
\usepackage[dvipdfmx]{graphicx}
\usepackage{color}

\allowdisplaybreaks
\numberwithin{equation}{section}

\newtheorem{thm}{Theorem}[section]

\newtheorem{rem}[thm]{Remark}
\newtheorem{ass}[thm]{Assumption}

\newcommand{\mcb}{\mathcal{B}}
\newcommand{\mcc}{\mathcal{C}}
\newcommand{\mcf}{\mathcal{F}}

\newcommand{\mcl}{\mathcal{L}}
\newcommand{\mcm}{\mathcal{M}}

\newcommand{\mfp}{\mathfrak{p}}
\newcommand{\mbbe}{\mathbb{E}}
\newcommand{\mbbh}{\mathbb{H}}
\newcommand{\mbbn}{\mathbb{N}}
\newcommand{\mbbr}{\mathbb{R}}
\newcommand{\mbbu}{\mathbb{U}}
\newcommand{\mbby}{\mathbb{Y}}
\newcommand{\mbbz}{\mathbb{Z}}
\newcommand{\mbX}{\mathbf{X}}

\newcommand{\del}{\delta}

\newcommand{\sig}{\sigma}
\newcommand{\ep}{\epsilon}
\newcommand{\D}{\Delta}
\newcommand{\Sig}{\Sigma}
\newcommand{\lam}{\lambda}
\newcommand{\Lam}{\Lambda}

\newcommand{\Gam}{\Gamma}
\newcommand{\p}{\partial}

\newcommand{\cil}{\xrightarrow{\mcl}} 
\newcommand{\cip}{\xrightarrow{P}} 

\newcommand{\argmin}{\mathop{\rm argmin}} 
\newcommand{\argmax}{\mathop{\rm argmax}}
\newcommand{\diag}{\mathop{\rm diag}} 
\newcommand{\tr}{\mathop{\rm tr}} 

\def\ds#1{\displaystyle{#1}} 

\def\nn{\nonumber}

\def\sumj{\sum_{j=1}^{n}}
\def\sumk{\sum_{k=1}^{K}}
\def\prk{\prod_{j=k+1}^{K}}

\def\pr{\mathbb{P}}
\def\E{\mathbb{E}}
\def\var{\mathrm{var}}

\def\dim{\mathrm{dim}}

\def\ult{\underline{\theta}}
\def\olt{\overline{\theta}}
\def\tz{\theta_{0}}
\def\tkz{\theta_{k,0}}

\def\tes{\hat{\theta}_{n}}
\def\tet{\tilde{\theta}_{n}}

\def\qbic{\mathrm{QBIC}_{n}}
\def\bic{\mathrm{BIC}_{n}}
\def\kl{\mathrm{KL}}
\def\bfn{\mathrm{BF}_{n}}

\title[Schwarz type model comparison for LAQ models]{Schwarz type model comparison for LAQ models}

\author[S. Eguchi]{Shoichi Eguchi}
\address[Corresponding author]{Graduate School of Mathematics, Kyushu University, 744 Motooka Nishi-ku Fukuoka 819-0395, Japan.}
\email{s-eguchi@math.kyushu-u.ac.jp}

\author[H. Masuda]{Hiroki Masuda}
\address{Faculty of Mathematics, Kyushu University, 744 Motooka Nishi-ku Fukuoka 819-0395, Japan}
\email{hiroki@math.kyushu-u.ac.jp}

\date{\today}
\keywords{Approximate Bayesian model comparison, Gaussian quasi-likelihood, locally asymptotically quadratic family, quasi-likelihood, Schwarz's criterion.}

\begin{document}
\setlength{\baselineskip}{4.5mm}

\maketitle

\begin{abstract}
For model-specification purpose, we study asymptotic behavior of the marginal quasi-log likelihood associated with a family of locally asymptotically quadratic (LAQ) statistical experiments. 
Our result entails a far-reaching extension of applicable scope of the classical approximate Bayesian model comparison due to Schwarz, with frequentist-view theoretical foundation. In particular, the proposed statistics can deal with both ergodic and non-ergodic stochastic-process models, where the corresponding $M$-estimator is of multi-scaling type and the asymptotic quasi-information matrix is random. Focusing on the ergodic diffusion model, we also deduce the consistency of the multistage optimal-model selection where we may select an optimal sub-model structure step by step, so that computational cost can be much reduced.
We illustrate the proposed method by the Gaussian quasi-likelihood for diffusion-type models in details, together with several numerical experiments.
\end{abstract}

\section{Introduction} \label{Intro}

The objective of this paper is Bayesian model comparison for a general class of statistical models, which includes various kinds of stochastic-process models that cannot be handled by preceding results. 
There are two classical principles of model selection: the Kullback-Leibler divergence (KL divergence) principle and the Bayesian one,  acted over Akaike information criterion (AIC, \cite{Aka73} and \cite{Aka74}) and Schwarz or Bayesian information criterion (BIC, \cite{Sch78}), respectively. A common knowledge is that there are no universal politic between AIC and BIC type statistics, and they are indeed used for different purposes. 
On the one hand, the AIC is a predictive model-selection criterion minimizing the KL divergence between prediction and true models, 
not intended to pick up the true model consistently even if it does exist in the candidate-model set. 
On the other hand, the BIC is used to look for better model description. 
The BIC for i.i.d. model is given by
\begin{equation}
\mathrm{BIC}_{n} = -2\ell_{n}(\hat{\theta}_{n}^{\mathrm{MLE}})+p\log n,
\nn
\end{equation}
where $\ell_{n}$, $\hat{\theta}_{n}^{\mathrm{MLE}}$, and $p$ denote 
the log-likelihood function, the maximum-likelihood estimator, and the dimension of the parameter space of the statistical model to be assessed, respectively. 
The model selection consistency via BIC type statistics has been studied by many authors in several different model setups, 
e.g. \cite{Boz87}, \cite{CasGirMarMor09}, and \cite{Nis84}, to mention just a few old ones. 

There also do exist many studies of the BIC methodology in the time series context. 
The underlying principles, such as maximization of posterior-selection probability, remain the same in this case. 
It should be mentioned is that \cite{CavNea99} demonstrated that derivation of the classical BIC could be generalized into 
general $\sqrt{n}$-consistent framework with constant asymptotic information. 
Their argument supposes the almost-sure behaviors of the likelihood characteristics, especially of the observed information matrix. 
Our stance is similar to theirs, but more general so as to subsume a wide spectrum of models that cannot be handled by \cite{CavNea99}. 
However, much less is known about theoretically guaranteed information criteria concerning sampled data from stochastic-process models; 
to mention some of them, we refer to \cite{Uch10}, \cite{UchYos01}, \cite{UchYos06}, \cite{UchYos16}, and also \cite{SeiKom07}.

Our primary interest is to extend the range of application of Schwarz's BIC to a large degree in a unified way, so as to be able to target a wide class of dependent-data models especially including 
the locally asymptotically mixed-normal family of statistical experiments. One may state the Bayesian principle of model selection among $\mcm_{1},\dots,\mcm_{M}$ is to choose the model that is most likely in terms of the posterior probability, which is typically done in terms of the logarithmic Bayes factor and the Kullback-Leibler divergence, through approximating the (expected) marginal quasi-log likelihood. 
Unfortunately, a mathematically-rigorous derivation of BIC type statistics is sometimes missing in the literature, especially when underlying model is non-ergodic. 
In this paper, we will focus on {\it locally asymptotically quadratic (LAQ)} statistical models, where asymptotic Fisher information matrix being possibly random. We will introduce the quasi BIC (QBIC) through the {\it stochastic} expansion of the {\it marginal quasi-likelihood}. Here, we use the terminology ``quasi'' to mean that the model may be misspecified in the sense that the candidate models may not include the true joint distribution of data sequence; see \cite{LvLiu14} for information criteria for a class of generalized linear models for independent data.
Our proof of the expansion essentially utilizes the polynomial type large deviation inequality (PLDI) of \cite{Yos11}; quite importantly, the asymptotic information matrix then may be random (i.e. suitably scaled observed information (random bilinear form) has a random limit in probability), especially enabling us to deal with non-ergodic models; indeed, random limiting information is quite common in the so-called non-ergodic statistics. Since a quasi-likelihood may be used even for semiparametric models where possibly infinite-dimensional nuisance element, so does our QBIC. 

Our asymptotic results are widely applicable enough to provide a unified way to deduce validity of approximate model assessment via BIC type statistics, and to cover a broad spectrum of quasi-likelihoods associated with dependent data models having the LAQ structure. 
Though we do not go into any detail here, the popular cointegration models (see \cite{Bos10} and the references therein) would be in the scope of QBIC as well. In particular, \cite{Kim98} clarified the ``correct-BIC'' form in the context of non-stationary time series models, where the observed information matrix is involved in the bias-correction term, resulting in an extension of the classical BIC.

We note that there are many other works on the model selection include the risk information criterion \cite{FosGeo94}, the generalized information criterion \cite{KonKit96}, the ``parametricness'' index \cite{LiuYan11}, and many extensions of AIC and BIC including \cite{CheChe08} and \cite{LvLiu14}.
We refer to \cite{BurAnd02}, \cite{ClaHjo08}, and \cite{KonKit08} for comprehensive accounts of information criteria, and also to \cite{DziCofLanLi12} for an illustration from practical point of view. 

\medskip

This paper is organized as follows. 
Section \ref{hm:sec_pre} describes the basic model setup and some related backgrounds.
In Section \ref{hm:sec_qbic} we will present the asymptotic expansions of the expected Bayes factor (equivalently, the marginal quasi-log likelihood and the Kullback-Leibler divergence), which provides us with a unified way to quantitative model comparison in a model-descriptive sense; the presentation contains a revised and extended version of \cite{EguMas15}.
In Section \ref{hm:sec_GQNLE}, we illustrate the proposed model selection method by the Gaussian quasi-likelihood estimation of ergodic diffusion process and volatility-parameter estimation for a class of continuous semimartingales, both based on high-frequency sampling; to the best of our knowledge, this is the first place that mathematically validates Schwarz's methodology of model comparison for high-frequency data from a stochastic process.
Section \ref{hm:qbic.msc} is devoted to the model-selection consistency with respect to the optimal model, which is naturally defined to be any candidate model minimizing the quasi-entropy associated with candidate quasi-likelihoods. When in particular the quasi-maximum likelihood estimator is of multi-scaling type, we prove the consistency of the multistage optimal-model selection procedure, where we partially select an optimal model structure step by step, resulting in a reduced computational cost; also, we will briefly mention how the argument can be extended to the general model setting.
Section \ref{Simu} give some numerical experiments supporting our asymptotic results. 
All the proofs are presented in Section \ref{hm:sec_Proofs}.

\section{Preliminaries}\label{hm:sec_pre}

\subsection{Basic model setup}

We begin with describing our basic Bayesian-model setup used throughout this paper. 
Denote by $\mbX_{n}$ an observation random variable defined on an underlying probability space $(\Omega,\mcf,\pr)$, 
and by $G_{n}(dx)=g_{n}(x)\mu_{n}(dx)$ the true distribution $\mcl(\mbX_{n})$, 
where $\mu_{n}$ is a $\sig$-finite dominating measure on a Borel state space of $\mbX_{n}$, 
that is, $G_{n}(dx)=\pr\circ\mbX_{n}^{-1}(dx)$.

Suppose that we are given a set of $M$ candidate Bayesian models $\mcm_{1},\dots,\mcm_{M}$:
\begin{equation}
\mcm_{m}=\big\{\left(\mfp_{m},\,\pi_{m,n}(\theta_{m}),
\,\mbbh_{m,n}(\cdot|\theta_{m})\right)\big|~\theta_{m}\in\Theta_{m}\big\},\qquad m=1,\dots,M,
\nn
\end{equation}
where the ingredients in each $\mcm_{m}$ are given as follows.

\begin{itemize}
\item $\mfp_{m}>0$ denotes the relative likeliness of the model-$\mcm_{m}$ occurrence among $\mcm_{1},\dots,\mcm_{M}$; 
we have $\sum_{m=1}^{M}\mfp_{m}=1$.
\item $\pi_{m,n}: \Theta_{m}\to(0,\infty)$ is the prior distribution $\mcl(\theta_{m})$ of $m$th-model parameter $\theta_{m}$, 
here defined to be a probability density function possibly depending on the sample size $n$, 
with respect to the Lebesgue density on a bounded convex domain $\Theta_{m}\subset\mbbr^{p_{m}}$.
\item The measurable function $x\mapsto \mbbh_{m,n}(x | \theta_{m})$ for each $\theta_{m}\in\Theta_{m}$ defines 
a logarithmic regular conditional probability density of $\mcl(\mbX_{n}|\theta_{m})$ with respect to $\mu_{n}(dx)$.
\end{itemize}
Each $\mcm_{m}$ may be misspecified in the sense that 
the true data generating model $g_{n}(x)$ does not belong to the family $\{\exp\{\mbbh_{m,n}(\cdot | \theta_{m})\}|\,\theta_{m}\in\Theta_{m}\}$; 
we will, however, assume suitable regularity conditions for the associated statistical random fields to have a suitable asymptotic behavior.

Concerning the model $\mcm_{m}$, the random function $\theta_{m}\mapsto \exp\{\mbbh_{m,n}(\mbX_{n} | \theta_{m})\}$, 
assumed to be a.s. well-defined, is referred to as the \textit{quasi-likelihood} of $\mcl(\mbX_{n}|\theta_{m})$. 
The {\it quasi-maximum likelihood estimator (QMLE)} $\hat{\theta}_{m,n}$ associated with $\mbbh_{m,n}$ is defined to be any maximizer of $\mbbh_{m,n}$:
\begin{equation}
\hat{\theta}_{m,n}\in\argmax_{\theta\in\bar{\Theta}_{m}}\mbbh_{m,n}(\theta).
\nn
\end{equation}
For simplicity 
we will assume the a.s. continuity of $\mbbh_{m,n}$ over the compact set $\bar{\Theta}_{m}$, so that $\hat{\theta}_{m,n}$ always exists.

\medskip

Our objective includes estimators of multi-scaling type, meaning that 
the components of $\hat{\theta}_{m,n}$ converges at different rates, which can often occur when considering high-frequency asymptotics. 
A typical example is the Gaussian quasi-likelihood estimation of ergodic diffusion process: see \cite{Kes97}, also Section \ref{GQMLE2}. 
Let $K_{m}\in\mbbn$ be a given number, which represents the number of the components having different convergence rates in $\mcm_{m}$, 
and assume that the $m$th-model parameter vector is divided into $K_{m}$ parts: 
\begin{equation}
\theta_{m}=(\theta_{m,1},\ldots,\theta_{m,K_{m}})\in \prod_{k=1}^{K_{m}}\Theta_{m,k}=\Theta_{m},
\nonumber
\end{equation}
with each $\Theta_{m,k}$ being a bounded convex domain in $\mbbr^{p_{m,k}}$, $k\in\{1,\ldots,K_{m}\}$, where $p_{m}=\sum_{k=1}^{K_{m}} p_{m,k}$. 
Then the QMLE in the $m$th model takes the form $\hat{\theta}_{m,n}=(\hat{\theta}_{m,1,n},\ldots,\hat{\theta}_{m,K_{m},n})$. 
The optimal value of $\theta_{m}$ associated with $\mbbh_{m,n}$, to be precisely defined later on, 
is denoted by $\theta_{m,0}=(\theta_{m,1,0},\ldots,\theta_{m,K_{m},0})$, $\theta_{m,k,0}\in\Theta_{m,k}$.
The rate matrix in the model $\mcm_{m}$ is then given in the form
\begin{equation}
A_{m,n}(\theta_{m,0})=\diag\left( a_{m,1,n}(\theta_{m,0})I_{p_{m,1}},\ldots,a_{m,K_{m},n}(\theta_{m,0})I_{p_{m,K_{m}}}\right),
\label{hm:nmat}
\end{equation}
where $I_{p}$ denotes the $p$-dimensional identity matrix and $a_{m,k,n}$ $(\theta_{m,0})$ is a deterministic sequence satisfying that 
\begin{equation}
a_{m,k,n}(\theta_{m,0})\to 0,\qquad a_{m,i,n}(\theta_{m,0}) / a_{m,j,n}(\theta_{m,0}) \to 0\quad(i<j),\qquad n\to\infty.
\label{hm:nmat.cond}
\end{equation}
The diagonality of $A_{m,n}(\tz)$ here is just for simplicity.

\medskip

Since we are allowing not only data dependency but also possibility of the model misspecification, 
we may deal with a wide range of quasi-likelihood $\mbbh_{m,n}$, even including semiparametric situations 
such as the Gaussian quasi-likelihood; see Section \ref{hm:sec_GQNLE} for related models.

\subsection{Bayesian model-selection principle}

The quasi-marginal distribution of $\mbX_{n}$ in the $m$th model $\mcm_{m}$ is given by the density
\begin{equation}
x\mapsto f_{m,n}(x):= \int_{\Theta_{m}}\exp\{\mbbh_{m,n}(x|\theta_{m})\}\pi_{m,n}(\theta_{m})d\theta_{m},
\nonumber
\end{equation}
which is sometimes referred to as the model evidence of $\mcm_{i}$. 
Typical reasoning in Bayesian principle of model selection in $\mcm_{1},\dots,\mcm_{M}$ is 
to choose the model that is most likely to occur in terms of 
the posterior model-selection probability, namely to choose the model maximizing
\begin{equation}
\log\bigg(\frac{f_{m,n}(x)\mfp_{m,n}}{\sum_{i=1}^{M}f_{i,n}(x)\mfp_{i,n}}\bigg)
= \log f_{m,n}(x) + \log\mfp_{m} - \log\bigg(\sum_{i=1}^{M}f_{i,n}(x)\mfp_{i}\bigg)
\nonumber
\end{equation}
over $m=1,\dots,M$. This is equivalent to find
\begin{equation}
m_{0} := \argmax_{m\le M}\left\{ \log f_{m,n}(x) + \log\mfp_{m} \right\}.
\nonumber
\end{equation}
Then one proceeds with suitable almost-sure ($\Omega\ni\omega$-wise) 
asymptotic expansion of the logarithm of the quasi-marginal likelihood $\log f_{m,n}(x)$ for $n\to\infty$ 
around a suitable estimator, a measurable function of $x=x_{n}$ for each $n$: 
when $\sqrt{n}(\hat{\theta}_{m,n}-\theta_{m,0})=O_{p}(1)$, the resulting form may be quite often given by
\begin{equation} 
\log f_{m,n}(x) + \log\mfp_{m} \approx \mbbh_{m,n}(\hat{\theta}_{m,n}) - \frac{p_{m}}{2}\log n + O(1)\quad \text{a.s.}
\label{hm:qm.BIC.expa}
\end{equation}
This is the usual scenario of derivation of the classical-BIC type statistics; 
see \cite{CavNea99} and \cite{LvLiu14} as well as \cite{Sch78}.

\medskip

We recall that the expansion \eqref{hm:qm.BIC.expa} is also used to approximate the Bayes factor. 
The logarithmic Bayes factor of $\mcm_{i}$ against $\mcm_{j}$ is defined by the (random) ratio of posterior and prior odds of model-selection probabilities: 
letting $\pr(\mcm_{i}|\mbX_{n})$ denote the posterior model-selection probability of the $i$th model, we have
\begin{equation}
\log\bfn(i,j) := \log\frac{\pr(\mcm_{i}|\mbX_{n})/\pr(\mcm_{j}|\mbX_{n})}{\mfp_{i}/\mfp_{j}}
=\log\frac{f_{i,n}(\mbX_{n})}{f_{j,n}(\mbX_{n})}.
\label{hm:logBFij}
\end{equation}
The Bayes factor measures ``gain'' in model-selection odds between $\mcm_{i}$ and $\mcm_{j}$ 
when observing $\mbX_{n}$. We relatively prefer $\mcm_{i}$ to $\mcm_{j}$ if $\log\bfn(i,j)>0$, and vice versa. 
A selected model via the Bayes factor minimizes the total error rates compounding false-positive and false-negative probabilities, while, different from the AIC, it has no theoretical implication for predictive performance of the selected model. For a more detailed account of the philosophy of the Bayes factor, we refer to \cite{LavSch99}.

\medskip

As was explained in \cite{LvLiu14}, we have yet another interpretation based on the Kullback-Leibler (KL) divergence 
between the true distribution $g_{n}$ and the $m$th quasi-marginal distribution $f_{m,n}$:
\begin{align}
\kl(f_{m,n};g_{n})&:=
-\int\bigg(\log\frac{f_{m,n}(x)}{g_{n}(x)}\bigg)g_{n}(x)\mu_{n}(dx)
\nn\\
&=\int\{\log g_{n}(x)\}g_{n}(x)\mu_{n}(dx) - \int\{\log f_{m,n}(x)\}g_{n}(x)\mu_{n}(dx);
\label{hm:KL}
\end{align}
recall that in the classical-AIC methodology we instead look at $\kl\{f_{m,n}(\cdot;\hat{\theta}_{m,n});g_{n}\}$ 
where $\hat{\theta}_{m,n}=\hat{\theta}_{m,n}(\tilde{\mbX}_{n})$ denotes the MLE in the $m$th \textit{correctly} specified model, constructed from an i.i.d. copy $\tilde{\mbX}_{n}$ of $\mbX_{n}$. 
Based on \eqref{hm:KL}, we regard the $m_{0}$th model relatively optimal among $\mcm_{1},\dots,\mcm_{M}$ where
\begin{equation}
\argmin_{m\le M}\kl(f_{m,n};g_{n}) = \argmax_{m\le M}\int\{\log f_{m,n}(x)\}g_{n}(x)\mu_{n}(dx).
\nonumber
\end{equation}
Comparison of $f_{i,n}$ and $f_{j,n}$ is equivalent to look at the sign of
\begin{equation}
\kl(f_{j,n};g_{n}) - \kl(f_{i,n};g_{n}) 
= \int\log\bigg(\frac{f_{i,n}(x)}{f_{j,n}(x)}\bigg)g_{n}(x)\mu_{n}(dx)
= \E\bigg\{\log\bigg(\frac{f_{i,n}(\mbX_{n})}{f_{j,n}(\mbX_{n})}\bigg)\bigg\}.
\label{hm:KL.diff}
\end{equation}
As was noted in \cite{LvLiu14}, it is important to notice that this reasoning remains valid even when any of candidate models does not coincide with the true model. 
We also refer to \cite{GouRob98} for another Bayesian variable selection device based on the KL projection. 

We will introduce a set of regularity conditions under which explicit statistics $\qbic^{\sharp,i}$ for each model $\mcm_{i}$, $i=1,\dots,M$ (see \eqref{hm:Fhat_def} below) satisfies the stochastic expansion
\begin{equation}
\log\bfn(i,j) = \frac{1}{2}(\qbic^{\sharp,j}-\qbic^{\sharp,i}) + o_{p}(1).
\nonumber
\end{equation}
In the classical treatment originating \cite{Sch78}, the almost-sure expansion was relevant, see Remark \ref{hm:rem_nrX}.

\subsection{Expected logarithmic Bayes factor}\label{hm:sec_ebf}

Comparing \eqref{hm:KL.diff} with \eqref{hm:logBFij}, we see that the expected Bayes factor is directly related to the KL divergence difference:
\begin{equation}
\E\left\{\log\bfn(i,j)\right\} = \kl(f_{j,n};g_{n}) - \kl(f_{i,n};g_{n}).
\nonumber
\end{equation}
For Bayesian model comparison in our model setting, 
we wish to estimate the quantity $\E\left\{\log\bfn(i,j)\right\}$ for each pair $(i,j)$. 
In Section \ref{hm:qbic.mc} we will derive statistics $\qbic^{\sharp,1}, \dots, \qbic^{\sharp,M}$ such that for each $i,j\in\{1,\dots,M\}$ with $i\ne j$:
\begin{equation}
\E\bigg(\bigg|\log\bfn(i,j) - \frac{1}{2}(\qbic^{\sharp,j}-\qbic^{\sharp,i})\bigg| \bigg) = o(1).
\label{hm:lfb.L1conv}
\end{equation}
In particular, it follows that the statistics $(\qbic^{\sharp,j}-\qbic^{\sharp,i})/2$ 
serves as an asymptotically unbiased estimator of the expected Bayes factor:
\begin{equation}
\E\bigg( \E\left\{\log\bfn(i,j)\right\} - \frac{1}{2}(\qbic^{\sharp,j}-\qbic^{\sharp,i}) \bigg) = o(1),
\nonumber
\end{equation}
or, of the raw (random) Bayes factor:
\begin{equation}
\E\bigg( \log\bfn(i,j) - \frac{1}{2}(\qbic^{\sharp,j}-\qbic^{\sharp,i}) \bigg) = o(1).
\nonumber
\end{equation}
Obviously it suffices for \eqref{hm:lfb.L1conv} to show that
\begin{equation}
\E\left\{\left|\qbic^{\sharp,i} - (-2\log f_{i,n}(\mbX_{n})) \right|\right\}=o(1)
\nonumber
\end{equation}
for each $i$ (see Theorem \ref{hm:thm_mc}). 
We are thus led to the basic rule about an optimal model $\mcm_{m_{0}}$ in the sense of approximate Bayesian model description:
\begin{equation}
m_{0}\in\argmin_{1\le i\le M}\qbic^{\sharp,i}.
\nonumber
\end{equation}

\medskip

\begin{rem}{\rm 
To perform a model comparison based on Bayesian prediction, we should replace the marginal likelihood 
$f_{m,n}$ in \eqref{hm:KL} by a Bayesian predictive model and also the ``$g_{n}d\mu_{n}$''-integral by suitable one. 
We refer to \cite{VehOja12} for an extensive review of Bayesian prediction, 
and also to \cite{SeiKom07} for a study in this direction for the LAMN models.
}\qed\end{rem}

\begin{rem}{\rm 
While we here focus on the finite model comparison among $M$ candidates, it would be possible to consider a \textit{continuum} of models, say $(\mcm_{\lam})_{\in\Lambda}$ for some (possibly uncountable) model-index set $\Lam$. This is relevant when considering continuous fine-tuning in regularization methods, e.g. \cite{BicLi06}. 
Although we do not treat such a setting here, it is readily expected that our claims remain valid in an analogous form.
}\qed\end{rem}

\section{Quasi-Bayesian information criterion}\label{hm:sec_qbic}

We here focus on a single model $\mcm_{m}$ and consider asymptotic expansion related to $\mbbh_{m,n}$. 
From now on we will omit the model index ``$m$'' from the notation without mention, 
simply denoting the prior density and the quasi-log likelihood by $\pi_{n}(\theta)$ and $\mbbh_{n}(\theta)$, respectively. 
The parameter $\theta\in\Theta\subset\mbbr^{p}$ is graded into $K$ parts, say
\begin{equation}
\theta=(\theta_{1},\dots,\theta_{K}),\qquad \theta_{k}\in\mbbr^{p_{k}}.
\nonumber
\end{equation}
Here we wrote $p=\sum_{k=1}^{K}p_{k}$. Let $\tz\in\Theta$ be a constant. We are thinking of situations where 
the contrast function $\mbbh_{n}$ provides an $M$-estimator $\tes$ such that the $A_{n}(\tz)(\tes-\tz)$ tends in distribution to a non-trivial asymptotic distribution.
The rate matrix $A_{n}(\tz)$ is of the form \eqref{hm:nmat} satisfying \eqref{hm:nmat.cond}:
\begin{equation}
A_{n}(\tz)=\diag\left( a_{1,n}(\theta_{0})I_{p_{1}},\dots,a_{K,n}(\theta_{0})I_{p_{K}}\right),
\nonumber
\end{equation}
where positive decreasing sequences $a_{k,n}(\theta_{0})$ such that $a_{k,n}^{-1}(\theta_{0})/a_{l,n}^{-1}(\theta_{0})\to 0$ for $k>l$; 
under a suitable continuity condition on $\theta\mapsto A_{n}(\theta)$, we have $\log|A_{n}(\tes)|=\sum_{k=1}^{K}p_{k}\log a_{k,n}(\hat{\theta}_{n})$; here and in what follows, $|A|:=\det(A)$ for a square matrix $A$. 
The statistical random field associated with $\mbbh_{n}$ is given by
\begin{equation}
\mbbz_{n}(u)=\mbbz_{n}(u;\tz):=\exp\left\{\mbbh_{n}(\theta_{0}+A_{n}(\tz)u)-\mbbh_{n}(\theta_{0})\right\},
\label{Z_def}
\end{equation}
which is defined on the admissible domain $\mbbu_{n}(\tz):=\{u\in\mbbr^{p};\tz+A_{n}(\tz)u\in\Theta\}$. 
The objective here is to deduce the asymptotic behavior of the marginal quasi-log likelihood function
\begin{equation}
f_{n}(\mbX_{n})=\log\bigg(\int_{\Theta}\exp\{\mbbh_{n}(\theta)\}\pi_{n}(\theta)d\theta\bigg),
\nn
\end{equation}
and then derive an extension of the classical BIC.

\subsection{Stochastic expansion}

We begin with the stochastic expansion of the marginal quasi-log likelihood function. 
The polynomial type large deviation estimates \cite{Yos11}, mentioned in the introduction, is a powerful tool 
for ensuring the $L^{q}(\pr)$-boundedness of scaled $M$- and Bayes estimators stemming from the quasi-likelihood $\mbbh_{n}$. 
As seen below, the PLDI argument can be effectively used also to verify the key Laplace-approximation type argument in a unified manner.

\medskip

Write $\p_{\theta}=\p/\p\theta$, and denote by $\theta^{j}$ the $j$th element of $\theta$ and by $A_{n,ii}(\tz)$ the $(i,i)$th element of $A_{n}(\tz)$ 
(i.e., $A_{n,ii}(\tz)=a_{j,n}(\theta_{0})$ for some $j\in\{1,\dots,K\}$).

\begin{ass}
$\mbbh_{n}(\theta)$ is of class $\mcc^3(\Theta)$ and satisfies the following conditions:
\begin{itemize}
\item[(i)] $\ds{\D_{n}=\D_{n}(\tz):=A_{n}(\theta_{0})\p_{\theta}\mbbh_{n}(\theta_{0})=O_{p}(1)}$;
\item[(ii)] $\ds{\Gam_{n}=\Gam_{n}(\tz):=-A_{n}(\theta_{0})\p_{\theta}^{2}\mbbh_{n}(\theta_{0})A_{n}(\theta_{0})=\Gam_{0}+o_{p}(1)}$ where $\pr(\Gam_{0}>0)=1$;
\item[(iii)] $\ds{\max_{i,j,k\in\{1,\ldots,p\}}\sup_{\theta}\left|A_{n,ii}(\tz)A_{n,jj}(\tz)A_{n,kk}(\tz)\p_{\theta^{i}}\p_{\theta^{j}}\p_{\theta^{k}}\mbbh(\theta)\right|=o_{p}(1)}$.
\end{itemize}
\label{Ass1}
\end{ass}

\medskip

Assumption \ref{Ass1} implicitly sets down the optimal value $\tz$; 
of course, as in the usual $M$-estimation theory (e.g. \cite[Chapter 5]{van98}) 
it is possible to put more specific conditions in terms of the uniform-in-$\theta$ limits of 
suitable scaled quasi-log likelihoods function, but we omit them. 
The quadratic form $\Gam_{0}$ is the asymptotic quasi-Fisher information matrix, which may be random. 
A truly random example is the volatility-parameter estimation of a continuous semimartingale (see Section \ref{Simu2}). 
In particular, Assumption \ref{Ass1} leads to the LAQ approximation of $\log\mbbz_{n}$:
\begin{equation}
\sup_{u\in A}\bigg|\log\mbbz_{n}(u) - \bigg( \D_{n}[u] -\frac{1}{2}\Gam_{0}[u,u] \bigg)\bigg| = o_{p}(1)
\label{hm:pw.LAQ}
\end{equation}
for each compact set $A\subset\mbbr^{p}$. 

\medskip

\begin{ass}
The prior density $\pi_{n}$ satisfies the following:
\begin{itemize}
\item[(i)] $\pi_{n}(\tz)>0$ for all $n$, and $\displaystyle\sup_{n}\sup_{\theta}\pi_{n}(\theta)<\infty$;
\item[(ii)] $\displaystyle\sup_{|u|<M}\big|\pi_{n}(\tz+A_{n}(\tz)u)-\pi_{n}(\tz)\big|\to 0$ as $n\to\infty$ for each $M>0$.
\end{itemize}
\label{Ass2}
\end{ass}

\medskip

\begin{ass}
For any $\epsilon>0$ there exist $M>0$ and $N$ such that
\begin{equation}
\sup_{n\geq N}\pr\bigg(\int_{\mbbu_{n}(\tz)\cap\{|u|\geq M\}}\mbbz_{n}(u)du>\epsilon\bigg)<\epsilon.
\nonumber
\end{equation}
\label{Ass3}
\end{ass}

\medskip

Thanks to Assumption \ref{Ass3}, we can consider a general LAQ models in a unified manner. 
Let us mention some sufficient conditions for the key assumption Assumption \ref{Ass3}. 
To this end we need to introduce further notation. Write
\begin{equation}
\Gam_{n}(\theta)=-A_{n}(\tz)\p_{\theta}^{2}\mbbh_{n}(\theta)A_{n}(\tz),
\nonumber
\end{equation}
and denote by $\lambda_{\min}(M)$ the smallest eigenvalues of a given matrix $M$. 
We write $\ult_{k}=(\theta_{1},\ldots,\theta_{k})$ and $\olt_{k}=(\theta_{k},\ldots,\theta_{K})$, with $\ult_{k,0}$ and $\olt_{k,0}$ in a similar manner. 
Let $u:=(u_{1},\ldots,u_{K})\in\mbbr^{p_{1}}\times\cdots\times\mbbr^{p_{K}}$. The $k$th random field is defined by
\begin{align*}
\mbbz_{n}^{k}(u_{k};\ult_{k-1},\theta_{k,0},\olt_{k+1})
=\exp\left\{\mbbh_{n}(\ult_{k-1},\theta_{k,0}+a_{k,n}(\theta_{0})u_{k},\olt_{k+1})
-\mbbh_{n}(\ult_{k-1},\theta_{k,0},\olt_{k+1})\right\}. 
\end{align*}
The random fields $\mbbz_{n}^{k}$ is designed to focus on the $k$th-graded parameters, when we have more than one rate of convergence, i.e. when $K\ge 2$ 
(We neglect symbols with index $K+1$ like $\olt_{K+1}$ and ones with index $0$ like $\ult_{0}$).

\medskip

\begin{thm}
Let Assumption \ref{Ass1} holds. Then, Assumption \ref{Ass3} follows if at least one of the following conditions holds:
\begin{itemize}
\item[(i)] There exist constants $L>1$ and $C_{L}>0$ such that
\begin{align}
\sup_{n}\pr\bigg(\sup_{(u_{k},\olt_{k+1})\in\{|u_{k}|\geq r \}\times\prk\Theta_{j}}
\mbbz_{n}^{k}(u_{k};\ult_{k-1,0},\tkz,\olt_{k+1})\geq e^{-r}\bigg)\leq\frac{C_{L}}{r^{L}}
\label{pldi.2}
\end{align}
for $r>0$ and $k=1,\ldots,K$;
\item[(ii)] We have
\begin{align}
\limsup_{\del\to 0}\limsup_{n\to\infty}\pr\bigg(\inf_{\theta\in\Theta}
\lambda_{\min}\big(\Gam_{n}(\theta)\big)<\del\bigg)=0.
\label{AssTh3}
\end{align}
\end{itemize}
\label{QlfTh3}
\end{thm}

\medskip

The proof of Theorem \ref{QlfTh3}(i) can be found in \cite[Theorem 6]{Yos11}, together with a general device for how to verify \eqref{pldi.2}; 
what is important in the proof is that, inside the probability, we are bounding the supremum of the random field from below by the quickly decreasing ``$e^{-r}$'' (see also Remark \ref{hm:rem_pldi}).
The proof of \ref{QlfTh3}(ii) is given in Section \ref{hm:sec_proof.sc.tpe}; the condition \eqref{AssTh3} is a sort of {\it global} non-degeneracy condition of the asymptotic information matrix. 
Since we are dealing with the integral-type functional, the non-degeneracy condition may not be {\it local} in $u$.

\begin{rem}{\rm 
As already mentioned in \eqref{hm:pw.LAQ}, Assumption \ref{Ass1} ensures the LAQ structure of the random field $\mbbz_{n}$, 
so that, in the verification of Assumption \ref{Ass3} the uniform-in-$\theta$ asymptotic non-degeneracy of 
the quasi-observed-information matrix $\Gam_{n}(\theta)$ plays a crucial role. 
In the literature, among others: the original \cite{Sch78} considered genuinely Bayesian situation, 
where data was regarded as \textit{non-random} quantities; \cite{CavNea99} proved the key Laplace approximation for the marginal log-likelihood 
under the assumption that the minimum eigenvalue of the possibly random observed information matrix 
is \textit{almost surely} bounded away from zero (and also infinity); \cite{LvLiu14} considered the quasi-likelihood estimation in the generalized linear models 
where the observed information matrix is \textit{non-random}. 
When attempting to directly follow such routes, in general we need to impose almost-sure type condition instead of \eqref{AssTh3}, 
such as the existence of $\del>0$ for which
\begin{align}
\pr\bigg(\limsup_{n\to\infty}\inf_{\theta\in\Theta}\lambda_{\min}\big(\Gam_{n}(\theta)\big)<\del\bigg)=0.
\nn
\end{align}
\label{hm:rem_nrX}
}\qed\end{rem}

\medskip

\begin{rem}{\rm 
For reference, let us mention the tail-probability estimate about the normalized estimator
\begin{equation}
\hat{u}_{n} := A_{n}(\tz)^{-1}(\tes-\tz).
\nonumber
\end{equation}
We can consider the stepwise probability estimates of $\log\mbbz_{n}(u)$ through successive applications of the PLDI result \cite{Yos11}. 
Namely, the following statement holds for a constant $L>0$: {\it if there exists a universal constant $C_{L}>0$ such that
\begin{align}
\sup_{n}\pr\bigg(\sup_{(u_{k},\olt_{k+1})\in\{|u_{k}|\geq r \}\times\prk\Theta_{j}}\mbbz_{n}^{k}(u_{k};\hat{\ult}_{k-1},\tkz,\olt_{k+1})\geq1\bigg)
\leq\frac{C_{L}}{r^{L}}
\label{pldi}
\end{align}
for all $r>0$ and $k=1,\ldots,K$, then $\hat{u}_{n}$ satisfies the estimate
\begin{equation}
\sup_{n}\pr\left(|\hat{u}_{n}| \ge r\right)\leq\frac{C_{L}}{r^{L}},\qquad r>0,
\label{tail.estimate.u-hat}
\end{equation}
which implies the tightness of $(\hat{u}_n)$, hence in particular $\tes\cip \tz$.} 
As in \cite[Proposition 2]{Yos11}, we can derive \eqref{tail.estimate.u-hat} as follows:
we have
\begin{align*}
\pr\left(|\hat{u}_{n}| \ge r\right)\le \sum_{k=1}^{K}
\pr\bigg(\left|a_{k,n}^{-1}(\theta_{0})(\hat{\theta}_{k,n}-\theta_{k,0})\right|\ge\frac{r}{K}\bigg),
\end{align*}
and each $\pr\left(|a_{k,n}^{-1}(\tz)(\hat{\theta}_{k,n}-\tkz)|\ge\frac{r}{K}\right)$ can be bounded by
\begin{align*}
& \pr\bigg(\sup_{\frac{r}{K}\leq|u_{k}|}\big\{\mbbh_{n}\big(\hat{\ult}_{k-1},\tkz+a_{k,n}(\tz)u_{k},\overline{\hat{\theta}}_{k+1}
\big)-\mbbh_{n}\big(\hat{\ult}_{k-1},\tkz,\overline{\hat{\theta}}_{k+1}\big)\big\}\geq0\bigg) \\
&\le\pr\bigg(\sup_{(u_{k},\olt_{k+1})\in\{\frac{r}{K}\leq|u_{k}|\}\times\prk\Theta_{j}}\mbbz_{n}^{k}(u_{k};\hat{\ult}_{k-1},\tkz,\olt_{k+1})\geq1\bigg)
\leq\frac{C_{L}}{r^{L}}K^{L}
\end{align*}
for all $n>0$ and $r>0$.
Sufficient conditions for the PLDI \eqref{pldi.2} and  (\ref{pldi}) to hold can be found in \cite[Theorem 2]{Yos11}. 
The asymptotic mixed normality of $\hat{u}_{n}$ then follows from 
suitable functional weak convergence of $\mbbz_{n}$ on compact sets to a quadratic random fields of suitable form, 
which often follows through a stable convergence in law of the random linear form $\D_{n}=A_{n}(\tz)\p_{\theta}\mbbh_{n}(\tz)$.
\label{hm:rem_pldi}
}\qed\end{rem}

\medskip

Now we can state the stochastic expansion.

\begin{thm}
Suppose that Assumptions \ref{Ass1} to \ref{Ass3} are satisfied and that $\tes\cip \tz$.
\begin{itemize}
\item[(i)] We have the asymptotic expansion
\begin{align*}
\log\bigg(\int_{\Theta}\exp\{\mbbh_{n}(\theta)\}\pi_{n}(\theta)d\theta\bigg)
&=\mbbh_{n}(\tz)+\sumk p_k\log a_{k,n}(\theta_{0})-\frac{1}{2}\log|\Gam_{0}|+\frac{p}{2}\log2\pi \\
&{}\qquad+\frac{1}{2}\|\Gam_{0}^{-\frac{1}{2}}\D_{n}\|^2+\log\pi_{n}(\tz)+o_{p}(1).
\end{align*}

\item[(ii)] If further $\log a_{k,n}(\hat{\theta}_{n})=\log a_{k,n}(\tz)+o_{p}(1)$ and $\log\pi_{n}(\tes)=\log\pi_{n}(\tz)+o_{p}(1)$, then
\begin{align*}
\log\bigg(\int_{\Theta}\exp\{\mbbh_{n}(\theta)\}\pi_{n}(\theta)d\theta\bigg)
&=\mbbh_{n}(\tes)+\sumk p_{k}\log a_{k,n}(\hat{\theta}_{n})+\frac{p}{2}\log2\pi \\
&{}\qquad-\frac{1}{2}\log\big|\Gam_{n}(\tes)\big|+\log\pi_{n}(\tes)+o_{p}(1) \\
&=\mbbh_{n}(\tes)+\frac{p}{2}\log2\pi-\frac{1}{2}\log\big|-\p_{\theta}^{2}\mbbh_{n}(\tes)\big|+\log\pi_{n}(\tes)+o_{p}(1).
\end{align*} 
\end{itemize}
\label{QlfTh2}
\end{thm}

\medskip

It follows from Theorem \ref{QlfTh2}(ii) that the statistics 
\begin{equation}
\qbic^{\sharp}:=-2\mbbh_{n}(\tes)+\log\big|-\p_{\theta}^{2}\mbbh_{n}(\tes)\big|-2\log\pi_{n}(\tes)-p\log2\pi
\label{hm:Fhat_def}
\end{equation}
is a consistent estimator of the marginal quasi-log likelihood function multiplied by ``$-2$''. 
Then, ignoring the $O_{p}(1)$ parts as usual, we define the {\it quasi-Bayesian information criterion (QBIC)} by
\begin{equation}
\qbic=\qbic(\mbX_{n}):=-2\mbbh_{n}(\tes)+\log\big|-\p_{\theta}^{2}\mbbh_{n}(\tes)\big|.
\label{QBIC_def1}
\end{equation}
As long as $\pi_{n}$ is not so dominant and $n$ is moderately large, using $\qbic$ instead of $\qbic^{\sharp}$ would be enough in practice. We compute QBIC for each candidate model, say $\qbic^{(1)},\dots,\qbic^{(M)}$, 
and then define the best model $\mcm_{m_{0}}$ in the sense of approximate Bayesian model description:
\begin{equation}
m_{0}=\argmin_{1\le m\le M}\qbic^{(m)}.
\nonumber
\end{equation}
In view of Assumption \ref{Ass1} we see that
\begin{equation}
\qbic=-2\mbbh_{n}(\tes)+2\sumk p_{k}\log a_{k,n}^{-1}(\hat{\theta}_{n})+O_{p}(1).
\label{QBIC_def2}
\end{equation}
Since the second term in the right-hand side diverges in probability, we could more simply define QBIC to be the the sum of the first two terms in the right-hand side of \eqref{QBIC_def2}. 
We may thus define {\it Schwarz's BIC} in our context by
\begin{equation}
\bic=-2\mbbh_{n}(\tes)+2\sumk p_{k}\log a_{k,n}^{-1}(\hat{\theta}_{n}).
\label{BIC_def}
\end{equation}
Note that in the classical case of single $\sqrt{n}$-scaling \eqref{QBIC_def2} reduces to the familiar form
\begin{equation}
\bic=-2\mbbh_{n}(\tes)+p \log n.
\label{hm:c.bic}
\end{equation}
The statistics $\qbic$ thus provides us with a far-reaching extension of derivation machinery of the classical BIC. 

Although the original definition (\ref{QBIC_def1}) has higher computational load than (\ref{QBIC_def2}), 
it enables us to incorporate a model-complexity bias correction taking the volume of observed information into account. 
In particular, to reflect data information for dependent-data models, (\ref{QBIC_def1}) would be more suitable than (\ref{QBIC_def2}) 
whose bias correction is only based on the rate of convergence. 

\medskip

\begin{rem}{\rm 
Making use of the observed information matrix \eqref{QBIC_def1} for regularization has been already mentioned in the literature; 
for example, \cite{Boz87}, \cite{Kas82}, and \cite{Scl87} contain such statistics for some variants of the AIC statistics. 
Further, it is worth mentioning that using the observed-information is a right way for some non-stationary models \cite{Kim98}.
}\end{rem}

\medskip

\begin{rem}{\rm 
At the beginning the prior model-selection probabilities $\pi_{1,n},\dots,\pi_{M,n}$ are to be set in a subjective manner. 
As usual, using the QBIC of the candidate models we may estimate the posterior model-selection probabilities 
in the data-driven manner through the quantities
\begin{equation}
\hat{\pi}_{m,n}(\mbX_{n}):=\frac{\mfp_{m}\exp\{-\qbic^{(m)}(\mbX_{n})/2\}}
{\sum_{l=1}^{M}\mfp_{l}\exp\{-\qbic^{(l)}(\mbX_{n})/2\}},\qquad m=1,\dots,M,
\nonumber
\end{equation}
or those with QBIC replaced by BIC.
}\qed\end{rem}

\medskip

\begin{rem}[Variants of QBIC]{\rm 
In practice we may conveniently consider several variants of the QBIC \eqref{QBIC_def1}. 
When $\Gam_{0}$ takes the form $\Gam_{0}=\diag(\Gam_{10},\dots,\Gam_{K0})$ with each $\Gam_{k0}\in\mbbr^{p_{k}}\otimes\mbbr^{p_{k}}$ being a.s. positive-definite, 
we may slightly simplify the form of the QBIC as follows. 
Since under Assumption \ref{Ass1} we have
\begin{equation}
-a_{k,n}(\tz)a_{l,n}(\tz)\p_{\theta_{k}}\p_{\theta_{l}}\mbbh_{n}(\tes)=o_{p}(1),\qquad k\ne l,
\nonumber
\end{equation}
and also taking logarithmic determinant of a positive-definite matrix is continuous, the basic asymptotic expansion becomes
\begin{align}
\log\bigg(\int_{\Theta}\exp\{\mbbh_{n}(\theta)\}\pi_{n}(\theta)d\theta\bigg)
&=\mbbh_{n}(\hat{\theta}_{n})-\frac{1}{2}\sum_{k=1}^{K}\log\big|-\p_{\theta_{k}}^{2}\mbbh_{n}(\hat{\theta}_{n})\big|+ O_{p}(1),
\nn
\end{align}
giving rise to the QBIC of the form
\begin{equation}
-2\mbbh_{n}(\hat{\theta}_{n})+\sum_{k=1}^{K}\log\big|-\p_{\theta_{k}}^{2}\mbbh_{n}(\hat{\theta}_{n})\big|.
\label{mQBIC_def4}
\end{equation}
This is the case if $A_{n}(\tz)(\tes-\tz) \cil MN(0,\Sig_{0})$, the symbol ``$MN$'' referring to ``mixed-normal'', 
we may have $\Sig_{0}=\diag(\Sig_{10},\dots,\Sig_{K0})$ with each $\Sig_{k0}\in\mbbr^{p_{k}}\otimes\mbbr^{p_{k}}$ being a.s. positive-definite. 
In Section \ref{GQMLE2}, we will deal with an example where the estimator is asymptotically normally distributed at two different rates, with the asymptotic covariance matrix of $\hat{u}_{n}$ being block diagonal. 

We may also consider finite-sample manipulations of QBIC without breaking its asymptotic behavior. 
For example, the problem caused by $|-\p_{\theta}^{2}\mbbh_{n}(\tes)|\le 0$ can be avoided by using
\begin{align}
&-2\mbbh_{n}(\tes) + I\left\{ |-\p_{\theta}^{2}\mbbh_{n}(\tes)|>0 \right\}\log\big|-\p_{\theta}^{2}\mbbh_{n}(\tes)\big| \nn\\
&{}\qquad+ I\left\{ |-\p_{\theta}^{2}\mbbh_{n}(\tes)|\le 0 \right\}\sumk p_{k}\log \left(a_{k,n}^{-2}(\hat{\theta}_{k,n})\right)
\nn
\end{align}
instead of \eqref{QBIC_def1}; obviously, the difference between this quantity and $\qbic$ is of $o_{p}(1)$. 
Further, we may use any $\hat{\Gam}_{n}$ such that $\hat{\Gam}_{n}\cip\Gam_{0}$:
\begin{align}
-2\mbbh_{n}(\tes) -2\log|A_{n}(\tes)| + \log|\hat{\Gam}_{n}|,
\nn
\end{align}
which would be convenient if $\hat{\Gam}_{n}$ is more likely to be stable than $\Gam_{n}(\tes)=-A_{n}(\tes)\p_{\theta}^{2}\mbbh_{n}(\tes)A_{n}(\tes)$; 
for example, if we beforehand know the specific form of $\Gam_{0}=\Gam_{0}(\theta)$, 
then it would be (numerically) more stable to use $\Gam_{0}(\tes)$ instead of $\Gam_{n}(\tes)$.
\label{hm:rem_variants}
}\qed\end{rem}

\subsection{Convergence of the expected values}\label{hm:qbic.mc}

From the frequentist point of view where $\mbX_{n}$ is regarded as a random element, 
it may be desirable to verify the convergence of {\it expected} marginal quasi-log likelihood, 
which follows from the asymptotic uniform integrability of the sequence 
\begin{equation}
\bigg\{ \bigg|-2\log\bigg(\int_{\Theta}\exp\{\mbbh_{n}(\theta)\}\pi_{n}(\theta)d\theta\bigg) - \qbic^{\sharp}\bigg| \bigg\}_{n}.
\nonumber
\end{equation}
In particular, $\qbic^{\sharp}$ will be ensured to be an asymptotically unbiased estimator of the expected logarithmic Bayes factor; see Section \ref{hm:sec_ebf}, in particular \eqref{hm:lfb.L1conv}.

Let us recall the notation $\D_{n}=A_{n}(\tz)\p_{\theta}\mbbh_{n}(\tz)$ and $\Gam_{n}(\theta)=-A_{n}(\tz)\p_{\theta}^{2}\mbbh_{n}(\theta)A_{n}(\tz)$. 
First we strengthen Assumptions \ref{Ass1} and \ref{Ass2}.

\begin{ass}
The random function $\mbbh_{n}$ is of class $\mcc^{3}(\Theta)$ a.s. and for every $r>0$
\begin{equation}
\sup_{n}\E\bigg( |\D_{n}|^{r} + \sup_{\theta}|\Gam_{n}(\theta)|^{r}
+ \sum_{i=1}^{p}\sup_{\theta}\left|A_{n}(\tz)\p_{\theta^{i}}\p_{\theta}^{2}\mbbh_{n}(\theta)A_{n}(\tz)\right|^{r}
\bigg)<\infty.
\nonumber
\end{equation}
\label{hm:A_mc1}
\end{ass}

\begin{ass}
In addition to Assumption \ref{Ass2}, we have $\ds{0<\inf_{n,\theta}\pi_{n}(\theta)\le\sup_{n,\theta}\pi_{n}(\theta)<\infty}$.
\label{hm:A_mc2}
\end{ass}

Next we strengthen Assumptions \ref{Ass3} by the following.

\begin{ass}
The exists an a.s. positive definite random matrix $\Gam_{0}$ such that $\Gam_{n}(\tz)\cip\Gam_{0}$, and for some $q>3p$ we have
\begin{equation}
\limsup_{n}\E\bigg( \sup_{\theta}\lam_{\min}^{-q}\left(\Gam_{n}(\theta)\right) \bigg)<\infty.
\nonumber
\end{equation}
\label{hm:A_mc3}
\end{ass}

The moment bounds in Assumption \ref{hm:A_mc3} was studied in \cite{ChaHuaIng13} and \cite{ChaIng11} for some time series models, 
with a view toward prediction. 
The integrability in Assumption \ref{hm:A_mc3} is related to the key index $\chi_{0}$ of \cite{UchYos13} in case of volatility estimation of continuous It\^o process.

Under Assumptions \ref{hm:A_mc1} and \ref{hm:A_mc3} we have 
$\lam_{\min}^{-q}(\Gam_{n}(\tz))\cip\lam_{\min}^{-q}(\Gam_{0})$ 
by the continuous mapping theorem, and also $\lam_{\min}^{-1}(\Gam_{0})\in L^{q}(\pr)$ 
as well as $\Gam_{0}\in\bigcap_{r>0}L^{r}(\pr)$.


Finally, we impose the boundedness of moments of the normalized estimator; see Remark \ref{hm:rem_pldi}.
\begin{ass}
$\ds{\sup_{n}\E(|\hat{u}_{n}|^{r})<\infty}$ for some $r>3$.
\label{hm:A_mc5}
\end{ass}

We can now state the $L^{1}(\pr)$-converge result.

\begin{thm}
If Assumption \ref{Ass3} and Assumptions \ref{hm:A_mc1} to \ref{hm:A_mc5} hold and $\log a_{k,n}(\hat{\theta}_{n})=\log a_{k,n}(\tz)+o_{p}(1)$ for $k=1,\dots, K$, then we have
\begin{equation}
\lim_{n\to\infty}\E\bigg\{\bigg|-2\log\bigg(\int_{\Theta}\exp\{\mbbh_{n}(\theta)\}\pi_{n}(\theta)d\theta\bigg) - \qbic^{\sharp}\bigg| \bigg\} =0.
\nonumber
\end{equation}
In particular, $\qbic^{\sharp}$ is an asymptotically unbiased estimator of the logarithm of the quasi-marginal likelihood. 
\label{hm:thm_mc}
\end{thm}

\section{Gaussian quasi-likelihood} \label{hm:sec_GQNLE}

This section is devoted to the Gaussian quasi-likelihood.

\subsection{General framework} \label{GQMLE1}
A general setting for the Gaussian quasi-likelihood estimation is described as follows. Let $\mbX_{n}=(X_{n,j})_{j=0}^{n}=(X_{n,0},\ldots,X_{n,n})$ be an array of random variables, where $X_{n,j}\in\mbbr$ for brevity. Let $\mcf_{n,j}:=\sig(X_{n,j};j\leq n)$ denote the $\sig$-field representing the data information at stage $j$ when the total number of data is $n$. The Gaussian quasi-likelihood (in the univariate case) is constructed as if the conditional distribution of $X_{n,j}$ given past information $\mcf_{n,j-1}$ is Gaussian, say
\begin{align*}
\mcl(X_{n,j}|X_{n,0},\ldots,X_{n,j-1})\thickapprox N(\mu_{n,j-1}(\theta),\sig_{n,j-1}(\theta)),
\end{align*}
where $\mu_{n,j-1}$ and $\sig_{n,j-1}$ are $\mcf_{n,j-1}$-measurable (predictable) random function on $\Theta$; most often,
\begin{align*}
\mu_{n,j-1}(\theta)=\E(X_{n,j}|\mcf_{n,j-1}),\qquad\sig_{n,j-1}(\theta)=\var(X_{n,j}|\mcf_{n,j-1}),
\end{align*}
where the conditional expectation and variance are taken under the image measure of $\mbX_{n}$ associated with the parameter value $\theta$. 
In what follows, we will suppress the subscript $``n"$. 

Because the quasi-likelihood is given by
\begin{align*}
\theta&\mapsto\sumj\log\frac{1}{\sqrt[]{2\pi\sig_{j-1}^{2}(\theta)}}\exp\left\{-\frac{1}{2\sig_{j-1}^{2}(\theta)}\big(X_{j}-\mu_{j-1}(\theta)\big)^{2}\right\} \\
&=(\mathrm{const.})+\left[-\frac{1}{2}\sumj\left\{\log\sig_{j-1}^{2}(\theta)+\frac{\big(X_{j}-\mu_{j-1}(\theta)\big)^{2}}{\sig_{j-1}^2(\theta)}\right\}\right],
\end{align*}
we may define the Gaussian quasi-likelihood function by
\begin{align*}
\mbbh_{n}(\theta)=-\frac{1}{2}\sumj\left\{\log\sig_{j-1}^{2}(\theta)+\frac{\big(X_{j}-\mu_{j-1}(\theta)\big)^{2}}{\sig_{j-1}^{2}(\theta)}\right\}.
\end{align*}
Then, supposing that $\mbbh_{n}$ and its partial derivatives can be continuously extended to the boundary $\p\Theta$, we define the Gaussian QMLE (GQMLE) by any maximizer of $\mbbh_{n}$ over $\bar{\Theta}$.

The Gaussian quasi-likelihood is designed to fit not full joint distribution but only conditional-mean and conditional-covariance structures. 
The simplest case is the location-parameter estimation by the sample mean in the i.i.d.-data setting, 
where $\sig^{2}_{j-1}(\theta)\equiv 1$ (set for brevity) and $\mu_{j-1}(\theta)=\theta$, namely the least-squares estimation without ``full'' specification of the underlying population distribution.
Although the GQMLE is not (possibly far from being) asymptotically efficient when the model is misspecified, the GQMLE quite often exhibits asymptotic (mixed-)normality under appropriate conditions even if the conditional distribution is deviating from being normal.


\subsection{Ergodic diffusion process}\label{GQMLE2} 

Let $\mbX_{n}=(X_{t_{j}})_{j=0}^{n}$ with $t_{j}=jh_{n}$, where $h_{n}$ is the discretization step and $nh_{n}=T_{n}$ and $X_{t}$ is a solution to the $d$-dimensional diffusion process defined by the stochastic differential equation
\begin{align*}
dX_{t}=a(X_{t})dt+b(X_{t})dw_{t},\quad t\in[0,T_{n}],~X_{0}=x_{0}.
\end{align*}
Here $a$ is an $\mbbr^{d}$-valued function defined on $\mbbr^{d}$, $b$ is an $\mbbr^{d}\otimes\mbbr^{d}$-valued function defined on $\mbbr^{d}$, $w_{t}$ is an $d$-dimensional standard Wiener process, and $x_{0}$ is a deterministic initial value.
We assume that $h_{n}\to0$, $T_{n}=nh_{n}\to\infty$, $nh_{n}^{2}\to0$ as $n\to0$ and that for some positive constant $\epsilon_{0}$, $nh_{n}\geq n^{\epsilon_{0}}$ for every large $n$. 
Let us consider the following stochastic differential equation as statistical model $\mcm_{m_{1},m_{2}}$:
\begin{align}
dX_{t}=a_{m_{2}}(X_{t},\theta_{m_{2}})dt+b_{m_{1}}(X_{t},\theta_{m_{1}})dw_{t},\quad t\in[0,T_{n}],\;X_{0}=x_{0}, \label{mod.erg}
\end{align}
where $a_{m_{2}}$ is an $\mbbr^{d}$-valued function defined on $\mbbr^{d}\times\Theta_{m_{2}}$, $b_{m_{1}}$ is an $\mbbr^{d}\otimes\mbbr^{d}$-valued function defined on $\mbbr^{d}\times\Theta_{m_{1}}$ and $(m_{1},m_{2})\in\{1,\ldots,M_{1}\}\times\{1,\ldots,M_{2}\}$; namely, we consider $M_{1}\times M_{2}$ models in total. 
In each model $\mcm_{m_{1},m_{2}}$, the coefficients $b_{m_{1}}$ and $a_{m_{2}}$ are assumed to be known up to the finite-dimensional parameter 
$\theta_{m_{1},m_{2}}:=(\theta_{m_{1}},\theta_{m_{2}})\in\Theta_{m_{1}}\times\Theta_{m_{2}}\subset\mbbr^{p_{m_{1}}}\times\mbbr^{p_{m_{2}}}$.
We focus on the case of correctly specified parametric coefficients: 
we assume that for each $m$ there exists the true value $(\theta_{m_{1},0},\theta_{m_{2},0})$ for which $b_{m_{1}}(\cdot,\theta_{m_{1},0})=b(\cdot)$ and $a_{m_{2}}(\cdot,\theta_{m_{2},0})=a(\cdot)$.

Below, we omit the model index ``$m_{1}$" and ``$m_{2}$" from the notation.
That is, the stochastic differential equation (\ref{mod.erg}) is expressed by
\begin{align*}
dX_{t}=a(X_{t},\theta_{2})dt+b(X_{t},\theta_{1})dw_{t},\;t\in[0,T_{n}],\;X_{0}=x_{0}.
\end{align*}
Let $B(x,\theta_{1}):=b(x,\theta_{1})b^{\prime}(x,\theta_{1})$ and $\D_{j}X:=X_{t_{j}}-X_{t_{j-1}}$.  
We obtain the quasi-likelihood function
\begin{align*}
\prod_{j=1}^{n}(2\pi h_{n})^{-\frac{d}{2}}\big|B(X_{t_{j-1}},\theta_{1})\big|^{-\frac{1}{2}}\exp\left\{-\frac{1}{2h_{n}}B(X_{t_{j-1}},\theta_{1})^{-1}\left[\big(\D_{j}X-h_{n}a(X_{t_{j-1}},\theta_{2})\big)^{\otimes2}\right]\right\},
\end{align*}
where $x^{\otimes 2}:=xx'$.
Then, up to an additive constant common to all the candidate models, the quasi-log likelihood function is given by
\begin{align}
\mbbh_{n}(\theta)&=-\frac{1}{2}\sumj\left\{\log\big|B(X_{t_{j-1}},\theta_{1})\big|+\frac{1}{h_{n}}B(X_{t_{j-1}},\theta_{1})^{-1}\left[\big(\D_{j}X-h_{n}a(X_{t_{j-1}},\theta_{2})\big)^{\otimes2}\right]\right\},\label{ErLf}
\end{align}
Let $A_{n}(\tz)=\diag\big(\frac{1}{\sqrt{n}}I_{p_{1}},\frac{1}{\sqrt{nh_{n}}}I_{p_{2}}\big)$ be the rate matrix.

We assume the following conditions \cite[Section 6]{Yos11}:

\begin{ass}
\begin{itemize}
\item[(i)] For some constant $C$,
\begin{align*}
& \sup_{\theta_{2}\in\Theta_{2}}|\p_{\theta_{2}}^{i}a(x,\theta_{2})|\leq C(1+|x|)^{C}\quad(0\leq i\leq4), \nn\\
& \sup_{\theta_{1}\in\Theta_{1}}|\p_{x}^{j}\p_{\theta_{1}}^{i}b(x,\theta_{1})|\leq C(1+|x|)^{C}\quad(0\leq i\leq4,0\leq j\leq2).
\end{align*}

\item[(ii)] $\displaystyle \inf_{|u|=1}\inf_{(x,\theta_{1})}B(x,\theta_{1})[u,u]>0$.

\item[(iii)] There exists a constant $C$ such that for every $x_{1},x_{2}\in\mbbr^{p}$,
\begin{align*}
\sup_{\theta_{2}\in\Theta_{2}}|a(x_{1},\theta_{2})-a(x_{2},\theta_{2})|+\sup_{\theta_{1}\in\Theta_{1}}|b(x_{1},\theta_{1})-b(x_{2},\theta_{1})|
\leq C|x_{1}-x_{2}|
\end{align*}

\item[(iv)] $X_{0}\in\bigcap_{p>0}L^{p}(\pr)$.

\end{itemize}
\label{Ass4}
\end{ass}

\begin{ass}
For some constant $a>0$,
\begin{align*}
\sup_{t\in\mbbr_{+}}\sup_{\substack{A\in\sigma[X_{r};r\leq t] \\ B\in\sigma[X_{r};r\geq t+h]}}\big|\pr(A\cap B)-\pr(A)\pr(B)\big|\leq a^{-1}e^{-ah}\quad(h>0).
\end{align*}
\label{Ass5}
\end{ass}

Under Assumption \ref{Ass5} ensures the ergodicity: 
there exists a unique invariant probability measure $\nu=\nu_{\tz}$ of $X_{t}$ such that
\begin{align*}
\frac{1}{T}\int_{0}^{T}g(X_{t})dt\cip\int_{\mbbr^{d}}g(x)\nu(dx)\quad(T\to\infty)
\end{align*}
for any bounded measurable function $g$.

\begin{ass}
There exists a positive constant $\chi>0$ such that $\mbby_{1,0}(\theta_{2})\leq-\chi|\theta_{1}-\theta_{1,0}|^{2}$ for all $\theta_{1}\in\Theta_{1}$, where
\begin{align*}
\mbby_{1,0}(\theta_{1})=-\frac{1}{2}\int_{\mbbr^{d}}\left\{\tr\big(B(x,\theta_{1})^{-1}B(x,\theta_{1,0})-I_{p}\big)+\log\frac{|B(x,\theta_{1})|}{|B(x,\theta_{1,0})|}\right\}\nu(dx).
\end{align*}
\label{Ass6}
\end{ass}

\begin{ass}
There exists a positive constant $\chi^{\prime}>0$ such that $\mbby_{2,0}(\theta_{2})\leq-\chi^{\prime}|\theta_{2}-\theta_{2,0}|^{2}$ for all $\theta_{2}\in\Theta_{2}$, where
\begin{align*}
\mbby_{2,0}(\theta_{2})=-\frac{1}{2}\int_{\mbbr^{d}}B(x,\theta_{1,0})^{-1}\left[\big(a(x,\theta_{2})-a(x,\theta_{2,0})\big)^{\otimes2}\right]\nu(dx).
\end{align*}
\label{Ass7}
\end{ass}

The partial derivatives of $\mbbh_{n}$ are given as follows: 
for $u_{1}\in{\mbbr^{m_{1}}}$ and $u_{2}\in\mbbr^{m_{2}}$,
\begin{align*}
\p_{\theta_{1}}^{2}\mbbh_{n}(\theta_{1},\theta_{2})[u_{1}^{\otimes2}]
&=-\frac{1}{2}\sumj\bigg\{\p_{\theta_{1}}^{2}\log\frac{|B(X_{t_{j-1}},\theta_{1})|}{|B(X_{t_{j-1}},\theta_{1,0})|}
\nn\\
&{}\qquad
+\frac{1}{h_{n}}\p_{\theta_{1}}^{2}B(X_{t_{j-1}},\theta_{1})^{-1}\left[u_{1}^{\otimes2},\big(\D_{j}X-h_{n}a(X_{t_{j-1}},\theta_{2})\big)^{\otimes2}\right]\bigg\}, \\
\p_{\theta_{2}}^{2}\mbbh_{n}(\theta_{1},\theta_{2})\left[u_{2}^{\otimes2}\right]
&=-\sumj B(X_{t_{j-1}},\theta_{1})^{-1}\bigg\{\big[\p_{\theta_{2}}a(X_{t_{j-1}},\theta_{2})[u_{2}],\p_{\theta_{2}}h_{n}a(X_{t_{j-1}},\theta_{2})[u_{2}]\big] \\
&\qquad\qquad-\left[\p_{\theta_{2}}^{2}a(X_{t_{j-1}},\theta_{2})\left[u_{2}^{\otimes2}\right],\D_{j}X-h_{n}a(X_{t_{j-1}},\theta_{2})\right]\bigg\},
\nn\\
\p_{\theta_{1}}\p_{\theta_{2}}\mbbh_{n}(\theta_{1},\theta_{2})[u_{1},u_{2}]
&=\sumj\p_{\theta_{1}}B(X_{t_{j-1}},\theta_{1})\left[u_{1},\p_{\theta_{2}}a(X_{t_{j-1}},\theta_{2})u_{2},\D_{j}X-h_{n}a(X_{t_{j-1}},\theta_{2})\right].
\end{align*}
Then, we obtain the corresponding QBIC as in the following theorem. The proof is given in Section \ref{hm:sec_Proofs}.

\begin{thm}
Suppose that Assumptions \ref{Ass4} to \ref{Ass7} are satisfied. Then, the assumptions in Theorem \ref{QlfTh2} are satisfied and the corresponding QBIC is given by
\begin{align*}
\qbic&=\sumj\left\{\log|B(X_{t_{j-1}},\hat{\theta}_{1,n})|+\frac{1}{h_{n}}B(X_{t_{j-1}},\hat{\theta}_{1,n})^{-1}\left[\big(\D_{j}X-h_{n}a(X_{t_{j-1}},\hat{\theta}_{2,n})\big)^{\otimes2}\right]\right\} \\
&\qquad\qquad+\log\left|-\left({\renewcommand\arraystretch{1}\begin{array}{ll}\p_{\theta_{1}}^{2}\mbbh_{n}(\tes)&\p_{\theta_{1}}\p_{\theta_{2}}\mbbh_{n}(\tes)\\ \p_{\theta_{1}}\p_{\theta_{2}}\mbbh_{n}(\tes)&\p_{\theta_{2}}^{2}\mbbh_{n}(\tes)\end{array}}\right)\right|.
\end{align*}
\label{ErTh1}
\end{thm}

In the present case of ergodic diffusion process, the convergence in probability
\begin{align}
\frac{1}{\sqrt{n^{2}h_{n}}}\p_{\theta_{1}}\p_{\theta_{2}}\mbbh_{n}(\tes)\cip0\quad(n\to\infty)
\label{er_cip}
\end{align}
is satisfied, so that
\begin{align*}
\log\left|-A_{n}(\tes)\p_{\theta}^{2}\mbbh_{n}(\tes)A_{n}(\tes)\right|
&=\log\left|{\renewcommand\arraystretch{1.5}\begin{array}{cc}-\frac{1}{n}\p_{\theta_{1}}^{2}\mbbh_{n}(\tes)&-\frac{1}{\sqrt{n^{2}h_{n}}}\p_{\theta_{1}}\p_{\theta_{2}}\mbbh_{n}(\tes)\\ -\frac{1}{\sqrt{n^{2}h_{n}}}\p_{\theta_{1}}\p_{\theta_{2}}\mbbh_{n}(\tes)^{\prime}&-\frac{1}{nh_{n}}\p_{\theta_{2}}^{2}\mbbh_{n}(\tes)\end{array}}\right| \\
&=\log\left|{\renewcommand\arraystretch{1.5}\begin{array}{cc}-\frac{1}{n}\p_{\theta_{1}}^{2}\mbbh_{n}(\tes)&0\\ 0&-\frac{1}{nh_{n}}\p_{\theta_{2}}^{2}\mbbh_{n}(\tes)\end{array}}\right|+o_{p}(1) \\
&=\log\big|A_{n}(\tes)\diag\big(-\p_{\theta_{1}}^{2}\mbbh_{n}(\tes),-\p_{\theta_{2}}^{2}\mbbh_{n}(\tes)\big)A_{n}(\tes)\big|+o_{p}(1).
\end{align*}
In the asymptotic framework, statistics $\hat{S}_{n}$ such that $\hat{S}_{n}$ is easier to compute and that $\hat{S}_{n}=\qbic+O_{p}(1)$ may be used as a variant of $\qbic$; recall \eqref{QBIC_def1} and \eqref{QBIC_def2}, and also Remark \ref{hm:rem_variants}.

\begin{thm}
Assume that Assumptions \ref{Ass4}-\ref{Ass7} hold, then the difference between the statistics
\begin{align*}
& \sumj\left\{\log|B(X_{t_{j-1}},\hat{\theta}_{1,n})|+\frac{1}{h_{n}}B(X_{t_{j-1}},\hat{\theta}_{1,n})^{-1}\left[\big(\D_{j}X-h_{n}a(X_{t_{j-1}},\hat{\theta}_{2,n})\big)^{\otimes2}\right]\right\} \\
&\qquad\qquad+\log|-\p_{\theta_{1}}^{2}\mbbh_{n}(\tes)|+\log|-\p_{\theta_{2}}^{2}\mbbh_{n}(\tes)|
\end{align*}
and the $\qbic$ given in Theorem \ref{ErTh1} is $o_{p}(1)$.
\label{ErTh2}
\end{thm}

The BIC corresponding to \eqref{BIC_def} takes the from
\begin{align*}
& \sumj\left\{\log|B(X_{t_{j-1}},\hat{\theta}_{1,n})|+\frac{1}{h_{n}}B(X_{t_{j-1}},\hat{\theta}_{1,n})^{-1}\left[\big(\D_{j}X-h_{n}a(X_{t_{j-1}},\hat{\theta}_{2,n})\big)^{\otimes2}\right]\right\} \\
&\qquad\qquad+p\log n + p_{2}\log h_{n},
\end{align*}
clarifying that the high frequency of data indeed has the significant impact through the diverging term ``$p_{2}\log h_{n}$'';  
one might formally set the formal-BIC term to be ``$p\log n$'', but it is incorrect for the present high-frequency data.

\begin{rem}{\rm
It follows from Kessler \cite{Kes97} that we have 
\begin{align*}
\big(\sqrt{n}(\hat{\theta}_{1,n}-\theta_{1,0}),\sqrt{nh_{n}}(\hat{\theta}_{2,n}-\theta_{2,0})\big)\cil N_{p}\big(0,\diag\big(\Gam_{1,0}(\theta_{1,0})^{-1},\Gam_{2,0}(\theta_{1,0},\theta_{2,0})^{-1}\big)\big),
\end{align*}
where
\begin{align*}
&\Gam_{1,0}(\theta_{1,0})[u_{1}^{\otimes2}]=\frac{1}{2}\int\tr\big\{B(x,\theta_{1,0})^{-1}\big(\p_{\theta_{1}}B(x,\theta_{1,0})\big)B(x,\theta_{1,0})^{-1}\big(\p_{\theta_{1}}B(x,\theta_{1,0})\big)[u_{1}^{\otimes2}]\big\}\nu(dx), \\
&\Gam_{2,0}(\theta_{1,0},\theta_{2,0})[u_{2}^{\otimes2}]=\int B(x,\theta_{1,0})^{-1}\big[\p_{\theta_{2}}a(x,\theta_{2,0})[u_{2}],\p_{\theta_{2}}a(x,\theta_{2,0})[u_{2}]\big]\nu(dx)
\end{align*}
for $u_{1}\in\mbbr^{m_{1}}$, $u_{2}\in\mbbr^{m_{2}}$.
We know from Gobet \cite{Gob02} that this GQMLE is asymptotically efficient in the sense of Haj\'{e}k-Le Cam.
\label{ErRem1}
}\qed\end{rem}

\begin{rem}[Separately convex example]{\rm 
Consider the following class of univariate stochastic differential equations:
\begin{align*}
dX_{t}=\bigg(\sum_{k=1}^{p_{2}}\theta_{2,k}a_{k}(X_{t})\bigg)dt
+\exp\bigg(\frac{1}{2}\sum_{\ell=1}^{p_{1}}\theta_{1,\ell}b_{\ell}(X_{t})\bigg)dw_{t},
\end{align*}
where $a_{k}$ and $b_{l}$ are known functions (basis functions). 
Write $\theta_{1}=(\theta_{1,1},\ldots,\theta_{1,p_{1}})^{\prime}$, $\theta_{2}=(\theta_{2,1},\ldots,\theta_{2,p_{2}})^{\prime}$, $a(x)=\big(a_{1}(x),\ldots,a_{p_{2}}(x)\big)^{\prime}$, $b(x)=\big(b_{1}(x),\ldots,b_{p_{1}}(x)\big)^{\prime}$.
Then the quasi-likelihood function is given by
\begin{align*}
\mbbh_{n}(\theta_{1},\theta_{2})&=-\frac{1}{2}\sumj\left\{\theta_{1}^{\prime}b(X_{t_{j-1}})
+\frac{1}{h_{n}}\big(\D_{j}X-h_{n}a(X_{t_{j-1}})'\theta_{2}\big)^{2}\exp\{-b(X_{t_{j-1}})'\theta_{1}\}\right\}.
\end{align*}
The corresponding joint QBIC of Theorem \ref{ErTh2} is given by
\begin{align*}
\qbic&=-2\mbbh_{n}(\hat{\theta}_{1,n},\hat{\theta}_{2,n})
+\log\bigg|h_{n}\sumj \exp\big\{-\hat{\theta}_{1,n}^{\prime}b(X_{t_{j-1}})\big\}a^{\otimes 2}(X_{t_{j-1}})\bigg| \\
&\qquad+\log\bigg|\frac{1}{2}\sumj\frac{1}{h_{n}}\exp\big\{-\hat{\theta}_{1,n}^{\prime}b(X_{t_{j-1}})\big\}\big(\D_{j}X-h_{n}\hat{\theta}_{2,n}^{\prime}a(X_{t_{j-1}})\big)^{2}b^{\otimes 2}(X_{t_{j-1}})\bigg|.
\end{align*}

Several adaptive-estimation methodologies for general parametric ergodic diffusions have been developed in the literature; see \cite{UchYos12} and \cite{KamUch15} as well as the references therein.
We here remark that, under mild conditions on the functions $a$ and $b$, the optimization may be made even simpler and more efficient by using an adaptive estimation strategy. 
This is because of the convexity of each of the random functions to be optimized: specifically, we first get an estimate $\hat{\theta}_{1,n}$ of $\theta_{1}$ by the convex random function
\begin{equation}
\mbbh_{1,n}(\theta_{1}):=-\frac{1}{2}\sumj\bigg(\theta_{1}^{\prime}b(X_{t_{j-1}})+\frac{1}{h_{n}}\big(\D_{j}X\big)^{2}\exp\{-\theta_{1}'b(X_{t_{j-1}})\}\bigg)
\nonumber
\end{equation} 
(regarding $a(x)\equiv 0$ in the original $\mbbh_{n}(\theta_{1},\theta_{2})$). 
Second, we get an estimate $\hat{\theta}_{2,n}$ by the explicit maxima of the convex random function
\begin{equation}
\mbbh_{n}(\hat{\theta}_{1,n},\theta_{2}) = -\frac{1}{2}\sumj \frac{1}{h_{n}}\big(\D_{j}X-h_{n}\theta_{2}'a(X_{t_{j-1}})\big)^{2}\exp\{-\hat{\theta}_{1,n}'b(X_{t_{j-1}})\}.
\nonumber
\end{equation}
This framework naturally provides us with an adaptive model-selection procedure, see Section \ref{hm:sec_msc.ergo.diff} for details.
\label{hm:ex.separate.convex.diffusion}
}\qed\end{rem}

\subsection{Volatility-parameter estimation for continuous semimartingale} \label{GQMLE3}

In this section, we deal with the stochastic regression model
\begin{align*}
dY_{t}=b_{t}dt+\sigma(X_{t},\theta)dw_{t},\quad t\in[0,T],
\end{align*}
where $w$ is an $r$-dimensional standard Wiener process, $b$ and $X$ are progressively measurable processes with values in $\mbbr^{m}$ and $\mbbr^{d}$, respectively, $\sigma$ is an $\mbbr^{m}\otimes\mbbr^{r}$-valued function defined on $\mbbr^{d}\times\Theta$ with $\Theta\in\mbbr^{p}$. Data set consists of discrete observations ${\bf X}_{n}=(X_{t_{j}},Y_{t_{j}})_{j=0}^{n}$ with $t_{j}=jh_{n}$, where $h_{n}=T/n$ with $T$ fixed. The process $b$ is completely unobservable and unknown. All processes are defined on a filtered probability space $\mcb:=(\Omega,\mcf,(\mcf_{t})_{t\le T},P)$.

Let $S(x,\theta):=\sigma(x,\theta)\sigma(x,\theta)^{\prime}$ and $\D_{j}Y:=Y_{j}-Y_{j-1}$. Then the quasi-likelihood function becomes
\begin{align*}
\mbbh_{n}(\theta)=-\frac{1}{2}\sumj\left\{\log\big|S(X_{t_{j-1}},\theta)\big|+\frac{1}{h_{n}}S(X_{t_{j-1}},\theta)^{-1}\left[(\D_{j}Y)^{\otimes2}\right]\right\}.
\end{align*}
The asymptotic distribution of $A_{n}(\tz)^{-1}(\tes-\tz)=\sqrt{n}(\tes-\tz)$ is mixed normal, i.e.
\begin{align*}
\sqrt{n}(\tes-\tz)\cil \Sigma_{\theta}^{-1/2}Z,
\end{align*}
where $\Sigma_{\theta}$ is a symmetric $p\times p$-matrix which is a.s. positive-definite, 
and $Z$ is a $p$-variate standard-normal random variable which is defined on an extension of $\mcb$ and is independent of $\mcb$, see \cite{GenJac93} and \cite{UchYos13}.

The QBIC is computed as
\begin{align*}
\qbic&=\sumj\left\{\log\big|S(X_{t_{j-1}},\hat{\theta}_{n})\big|+\frac{1}{h_{n}}S^{-1}(X_{t_{j-1}},\hat{\theta}_{n})\left[(\D_{j}Y)^{\otimes2}\right]\right\}
\nn\\
&{}\qquad
+\log\Bigg|\frac{1}{2}\sumj\Bigg\{\p_{\theta}^{2}\log\big|S(X_{t_{j-1}},\hat{\theta}_{n})\big|
+\frac{1}{h_{n}}\p_{\theta}^{2}(S^{-1})(X_{t_{j-1}},\hat{\theta}_{n})\left[(\D_{j}Y)^{\otimes2}\right]\Bigg\}\Bigg|.
\end{align*}
Let us consider the conditions for the QBIC to be valid, when $m=r=1$ with $\sigma(x,\theta)=\exp(x^{\prime}\theta/2)$. 
The quasi-likelihood function is then given by
\begin{align}
\mbbh_{n}(\theta)=-\frac{1}{2}\sumj\left\{X_{t_{j-1}}^{\prime}\theta+\frac{1}{h_{n}}(\D_{j}Y)^{2}\exp(-X_{t_{j-1}}^{\prime}\theta)\right\},
\label{VolaLf}
\end{align}
with $-\p_{\theta}^{2}\mbbh_{n}(\theta)=\frac{1}{2}\sumj \frac{(\D_{j}Y)^{2}}{h_{n}}\exp\left(-X_{t_{j-1}}^{\prime}\theta\right)X_{t_{j-1}}X_{t_{j-1}}^{\prime} \ge 0$ a.s.

\begin{ass}
\begin{itemize}
\item[(i)] $\forall q>0$, $\E(|X_{0}|^{q})<\infty$.
\item[(ii)] $\forall q>0$, $\exists C>0$, $\forall s,t\in[0,T]$, $\E(|X_{t}-X_{s}|^{q})<C|t-s|^{q/2}$.
\item[(iii)] $\forall q>0$, $\displaystyle\sup_{0\leq t\leq T}\E(|b_{t}|^{q})<\infty$.
\end{itemize}
\label{Ass8}
\end{ass}

\begin{ass}
\begin{itemize}
\item[(i)] $\displaystyle\sup_{\omega\in\Omega}\sup_{t\leq T}|X_{t}|<\infty$.
\item[(ii)] $\forall L>0$, $\exists C_{L}>0$, $\forall r>0$, $\displaystyle \pr\left\{\lambda_{\min}\bigg(\int_{0}^{T}X_{t}X_{t}^{\prime}dt\bigg)\leq\frac{1}{r}\right\}\leq\frac{C_{L}}{r^{L}}$.
\end{itemize}
\label{Ass9}
\end{ass}

It will be seen that Assumptions \ref{Ass8} and \ref{Ass9} ensure Assumption \ref{Ass1} and inequality (\ref{pldi.2}).

\begin{thm}
Let Assumptions \ref{Ass8} and \ref{Ass9} hold. Then, the assumptions in Theorem \ref{QlfTh2} are satisfied and the corresponding QBIC is given by
\begin{align*}
\qbic&=\sumj\left\{X_{t_{j-1}}^{\prime}\hat{\theta}_{n}+\frac{1}{h_{n}}(\D_{j}Y)^{2}\exp(-X_{t_{j-1}}^{\prime}\hat{\theta}_{n})\right\} \nn\\
&{}\qquad
+\log\bigg|\frac{1}{2h_{n}}\sumj (\D_{j}Y)^{2}\exp(-X_{t_{j-1}}^{\prime}\hat{\theta}_{n})X_{t_{j-1}}X_{t_{j-1}}^{\prime}\bigg|.
\end{align*}
\label{VolaCor1}
\end{thm}

\section{Model-selection consistency}\label{hm:qbic.msc}

As long as concerned with good prediction performance, model-selection consistency itself does not matter in an essential way. 
Given s model set, it does when attempting to find the one ``closest'' (in the sense of KL divergence) to the true data-generating model structure itself as much as possible; 
for example, estimation of daily integrated volatility in econometrics would be the case, 
for econometricians usually builds up daily-volatility prediction model through a time series model such as, among others, ARFIMA models. This section is devoted to studying the validity of model-selection consistency in our general setting. In particular, we propose an adaptive (stepwise) model selection strategy when we have more than one scaling rate. 
We start with a single-norming case. Then, before moving on the multi-scaling case, we look at the case of ergodic diffusions 
since it well illustrates the proposed method.

\subsection{Single-scaling case}

We first consider cases where
\begin{equation}
a_{n} = a_{m,k,n}(\theta_{0}) \to 0
\nonumber
\end{equation}
for each $m\in\{1,\ldots,M\}$ and $k\in\{1,\ldots,K_{m}\}$. Suppose that there exists a random function $\mbbh_{m,0}$ such that
\begin{align}
a_{n}^{2}\mbbh_{m,n}(\theta_{m})\cip\mbbh_{m,0}(\theta_{m}) \label{LLN1}
\end{align}
as $n\to\infty$. 
Moreover, we assume that the optimal parameter $\theta_{m,0}$ in the model $\mcm_{m}$ is the unique maximizer of $\mbbh_{m,0}$:
\begin{align*}
\{\theta_{m,0}\}=\argmax_{\theta_{m}\in\Theta_{m}}\mbbh_{m,0}(\theta_{m}).
\end{align*}
If $m_{0}$ satisfies
\begin{align*}
\{m_{0}\}=\argmin_{m\in\mathfrak{M}}\dim(\Theta_{m}),
\end{align*}
where $\mathfrak{M}=\argmax_{1\leq m\leq M}\mbbh_{m,0}(\theta_{m,0})$, we say that $\mcm_{m_{0}}$ is the {\it optimal model}.
That is, the optimal model is, if exists, an element of the optimal model set $\mathfrak{M}$ which has the smallest dimension.

Let $\Theta_{i}\subset\mbbr^{p_{i}}$ and $\Theta_{j}\subset\mbbr^{p_{j}}$ be the parameter space associated with $\mcm_{i}$ and $\mcm_{j}$, respectively.
We say that $\Theta_{i}$ is nested in $\Theta_{j}$ when $p_{i}<p_{j}$ and there exist a matrix $F\in\mbbr^{p_{j}\times p_{i}}$ with $F^{\prime}F=I_{p_{i}\times p_{i}}$ as well as a $c\in\mbbr^{p_{j}}$ such that $\mbbh_{i,n}(\theta_{i})=\mbbh_{j,n}(F\theta_{i}+c)$ for all $\theta_{i}\in\Theta_{i}$.
That is, when $\Theta_{i}$ is nested in $\Theta_{j}$, any model given by a parameter in $\Theta_{i}$ can also be generated by a parameter in $\Theta_{j}$, so that $\mcm_{j}$ includes $\mcm_{i}$.

\begin{thm}
Assume that $\mcm_{m_{0}}$ is the optimal model. 
Let $m\in\{1,\ldots,M\}\backslash\{m_{0}\}$, and let Assumptions \ref{Ass1} to \ref{Ass3} hold, and suppose that either
\begin{itemize}
\item[{\rm (i)}] $\Theta_{m_{0}}$ is nested in $\Theta_{m}$, or
\item[{\rm (ii)}] $\mbbh_{m,0}(\theta_{m})\neq\mbbh_{m_{0},0}(\theta_{m_{0},0})$ a.s. for any $\theta_{m}\in\Theta_{m}$.
\end{itemize}
Then we have
\begin{align}
& \lim_{n\to\infty}\pr\left(\qbic^{(m_{0})}-\qbic^{(m)}<0\right)=1, \label{msc1.1} \\
& \lim_{n\to\infty}\pr\left(\bic^{(m_{0})}-\bic^{(m)}<0\right)=1. \label{msc1.2}
\end{align}
\label{th.msc}
\end{thm}

This theorem indicates that the probability that QBIC and BIC choose the optimal model tends to $1$ as $n\to\infty$.

\subsection{Multi-scaling case: adaptive model comparison} \label{mul.msc}

For simplicity of exposition we consider the two-scaling case, i.e. $K=2$.
We propose a multi-step model-selection procedure, which seems natural and more effective especially when an adaptive estimation procedure is possible in such a way that we can estimate the first component of $\theta_{m,1}$ without knowledge of the second one $\theta_{m,2}$. That is to say, it should be possible to select an optimal ``partial'' model structure associated with $\theta_{m,1}$, with regarding that with $\theta_{m,2}$ as a nuisance element. 

Suppose that the full model is ``decomposed'' into two parts, each consisting of $M_{1}$ and $M_{2}$ candidates, 
resulting in $M_{1}\times M_{2}$ models in total. Write $(\mcm_{m_{1},m_{2}})_{m_{1}\le M_{1}; m_{2}\le M_{2}}$ for the set of all the candidate models.
We are given the ``full'' quasi-log likelihood function $\mbbh_{m_{1},m_{2},n}(\theta_{m_{1}},\theta_{m_{2}})$.
Roughly speaking, we goes as follows:
\begin{itemize}
\item First, introducing an auxiliary quasi-log likelihood associated with the first-component parameter $\theta_{1}$ ({\it without} involving $\theta_{2}$) 
and compare the corresponding (Q)BICs to select one optimal index, say $m^{\ast}_{1,n}\in\{1,\dots,M_{1}\}$, 
reducing the model-candidate set from $(\mbbh_{m_{1},m_{2},n})_{m_{1},m_{2}}$ to $(\mbbh_{m^{\ast}_{1,n},m_{2},n})_{m_{2}}$;
\item Second, based on the ``partly optimized'' quasi-log likelihoods $\mbbh_{m^{\ast}_{1,n},1,n},\dots,\mbbh_{m^{\ast}_{1,n},M_{2},n}$, 
we find a second-stage optimal index $m^{\ast}_{2}\in\{1,\dots,M_{2}\}$ through the (Q)BIC again;
\item Finally we pick the model $\mcm_{m_{1}^{\ast},m_{2}^{\ast}}$ as our final optimal model.
\end{itemize}
This adaptive procedure apparently reduces the computational cost (the number of comparison) to much extent compared with the joint-(Q)BIC case, 
i.e. from ``$O(M_{1}\times M_{2})$'' to ``$O(M_{1}+M_{2})$''; needless to say, the amount of reduction becomes larger for $K\ge 3$. 

\begin{rem}{\rm 
I would not be essential in the above argument that the final step is based on the original quasi-log likelihood $\mbbh_{m_{1},m_{2},n}$. 
What is essential for the model-selection consistency is that at each stage we have a suitable auxiliary quasi-likelihood function based on which we can estimate the suitably separated optimal model.
We here do not go into this direction.
}\qed\end{rem}

To be specific, we here focus on the ergodic diffusion discussed in Section \ref{GQMLE2}, and then briefly mention the general case in Section \ref{hm:sec_msc.general}.

\subsubsection{Example: ergodic diffusion}\label{hm:sec_msc.ergo.diff}

Here we consider the same setting as in Section \ref{GQMLE2}, that is, the model $\mcm_{m_{1},m_{2}}$ is given by \eqref{mod.erg}:
\begin{align}
dX_{t}=a_{m_{2}}(X_{t},\theta_{m_{2}})dt+b_{m_{1}}(X_{t},\theta_{m_{1}})dw_{t},\quad t\in[0,T_{n}],\;X_{0}=x_{0}.
\nn
\end{align}
Let $B_{m_{1}}(x,\theta_{m_{1}}):=b_{m_{1}}(x,\theta_{m_{1}})b_{m_{1}}(x,\theta_{m_{1}})^{\prime}$.
Up to an additive constant term, the quasi-likelihood function $\mbbh_{m_{1},m_{2},n}$ based on the local-Gauss approximation is given by
\begin{align}
\mbbh_{m_{1},m_{2},n}(\theta_{m_{1},m_{2}})&=-\frac{1}{2}\sumj\bigg\{\log\big|B_{m_{1}}(X_{t_{j-1}},\theta_{m_{1}})\big| \nn\\
&\qquad\qquad+\frac{1}{h_{n}}B_{m_{1}}(X_{t_{j-1}},\theta_{m_{1}})^{-1}\left[\big(\D_{j}X-h_{n}a_{m_{2}}(X_{t_{j-1}},\theta_{m_{2}})\big)^{\otimes2}\right]\bigg\}.
\label{hm:sqlf.ed}
\end{align}
Then,
\begin{align}
\frac{1}{n}\mbbh_{m_{1},m_{2},n}(\theta_{m_{1},m_{2}})&\cip-\frac{1}{2}\int_{\mbbr^{d}}\left[\tr\left\{B(x)B_{m_{1}}(x,\theta_{m_{1}})^{-1}\right\}+\log|B_{m_{1}}(x,\theta_{m_{1}})|\right]\nu(dx) \notag\\
&=:\mbbh_{m_{1},0}^{1}(\theta_{m_{1}}) \label{LNN2}
\end{align}
uniformly in $\theta_{m_{1}}$, where $B(x):=b(x)b(x)^{\prime}$.
We assume that the optimal parameter $\theta_{m_{1},0}$ and $m_{1,0}$ satisfy
\begin{align*}
\{\theta_{m_{1},0}\}&=\argmax_{\theta_{m_{1}}\in\Theta_{m_{1}}}\mbbh_{m_{1},0}^{1}(\theta_{m_{1}}), \\
\{m_{1,0}\}&=\argmin_{m_{1}\in\mathfrak{M}_{1}}\dim(\Theta_{m_{1}}),
\end{align*}
respectively. Here $\mathfrak{M}_{1}=\argmax_{1\leq m_{1}\leq M_{1}}\mbbh_{m_{1},0}^{1}(\theta_{m_{1},0})$.
Furthermore, assume that 
\begin{align}
& \frac{1}{nh_{n}}\big\{ \mbbh_{m_{1},m_{2},n}(\theta_{m_{1},m_{2}}) - \mbbh_{m_{1},m_{2},n}(\theta_{m_{1}}, \theta_{m_{2},0}) \big\} \nn\\
&\cip-\frac{1}{2}\int_{\mbbr^{d}}B_{m_{1}}(x,\theta_{m_{1}})^{-1}\left[\big(a(x)-a_{m_{2}}(x,\theta_{m_{2}})\big)^{\otimes2}\right]\nu(dx)
=:\mbbh_{m_{1},m_{2},0}(\theta_{m_{1},m_{2}}) \label{LLN3}
\end{align}
uniformly in $\theta_{m_{1},m_{2}}$, and that the optimal parameter $\theta_{m_{2},0}$ is the unique maximizer of $\mbbh_{m_{1,0},m_{2},0}$:
\begin{align*}
\{\theta_{m_{2},0}\}=\argmax_{\theta_{m_{2}}\in\Theta_{m_{2}}}\mbbh_{m_{1,0},m_{2},0}(\theta_{m_{1,0},0},\theta_{m_{2}}).
\end{align*}
If $m_{2,0}$ satisfies
\begin{align*}
\{m_{2,0}\}=\argmin_{m_{2}\in\mathfrak{M}_{2}}\dim(\Theta_{m_{2}}),
\end{align*}
where $\mathfrak{M}_{2}=\argmax_{1\leq m_{2}\leq M_{2}}\mbbh_{m_{1,0},m_{2},0}(\theta_{m_{1,0},0},\theta_{m_{2},0})$, we say that $\mcm_{m_{1,0},m_{2,0}}$ is the optimal model.
Since we consider a set of correctly specified models, it holds that $b_{m_{1,0}}(\cdot,\theta_{m_{1,0},0})=b(\cdot)$ and $a_{m_{2,0}}(\cdot,\theta_{m_{2,0},0})=a(\cdot)$.

Let $\Theta_{i_{1}}\times\Theta_{i_{2}}\subset\mbbr^{p_{i_{1}}}\times\mbbr^{p_{i_{2}}}$, $\Theta_{j_{1}}\times\Theta_{j_{2}}\subset\mbbr^{p_{j_{1}}}\times\mbbr^{p_{j_{2}}}$ be the parameter spaces associated with $\mcm_{i_{1},i_{2}}$ and $\mcm_{j_{1},j_{2}}$, respectively.
If $p_{i_{1}}<p_{j_{1}}$ and there exists a matrix $F_{1}\in\mbbr^{p_{j_{1}}\times p_{i_{1}}}$ with $F_{1}^{\prime}F_{1}=I_{p_{i_{1}}\times p_{i_{1}}}$ as well as a $c_{1}\in\mbbr^{p_{j_{1}}}$ such that $\mbbh_{i_{1},m_{2},n}(\theta_{i_{1}},\theta_{m_{2}})=\mbbh_{j_{1},m_{2},n}(F_{1}\theta_{i_{1}}+c_{1},\theta_{m_{2}})$ for all $\theta_{i_{1}}\in\Theta_{i_{1}}$ and $m_{2}\in\{1,\ldots,M_{2}\}$, we say that $\Theta_{i_{1}}$ is nested in $\Theta_{j_{1}}$.
It is defined in a similar manner that $\Theta_{i_{2}}$ is nested in $\Theta_{j_{2}}$.

\medskip

First we consider the \textit{joint QBIC} of \eqref{mQBIC_def4}:
\begin{align*}
\qbic^{(m_{1},m_{2})}&=-2\mbbh_{m_{1},m_{2},n}\big(\hat{\theta}_{m_{1},m_{2},n}\big) \\
&\qquad+\log\big|-\p_{\theta_{m_{1}}}^{2}\mbbh_{m_{1},m_{2},n}(\hat{\theta}_{m_{1},m_{2},n})\big|+\log\big|-\p_{\theta_{m_{2}}}^{2}\mbbh_{m_{1},m_{2},n}(\hat{\theta}_{m_{1},m_{2},n})\big|.
\end{align*}
If $(m_{1,n}^{\ast},m_{2,n}^{\ast})=\argmin_{(m_{1},m_{2})\in\{1,\ldots,M_{1}\}\times\{1,\ldots,M_{2}\}}\qbic^{(m_{1},m_{2})}$, we choose the model $\mcm_{m_{1,n}^{\ast},m_{2,n}^{\ast}}$, which we again call the {\it optimal model}, as the optimal model among the candidate models. The details of $\qbic^{(m_{1},m_{2})}$ is given in Theorems \ref{ErTh1} and \ref{ErTh2}.
Likewise, we define
\begin{align*}
\bic^{(m_{1},m_{2})}&=-2\mbbh_{m_{1},m_{2},n}\big(\hat{\theta}_{m_{1},m_{2},n}\big) 
+ p_{m_{1}}\log n + p_{m_{2}}\log T_{n}.
\end{align*}

\begin{thm}
Assume that $\mcm_{m_{1,0},m_{2,0}}$ is the optimal model. Let $(m_{1},m_{2})\in(\{1,\ldots,M_{1}\}\backslash\{m_{1,0}\})\times(\{1,\ldots,M_{2}\}\backslash\{m_{2,0}\})$. 
Let Assumptions \ref{Ass4} to \ref{Ass7} hold and suppose that $\Theta_{m_{1,0}}$ and $\Theta_{m_{2,0}}$ are nested in $\Theta_{m_{1}}$ and $\Theta_{m_{2}}$, respectively.
Then we have
\begin{align*}
\lim_{n\to\infty}\pr\left(\qbic^{(m_{1,0},m_{2,0})}-\qbic^{(m_{1},m_{2})}<0\right)=1,
\end{align*}
and the same statement with ``{\rm QBIC}'' replaced by ``{\rm BIC}''.
\label{th.msc2}
\end{thm}

\medskip

Next we turn to the \textit{two-step QBIC}. 
In the present case, we apply the previous single-scaling result twice for the single true data-generating model. 
First, we focus on the diffusion coefficient, which we can estimate more quickly than the drift one. 
Under suitable conditions, $\qbic^{(m_{1})}$ and $\bic^{(m_{1})}$ are given by
\begin{align*}
\qbic^{(m_{1})}&=-2\mbbh_{m_{1},n}^{1}(\hat{\theta}_{m_{1},n})
+\log\big|-\p_{\theta_{m_{1}}}^{2}\mbbh_{m_{1},n}^{1}(\hat{\theta}_{m_{1},n})\big|, \\
\bic^{(m_{1})}&=-2\mbbh_{m_{1},n}^{1}(\hat{\theta}_{m_{1},n})+p_{m_{1}}\log n,
\end{align*}
where $\mbbh_{m_{1},n}^{1}$ is defined by the joint quasi-likelihood \eqref{hm:sqlf.ed} with $a_{m_{2}}$ being null:
\begin{equation}
\mbbh^{1}_{m_{1},n}(\theta_{m_{1},m_{2}}) = -\frac{1}{2}\sumj\bigg\{\log\big|B_{m_{1}}(X_{t_{j-1}},\theta_{m_{1}})\big|
+\frac{1}{h_{n}}B_{m_{1}}(X_{t_{j-1}},\theta_{m_{1}})^{-1}\big[\big(\D_{j}X\big)^{\otimes2}\big]\bigg\},
\nn
\end{equation}
and where $\hat{\theta}_{m_{1},n}$ is the QMLE associated with $\mbbh_{m_{1},n}^{1}$. Note that we can write
\begin{equation}
\frac{1}{n}\mbbh_{n}(\theta) = \frac{1}{n}\mbbh^{1}_{n}(\theta_{1}) + \del^{1}_{n}(\theta)
\nonumber
\end{equation}
with $\sup_{\theta}|\del^{1}_{n}(\theta)|\cip 0$.
We proceed as follows. 
\begin{itemize}
\item First, assuming that $\Theta_{i_{1}}$ is nested in $\Theta_{j_{1}}$, i.e. $\mbbh_{i_{1},n}^{1}(\theta_{i_{1}})=\mbbh_{j_{1},n}^{1}(F_{1}\theta_{i_{1}}+c_{1})$, 
we set $\{m_{1,n}^{\ast}\}=\argmin_{1\leq m_{1}\leq M_{1}}\qbic^{(m_{1})}$.
\item Next, we consider the stochastic differential equation
\begin{align}
dX_{t}=a_{m_{2}}(X_{t},\theta_{m_{2}})dt+b_{m_{1,n}^{\ast}}(X_{t},\hat{\theta}_{m_{1,n}^{\ast},n})dw_{t}. \label{drif.part}
\end{align}
Assuming that $\{\hat{\theta}_{m_{2},n}\}=\argmax_{\theta_{m_{2}}\in\Theta_{m_{2}}}\mbbh_{m_{1,n}^{\ast},m_{2},n}(\hat{\theta}_{m_{1,n}^{\ast},n},\theta_{m_{2}})$ and that $\{m_{2,n}^{\ast}\}=\argmin_{1\leq m_{2}\leq M_{2}}$ $\qbic^{(m_{2}|m_{1,n}^{\ast})}$, where
\begin{align*}
\qbic^{(m_{2}|m_{1,n}^{\ast})}=-2\mbbh_{m_{1,n}^{\ast},m_{2},n}\big(\hat{\theta}_{m_{1,n}^{\ast},n},\hat{\theta}_{m_{2},n}\big)+\log\big|-\p_{\theta_{m_{2}}}^{2}\mbbh_{m_{1,n}^{\ast},m_{2},n}\big(\hat{\theta}_{m_{1,n}^{\ast},n},\hat{\theta}_{m_{2},n}\big)\big|,
\end{align*}
we select the model $\mcm_{m_{1,n}^{\ast},m_{2,n}^{\ast}}$ as the final best model.
\end{itemize}
When we use BIC, the best model is selected by a similar procedure and $\bic^{(m_{2}|m_{1,n}^{\ast})}$ is given by
\begin{align*}
\bic^{(m_{2}|m_{1,n}^{\ast})}=-2\mbbh_{m_{1,n}^{\ast},m_{2},n}\big(\hat{\theta}_{m_{1,n}^{\ast}},\hat{\theta}_{m_{2},n}\big)+p_{m_{2}}\log T_{n}.
\end{align*}

Joint (Q)BIC and two-step (Q)BIC may select different models for a fixed sample size, however, the model-selection consistency property is asymptotically shared.

\begin{thm}
Assume that $\mcm_{m_{1,0},m_{2,0}}$ is the optimal model. Let $(m_{1},m_{2})\in(\{1,\ldots,M_{1}\}\backslash\{m_{1,0}\})\times(\{1,\ldots,M_{2}\}\backslash\{m_{2,0}\})$.
Assume that Assumptions \ref{Ass4}-\ref{Ass7} hold and that $\Theta_{m_{1,0}}$ and $\Theta_{m_{2,0}}$ are nested in $\Theta_{m_{1}}$ and $\Theta_{m_{2}}$, respectively.
Then we have
\begin{align*}
&\lim_{n\to\infty}\pr\left(\qbic^{(m_{1,0})}-\qbic^{(m_{1})}<0\right)=1, \\
&\lim_{n\to\infty}\pr\left(\qbic^{(m_{2,0}|m_{1,n}^{\ast})}-\qbic^{(m_{2}|m_{1,n}^{\ast})}<0\right)=1,
\end{align*}
and the same statements with ``{\rm QBIC}'' replaced by ``{\rm BIC}''.
\label{th.msc3}
\end{thm}

\subsubsection{Remark on general case}\label{hm:sec_msc.general}

Without essential change, we may follow the same scenario as in the previous section for general LAQ models under the setting described in the beginning of this section: instead of the original ``full'' quasi-likelihood $\mbbh_{n}(\theta)$, we solely look at some ``auxiliary'' random fields $\theta_{1}\mapsto\mbbh^{1}_{n}(\theta_{1})$ and $\theta_{2}\mapsto\mbbh_{n}(\hat{\theta}_{1,n},\theta_{2})$ in this order, based on which we successively define the two-step QMLE as $\hat{\theta}_{1,n}\in\argmax\mbbh_{n}^{1}$ and $\hat{\theta}_{2,n}\in\argmax_{\theta_{2}}\mbbh_{n}^{2}(\hat{\theta}_{1,n},\theta_{2})$. More specifically, we let Assumptions \ref{Ass1} to \ref{Ass3} hold and further assume that there exist positive sequences $a_{i,n}(\tz)\to 0$ ($i=1,2$) such that $a_{1,n}(\tz)/a_{2,n}(\tz)\to 0$ and random functions $\mbby^{1}_{0}(\theta_{1})$ and $\mbby^{2}_{0}(\theta)$, for which:
\begin{enumerate}
\item $a_{1,n}^{2}(\tz)\mbbh_{n}(\theta) = a_{1,n}^{2}\mbbh^{1}_{n}(\theta_{1}) + \del^{1}_{n}(\theta)$ where $\sup_{\theta}|\del^{1}_{n}(\theta)|\cip 0$;
\item $\sup_{\theta_{1}}|a_{1,n}(\tz)^{2}\{\mbbh^{1}_{n}(\theta_{1}) - \mbbh^{1}_{n}(\theta_{1,0})\}-\mbby^{1}_{0}(\theta_{1})|\cip 0$, where $\{\theta_{1,0}\}=\argmax \mbby^{1}_{0}$ a.s.;
\item $\sup_{\theta}|a_{2,n}(\tz)^{2}\{\mbbh_{n}(\theta) - \mbbh_{n}(\theta_{1},\theta_{2,0})\}-\mbby^{2}_{0}(\theta)|\cip 0$, where $\{\theta_{2,0}\}=\argmax_{\theta_{2}}\mbby^{2}_{0}(\theta_{1,0},\theta_{2})$ a.s.
\end{enumerate}
Under these conditions we may deduce the consistency $\tes=(\hat{\theta}_{1,n},\hat{\theta}_{2,n})\cip\tz$, which combined with Assumption \ref{Ass1} gives the tightness $A_{n}^{-1}(\tz)(\tes-\tz)=O_{p}(1)$; in this case, the LAQ structure \eqref{hm:pw.LAQ} typically takes the form
\begin{equation}
\sup_{u=(u_{1},u_{2})\in A}\bigg|\log\mbbz_{n}(u) - \bigg( \D_{1,n}[u_{1}]+\D_{2,n}[u_{2}] -\frac{1}{2}(\Gam_{1,0}(\theta_{1,0})[u_{1},u_{1}]+\Gam_{2,0}(\tz)[u_{2},u_{2}]) \bigg)\bigg| = o_{p}(1)
\nn
\end{equation}
for each compact $A\in\mbbr^{p}$, with $\D_{1}=:(\D_{1,n}, \D_{2,n})$ and $\Gam_{0}=:\diag\{\Gam_{1,0}(\theta_{1,0}),\,\Gam_{2,0}(\tz)\}$.
Through a similar way to the proof of Theorems \ref{th.msc}, \ref{th.msc2} and \ref{th.msc3}, we can prove the model-selection consistency.


\section{Simulation results} \label{Simu}

We here conduct a number of simulations to observe finite-sample performance of (Q)BIC proposed in this paper. While what to be looked at is quasi-Bayes factors for candidate models, for conciseness we focus on the selection frequency as well as the estimation performance of the quasi-maximum likelihood estimates.
In Section \ref{Simu1}, we use the R package {\tt yuima} (Brouste {\it et al.} \cite{YUIMA14}) for generating data and estimating the parameter.
Moreover, we use {\tt optim} at software R to estimate the parameters in Sections \ref{Simu2} and \ref{Simu3}.
We set the initial value in numerical optimization to be random numbers drawn from uniform distribution $U(\theta_{i,0}-0.5,\,\theta_{i,0}+0.5)$, $1\le i\le p$. 
All the SDE coefficients considered here are of the partially convex type mentioned in Remark \ref{hm:ex.separate.convex.diffusion}.


\subsection{Ergodic diffusion process} \label{Simu1}

The data $\mbX_{n}=(X_{t_{j}})_{j=0}^{n}$ with $t_{j}=jn^{-2/3}$ and the number of data $n$ are obtained from the true model given by 
\begin{align*}
dX_{t}=-X_{t}dt+\exp\Big\{\frac{1}{2}(-2\cos X_{t}+1)\Big\}dw_{t},\quad t\in[0,T_{n}],\quad X_{0}=1,
\end{align*}
where $w$ is an $1$-dimensional standard Wiener process and $T_{n}=n^{1/3}$.
We consider the following as the diffusion part:
\begin{align*}
&{\bf Diff}\;{\bf 1:}\exp\Big\{\frac{1}{2}(\theta_{11}\cos X_{t}+\theta_{12}\sin X_{t}+\theta_{13})\Big\};
\;{\bf Diff}\;{\bf 2:}\exp\Big\{\frac{1}{2}(\theta_{11}\cos X_{t}+\theta_{12}\sin X_{t})\Big\}; \\
&{\bf Diff}\;{\bf 3:}\exp\Big\{\frac{1}{2}(\theta_{11}\cos X_{t}+\theta_{13})\Big\};
\;{\bf Diff}\;{\bf 4:}\exp\Big\{\frac{1}{2}(\theta_{12}\sin X_{t}+\theta_{13})\Big\}; \\
&{\bf Diff}\;{\bf 5:}\exp\Big\{\frac{1}{2}\theta_{11}\cos X_{t}\Big\};
\;{\bf Diff}\;{\bf 6:}\exp\Big\{\frac{1}{2}\theta_{12}\sin X_{t}\Big\};
\;{\bf Diff}\;{\bf 7:}\exp\Big\{\frac{1}{2}\theta_{13}\Big\}.
\end{align*}
Moreover, we assume the following as the drift part: 
\begin{align*}
{\bf Drif}\;{\bf 1:}\;\theta_{21}X_{t}+\theta_{22};
\;{\bf Drif}\;{\bf 2:}\;\theta_{21}X_{t};
\;{\bf Drif}\;{\bf 3:}\;\theta_{22}.
\end{align*}
The candidate models is given by the combination of the diffusion part and drift part.
For example, in the case of Diff 1 and Drif 1, we consider the statistical model
\begin{align*}
dX_{t}=(\theta_{21}X_{t}+\theta_{22})dt+\exp\Big\{\frac{1}{2}(\theta_{11}\cos X_{t}+\theta_{12}\sin X_{t}+\theta_{13})\Big\}dw_{t}.
\end{align*}
That is, the true model consists of Diff 3 and Drif 2. 

We compare model-selection frequency through QBIC, BIC, and the contrast-based information criterion (CIC), which is an AIC-type criterion introduced by \cite{Uch10} under the rapidly increasing experimental design $nh_{n}\to 0$ 
(see also Fujii and Uchida \cite{FujUch14} for CIC under a weaker sampling-design condition $nh^{q}\to 0$ for some $q\ge 2$). 
We simulate the number of the model selected by using joint QBIC, joint BIC, two-step QBIC, two-step BIC and CIC among the candidate models based on 1000 sample paths. The simulations are done for each $n=1000,3000,5000$.

Tables \ref{ms.joint} and \ref{ms.two-step} summarize the comparison results of the model-selection frequencies; they show quite similar tendencies, in particular, the frequencies that the model defined by Diff 3 and Drif 2 is selected by QBIC and BIC become larger as $n$ increases. Also noted is that BIC often takes values between QBIC and CIC; in particular, 
QBIC chooses the full model consisting of Diff 1 and Drif 1 more frequently than BIC. Moreover, joint (Q)BIC gets close to two-step (Q)BIC as $n$ increases.

It is worth mentioning that computation time of joint (Q)BIC was overall about twice of that of two-step (Q)BIC. This superiority of the two-step (Q)BIC should become more significant for higher-dimensional models.

\begin{table}[t]
\begin{center}
\caption{\footnotesize The number of models selected by joint QBIC, joint BIC and CIC in Section \ref{Simu1} over 1000 simulations for various $n$}
\begin{tabular}{| l | c || c c c c c c c |} \hline
\multicolumn{1}{| l |}{} & \multicolumn{1}{c ||}{Criteria} & \multicolumn{7}{c|}{$n=1000$} \\
 & & Diff 1 & Diff 2 & Diff $3^{\ast}$ & Diff 4 & Diff 5 & Diff 6 & Diff 7 \\ \hline
 & QBIC & 7 & 8 & 109 & 1 & 15 & 0 & 1 \\
Drif 1 & BIC & 0 & 20 & 105 & 1 & 49 & 0 & 2 \\ 
 & CIC & 25 & 23 & 136 & 3 & 19 & 0 & 2 \\ \hline
 & QBIC & 19 & 17 & {\bf 741} & 0 & 76 & 0 & 1 \\
Drif $2^{\ast}$ & BIC & 1 & 22 & {\bf 523} & 0 & 248 & 0 & 1 \\ 
 & CIC & 92 & 43 & {\bf 559} & 0 & 73 & 0 & 1 \\ \hline
 & QBIC & 0 & 0 & 5 & 0 & 0 & 0 & 0 \\
Drif 3 & BIC & 0 & 0 & 28 & 0 & 0 & 0 & 0 \\ 
 & CIC & 5 & 0 & 19 & 0 & 0 & 0 & 0 \\ \hline\hline
\multicolumn{1}{| l |}{} & \multicolumn{1}{c ||}{} & \multicolumn{7}{c|}{$n=3000$} \\
 & & Diff 1 & Diff 2 & Diff $3^{\ast}$ & Diff 4 & Diff 5 & Diff 6 & Diff 7 \\ \hline
 & QBIC & 1 & 2 & 102 & 0 & 0 & 0 & 0 \\
Drif 1 & BIC & 0 & 2 & 126 & 0 & 10 & 0 & 0 \\ 
 & CIC & 24 & 5 & 173 & 0 & 2 & 0 & 0 \\ \hline
 & QBIC & 12 & 4 & {\bf 867} & 0 & 12 & 0 & 0 \\
Drif $2^{\ast}$ & BIC & 1 & 4 & {\bf 786} & 0 & 63 & 0 & 0 \\ 
 & CIC & 110 & 6 & {\bf 667} & 0 & 7 & 0 & 0 \\ \hline
 & QBIC & 0 & 0 & 0 & 0 & 0 & 0 & 0 \\
Drif 3 & BIC & 0 & 0 & 8 & 0 & 0 & 0 & 0 \\ 
 & CIC & 0 & 0 & 6 & 0 & 0 & 0 & 0 \\ \hline\hline
\multicolumn{1}{| l |}{} & \multicolumn{1}{c ||}{} & \multicolumn{7}{c|}{$n=5000$} \\
 & & Diff 1 & Diff 2 & Diff $3^{\ast}$ & Diff 4 & Diff 5 & Diff 6 & Diff 7 \\ \hline
 & QBIC & 1 & 0 & 80 & 0 & 0 & 0 & 0 \\
Drif 1 & BIC & 0 & 0 & 113 & 0 & 3 & 0 & 0 \\ 
 & CIC & 30 & 1 & 166 & 0 & 2 & 0 & 0 \\ \hline
 & QBIC & 16 & 0 & {\bf 900} & 0 & 3 & 0 & 0 \\
Drif $2^{\ast}$ & BIC & 1 & 0 & {\bf 863} & 0 & 20 & 0 & 0 \\ 
 & CIC & 135 & 0 & {\bf 666} & 0 & 7 & 0 & 0 \\ \hline
 & QBIC & 0 & 0 & 0 & 0 & 0 & 0 & 0 \\
Drif 3 & BIC & 0 & 0 & 8 & 0 & 0 & 0 & 0 \\ 
 & CIC & 0 & 0 & 0 & 0 & 0 & 0 & 0 \\ \hline
\end{tabular}
\label{ms.joint}
\end{center}
\end{table}

\begin{table}[t]
\begin{center}
\caption{\footnotesize The number of models selected by two-step QBIC and two-step BIC in Section \ref{Simu1} over 1000 simulations for various $n$}
\begin{tabular}{| l | c || c c c c c c c |} \hline
\multicolumn{1}{| l |}{} & \multicolumn{1}{c ||}{Criteria} & \multicolumn{7}{c|}{$n=1000$} \\
 & & Diff 1 & Diff 2 & Diff $3^{\ast}$ & Diff 4 & Diff 5 & Diff 6 & Diff 7 \\ \hline
\multirow{2}{*}{Drif 1}  & QBIC & 4 & 3 & 108 & 0 & 0 & 0 & 0 \\ 
 & BIC & 1 & 6 & 120 & 0 & 41 & 0 & 0 \\ \hline
\multirow{2}{*}{Drif $2^{\ast}$} & QBIC & 19 & 10 & {\bf 798} & 0 & 45 & 0 & 0 \\ 
 & BIC & 3 & 16 & {\bf 588} & 0 & 199 & 0 & 0 \\ \hline
\multirow{2}{*}{Drif 3} & QBIC & 0 & 0 & 1 & 0 & 0 & 0 & 0 \\ 
 & BIC & 0 & 0 & 26 & 0 & 0 & 0 & 0 \\ \hline\hline
\multicolumn{1}{| l |}{} & \multicolumn{1}{c ||}{Criteria} & \multicolumn{7}{c|}{$n=3000$} \\
 & & Diff 1 & Diff 2 & Diff $3^{\ast}$ & Diff 4 & Diff 5 & Diff 6 & Diff 7 \\ \hline
\multirow{2}{*}{Drif 1} & QBIC & 6 & 0 & 77 & 0 & 0 & 0 & 0 \\ 
 & BIC & 0 & 3 & 111 & 0 & 3 & 0 & 0 \\ \hline
\multirow{2}{*}{Drif $2^{\ast}$} & QBIC & 19 & 1 & {\bf 892} & 0 & 4 & 0 & 0 \\ 
 & BIC & 1 & 1 & {\bf 836} & 0 & 36 & 0 & 0 \\ \hline
\multirow{2}{*}{Drif 3} & QBIC & 0 & 0 & 1 & 0 & 0 & 0 & 0 \\ 
 & BIC & 0 & 0 & 9 & 0 & 0 & 0 & 0  \\ \hline\hline
\multicolumn{1}{| l |}{} & \multicolumn{1}{c ||}{Criteria} & \multicolumn{7}{c|}{$n=5000$} \\
 & & Diff 1 & Diff 2 & Diff $3^{\ast}$ & Diff 4 & Diff 5 & Diff 6 & Diff 7 \\ \hline
\multirow{2}{*}{Drif 1} & QBIC & 1 & 0 & 80 & 0 & 0 & 0 & 0 \\ 
 & BIC & 0 & 3 & 115 & 0 & 1 & 0 & 0 \\ \hline
\multirow{2}{*}{Drif $2^{\ast}$} & QBIC & 14 & 0 & {\bf 904} & 0 & 1 & 0 & 0 \\ 
 & BIC & 2 & 0 & {\bf 864} & 0 & 18 & 0 & 0 \\ \hline
\multirow{2}{*}{Drif 3} & QBIC & 0 & 0 & 0 & 0 & 0 & 0 & 0 \\ 
 & BIC & 0 & 0 & 0 & 0 & 0 & 0 & 0 \\ \hline
\end{tabular}
\label{ms.two-step}
\end{center}
\end{table}


\subsection{Volatility-parameter estimation for continuous semimartingale} \label{Simu2}
Let $(X_{t_{j}},Y_{t_{j}})_{j=0}^{n}$ be a data set with $t_{j}=j/n$ and the number of data $n$. We simulate 1000 data sets from the stochastic integral equation
\begin{align*}
dY_{t}=\exp\Big(\frac{1}{2}X_{t}^{\prime}\tz\Big)dw_{t}
=\exp\left\{\frac{1}{2}(-2X_{2,t}+3X_{3,t})\right\}dw_{t},\quad t\in[0,1],\quad Y_{0}=0,
\end{align*}
where $X_{t}=(X_{1,t},X_{2,t},X_{3,t})^{\prime}$, the true parameter $\tz=(0,-2,3)^{\prime}$
and $w$ is an $1$-dimensional standard Wiener process. 
We consider the following models:
\begin{align*}
&{\bf Model}\;{\bf 1:}\;dY_{t}=\exp\Big\{\frac{1}{2}(\theta_{1}X_{1,t}+\theta_{2}X_{2,t}+\theta_{3}X_{3,t})\Big\}dw_{t};
\nn\\
&{\bf Model}\;{\bf 2:}\;dY_{t}=\exp\Big\{\frac{1}{2}(\theta_{1}X_{1,t}+\theta_{2}X_{2,t})\Big\}dw_{t}; \; 
{\bf Model}\;{\bf 3:}\;dY_{t}=\exp\Big\{\frac{1}{2}(\theta_{1}X_{1,t}+\theta_{3}X_{3,t})\Big\}dw_{t}; \nn\\
&{\bf Model}\;{\bf 4:}\;dY_{t}=\exp\Big\{\frac{1}{2}(\theta_{2}X_{2,t}+\theta_{3}X_{3,t})\Big\}dw_{t}; \; 
{\bf Model}\;{\bf 5:}\;dY_{t}=\exp\Big\{\frac{\theta_{1}}{2}X_{1,t}\Big\}dw_{t}; \nn\\
&{\bf Model}\;{\bf 6:}\;dY_{t}=\exp\Big\{\frac{\theta_{2}}{2}X_{2,t}\Big\}dw_{t}; \; 
{\bf Model}\;{\bf 7:}\;dY_{t}=\exp\Big\{\frac{\theta_{3}}{2}X_{3,t}\Big\}dw_{t}.
\end{align*}
Then the true model is Model 4. Note that Models 2, 3, 5, 6 and 7 are misspecified models.

For each model, the maximum likelihood estimator is obtained from the quasi-likelihood \eqref{VolaLf}.
In the Model 1 (full model), the statistics QBIC, BIC and formal AIC (fAIC) are given by
\begin{align*}
&\qbic=\sumj\bigg\{\big(\hat{\theta}_{1}X_{1,t_{j-1}}+\hat{\theta}_{2}X_{2,t_{j-1}}+\hat{\theta}_{3}X_{3,t_{j-1}}\big)
\nn\\
&{}\qquad
+n(\D_{j}Y)^{2}\exp\big(-\hat{\theta}_{1}X_{1,t_{j-1}}-\hat{\theta}_{2}X_{2,t_{j-1}}-\hat{\theta}_{3}X_{3,t_{j-1}}\big)\bigg\} \\
&\qquad\qquad+\log\left|\frac{n}{2}\sumj(\D_{j}Y)^{2}\exp\big(-\hat{\theta}_{1}X_{1,t_{j-1}}-\hat{\theta}_{2}X_{2,t_{j-1}}-\hat{\theta}_{3}X_{3,t_{j-1}}\big)X_{t_{j-1}}X_{t_{j-1}}^{\prime}\right|, \\
&\bic=\sumj\bigg\{\big(\hat{\theta}_{1}X_{1,t_{j-1}}+\hat{\theta}_{2}X_{2,t_{j-1}}+\hat{\theta}_{3}X_{3,t_{j-1}}\big)
\nn\\
&{}\qquad+n(\D_{j}Y)^{2}\exp\big(-\hat{\theta}_{1}X_{1,t_{j-1}}-\hat{\theta}_{2}X_{2,t_{j-1}}-\hat{\theta}_{3}X_{3,t_{j-1}}\big)\bigg\}+3\log n, \\
&\mathrm{fAIC}_{n}=\sumj\bigg\{\big(\hat{\theta}_{1}X_{1,t_{j-1}}+\hat{\theta}_{2}X_{2,t_{j-1}}+\hat{\theta}_{3}X_{3,t_{j-1}}\big)
\nn\\
&{}\qquad
+n(\D_{j}Y)^{2}\exp\big(-\hat{\theta}_{1}X_{1,t_{j-1}}-\hat{\theta}_{2}X_{2,t_{j-1}}-\hat{\theta}_{3}X_{3,t_{j-1}}\big)\bigg\}
+3\times2,
\end{align*}
where $\hat{\theta}_{1}$, $\hat{\theta}_{2}$ and $\hat{\theta}_{3}$ are the associated quasi-maximum likelihood estimators.


\subsubsection{Non-random covariate process} \label{Simu2.1}
First we set
\begin{equation}
X_{t_{j}}=\bigg(1,\, \cos\bigg(\frac{2j\pi}{n}\bigg),\, \sin\bigg(\frac{2j\pi}{n}\bigg)\bigg)^{\prime}, \quad j=0,1,\dots,n.
\nonumber
\end{equation}
We readily get
\begin{align*}
\int_{0}^{T}X_{t}X_{t}^{\prime}dt&=\int_{0}^{1}\left(\begin{array}{ccc}
1 & \cos(2t\pi) & \sin(2t\pi) \\ 
\cos(2t\pi) & \cos^{2}(2t\pi) & \cos(2t\pi)\sin(2t\pi) \\ 
\sin(2t\pi) & \cos(2t\pi)\sin(2t\pi) & \sin^{2}(2t\pi)
\end{array}\right)dt
=\left(\begin{array}{ccc}
1 & 0 & 0 \\ 
0 & 1/2 & 0 \\ 
0 & 0 & 1/2
\end{array}\right),
\end{align*}
so that $\det\big(\int_{0}^{T}X_{t}X_{t}^{\prime}dt\big)=\frac{1}{4}$.

In Table \ref{model_sele2}, Model 4 is selected with high frequency as the best model for all cases; 
the model under consideration is completely Gaussian, and the performance are pretty good. 
fAIC tends to choose a model larger than the true one even for large sample size, while QBIC and BIC do show the model-selection consistency. These phenomena are common in model selection, in particular, it should be noted that the model-selection inconsistency is not a defect as the AIC methodology is not intended to estimate the true model consistently.

Table \ref{model_esti2} summarizes the mean and the standard deviation of estimators in each model. In the case of Model 4, 
the estimators get closer to the true value and the standard deviation become smaller when the sample size become larger.

\begin{table}[t]
\begin{center}
\caption{\footnotesize The number of models selected by QBIC, BIC and fAIC in Section \ref{Simu2.1} over 1000 simulations for various $n$ (1-7 express the model labels, and the true model is model 4)}
\begin{tabular}{ | l || c c c c c c c | c c c c c c c | c c c c c c c |} \hline
\multicolumn{1}{| l ||}{Criterion} & \multicolumn{7}{c|}{$n=50$} & \multicolumn{7}{c|}{$n=100$} & \multicolumn{7}{c|}{$n=200$} \\
& 1 & 2 & 3 & $4^{\ast}$ & 5 & 6 & 7 & 1 & 2 & 3 & $4^{\ast}$ & 5 & 6 & 7 & 1 & 2 & 3 & $4^{\ast}$ & 5 & 6 & 7 \\ \hline
QBIC & 74 & 0 & 0 & {\bf 925} & 0 & 0 & 0 & 57 & 0 & 0 & {\bf 943} & 0 & 0 & 0 & 37 & 0 & 0 & {\bf 963} & 0 & 0 & 0 \\
BIC & 67 & 0 & 0 & {\bf 933} & 0 & 0 & 0 & 39 & 0 & 0 & {\bf 961} & 0 & 0 & 0 & 25 & 0 & 0 & {\bf 975} & 0 & 0 & 0 \\
fAIC & 183 & 0 & 0 & {\bf 817} & 0 & 0 & 0 & 178 & 0 & 0 & {\bf 822} & 0 & 0 & 0 & 179 & 0 & 0 & {\bf 821} & 0 & 0 & 0 \\ \hline
\end{tabular}
\label{model_sele2}
\end{center}
\end{table}

\begin{table}[t]
\begin{center}
\caption{\footnotesize The mean and the standard deviation (s.d.) of the estimator $\hat{\theta}_{1}$, $\hat{\theta}_{2}$ and $\hat{\theta}_{3}$ in Section \ref{Simu2.1} for various $n$ (1-7 express the models, and the true parameter $\tz=(0,-2,3)$)}
\smallskip
\begin{tabular}{| l | l || c c c | c c c | c c c |} \hline
\multicolumn{1}{|l|}{} & \multicolumn{1}{l ||}{} & \multicolumn{3}{c|}{$n=50$} & \multicolumn{3}{c|}{$n=100$} & \multicolumn{3}{c|}{$n=200$} \\
 & & $\hat{\theta}_{1}$ & $\hat{\theta}_{2}$ & $\hat{\theta}_{3}$ & $\hat{\theta}_{1}$ & $\hat{\theta}_{2}$ & $\hat{\theta}_{3}$ & $\hat{\theta}_{1}$ & $\hat{\theta}_{2}$ & $\hat{\theta}_{3}$ \\ \hline
1 & mean & -0.0564 & -1.8312 & 3.1129 & -0.0218 & -1.8972 & 3.0748 & -0.0183 & -1.9586 & 3.02986 \\
 & s.d. & 0.2086 & 0.2879 & 0.2978 & 0.1509 & 0.2006 & 0.2052 & 0.1026 & 0.1461 & 0.1442 \\ \hline 
2 & mean & 1.5872 & -1.8314 & -- &1.6145 & -1.8839 & -- & 1.5834 & -1.9602 & --  \\
 & s.d. & 0.3579 & 0.4852 & -- & 0.2497 & 0.3505 & -- & 0.1750 & 0.2525 & -- \\ \hline 
3 & mean & 0.6473 & -- & 3.1168 & 0.7259 & -- & 3.0660 & 0.7734 & -- & 3.0328 \\
 & s.d. & 0.2829 & -- & 0.3981 & 0.2054 & -- & 0.2705 & 0.1486 & -- & 0.1966 \\ \hline 
$4^{\ast}$ & mean & -- & -1.8312 & 3.1129 & -- & -1.8972 & 3.0748 & -- & -1.9586 & 3.0299 \\
 & s.d. & -- & 0.2879 & 0.2978 & -- & 0.2006 & 0.2052 & -- & 0.1461 & 0.1441 \\ \hline 
5 & mean & 0.4871 & -- & -- & 0.5045 & -- & -- & 0.4915 & -- & -- \\
 & s.d. & 0.2858 & -- & -- & 0.2866 & -- & -- & 0.2948 & -- & -- \\ \hline 
6 & mean & -- & -1.892 & -- & -- & -1.916 & -- & -- & -1.9726 & -- \\
 & s.d. & -- & 0.3427 & -- & -- & 0.2814 & -- & -- & 0.2157 & -- \\ \hline 
7 & mean & -- & -- & 3.0867 & -- & -- & 3.0498 & -- & -- & 3.0262 \\
 & s.d. & -- & -- & 0.3026 & -- & -- & 0.2267 & -- & -- & 0.1716 \\ \hline 
\end{tabular}
\label{model_esti2}
\end{center}
\end{table}


\subsubsection{Random covariate process} \label{Simu2.2}

We here consider the following two cases:
\begin{itemize}
\item[{\rm (i--1)}] $X_{t}=(X_{1,t},X_{2,t},X_{3,t})^{\prime}=\big(1,\, \cos(B_{t}),\, \sin(B_{t})\big)^{\prime}$;
\item[{\rm (i--2)}] $X_{t}=(X_{1,t},X_{2,t},X_{3,t})^{\prime}=\big(10,\, \cos(B_{t}),\, \sin(B_{t})\big)^{\prime}$,
\end{itemize}
where $B$ is a one-dimensional standard Wiener process.
For data generation, we use the 3-dimensional stochastic differential equation for $(X_{2,t},X_{3,t},Y_{t})$
\begin{align*}
d\begin{pmatrix}
X_{2,t} \\ X_{3,t} \\ Y_{t}
\end{pmatrix}
=-\frac{1}{2}\begin{pmatrix}
X_{2,t} \\ X_{3,t} \\ 0
\end{pmatrix}dt + 
\begin{pmatrix}
-X_{3,t} & 0 \\ X_{2,t} & 0 \\ 0 & \exp\big\{\frac{1}{2}(-2X_{2,t}+3X_{3,t})\big\}
\end{pmatrix}
d\binom{B_{t}}{w_{t}}.
\end{align*}

Tables \ref{model_sele3} and \ref{model_sele4} summarize the comparison results of the frequency of the model selection.
In the case of (i--1) (Table \ref{model_sele3}), the probability that a full model is chosen by QBIC seems to be too high when the sample size is small.
This phenomenon in QBIC would be caused by the problem that $|-\p_{\theta}^{2}\mbbh_{n}(\tes)|\approx 0$; we did observe that the values of the determinant in simulations were so small.
Nevertheless, judging from the whole of Table \ref{model_sele3}, tendencies of QBIC, BIC and fAIC for $n\to\infty$ are the same as Section \ref{Simu2.1}.
In (i--2) case (Table \ref{model_sele4}), QBIC tends to perform better than BIC and fAIC for all $n$; indeed, in this case we observed that the values of $|-\p_{\theta}^{2}\mbbh_{n}(\tes)|$ were far from zero.
Moreover, the true model was selected by using QBIC with high probability even for small sample size.

Tables \ref{model_esti3} and \ref{model_esti4} show that a tendency of the estimators for Model 1 and Model 4 is analogous to non-random case.
As is easily expected from the result \cite{UchYos11} concerning parametric estimation of a diffusion with misspecified coefficients, 
we need to let $T_{n}\to\infty$ in order to consistently estimate optimal parameter values.

\begin{table}[t]
\begin{center}
\caption{\footnotesize The number of models selected by QBIC, BIC and fAIC in Section \ref{Simu2.2} (i--1) over 1000 simulations for various $n$ (1-7 express the models, and the true model is model 4)}
\begin{tabular}{| l || c c c c c c c | c c c c c c c | c c c c c c c |} \hline
\multicolumn{1}{| l ||}{Criterion} & \multicolumn{7}{c|}{$n=200$} & \multicolumn{7}{c|}{$n=500$} & \multicolumn{7}{c|}{$n=1000$} \\
& 1 & 2 & 3 & $4^{\ast}$ & 5 & 6 & 7 & 1 & 2 & 3 & $4^{\ast}$ & 5 & 6 & 7 & 1 & 2 & 3 & $4^{\ast}$ & 5 & 6 & 7 \\ \hline
QBIC & 831 & 0 & 5 & {\bf 164} & 0 & 0 & 0 & 657 & 1 & 8 & {\bf 334} & 0 & 0 & 0 & 500 & 0 & 7 & {\bf 493} & 0 & 0 & 0 \\ 
BIC & 8 & 29 & 234 & {\bf 729} & 0 & 0 & 0 & 8 & 5 & 141 & {\bf 846} & 0 & 0 & 0 & 5 & 0 & 117 & {\bf 878} & 0 & 0 & 0 \\ 
fAIC & 75 & 24 & 224 & {\bf 677} & 0 & 0 & 0 & 107 & 4 & 132 & {\bf 757} & 0 & 0 & 0 & 129 & 0 & 105 & {\bf 766} & 0 & 0 & 0 \\ \hline\hline
\multicolumn{1}{| l ||}{Criterion} & \multicolumn{7}{c|}{$n=3000$} & \multicolumn{7}{c|}{$n=5000$} & \multicolumn{7}{c|}{$n=10000$} \\
& 1 & 2 & 3 & $4^{\ast}$ & 5 & 6 & 7 & 1 & 2 & 3 & $4^{\ast}$ & 5 & 6 & 7 & 1 & 2 & 3 & $4^{\ast}$ & 5 & 6 & 7 \\ \hline
QBIC & 250 & 0 & 7 & {\bf 743} & 0 & 0 & 0 & 217 & 0 & 8 & {\bf 775} & 0 & 0 & 0 & 123 & 0 & 3 & {\bf 874} & 0 & 0 & 0 \\ 
BIC & 0 & 0 & 43 & {\bf 957} & 0 & 0 & 0 & 4 & 0 & 40 & {\bf 956} & 0 & 0 & 0 & 4 & 0 & 8 & {\bf 988} & 0 & 0 & 0 \\ 
fAIC & 111 & 0 & 38 & {\bf 851} & 0 & 0 & 0 & 156 & 0 & 30 & {\bf 814} & 0 & 0 & 0 & 153 & 0 & 5 & {\bf 842} & 0 & 0 & 0 \\ \hline
\end{tabular}
\label{model_sele3}
\end{center}
\end{table}

\begin{table}[t]
\begin{center}
\caption{\footnotesize The mean and the standard deviation (s.d.) of the estimator $\hat{\theta}_{1}$, $\hat{\theta}_{2}$ and $\hat{\theta}_{3}$ in Section \ref{Simu2.2} (i--1) for various $n$ (1-7 express the models, and the true parameter $\theta_{0}=(0,-2,3)$)}
\smallskip
\begin{tabular}{| l | l || c c c | c c c | c c c |} \hline
\multicolumn{1}{|l|}{} & \multicolumn{1}{l ||}{} & \multicolumn{3}{c|}{$n=200$} & \multicolumn{3}{c|}{$n=500$} & \multicolumn{3}{c|}{$n=1000$} \\
 & & $\hat{\theta}_{1}$ & $\hat{\theta}_{2}$ & $\hat{\theta}_{3}$ & $\hat{\theta}_{1}$ & $\hat{\theta}_{2}$ & $\hat{\theta}_{3}$ & $\hat{\theta}_{1}$ & $\hat{\theta}_{2}$ & $\hat{\theta}_{3}$ \\ \hline
1 & mean & -0.2910 & -1.7237 & 2.9497 & -0.0683 & -1.9422 & 2.9842 & -0.0458 & -1.9548 & 2.9974 \\ 
 & s.d. & 2.2703 & 2.2300 & 0.8309 & 1.5636 & 1.5293 & 0.5462 & 1.0542 & 1.0335 & 0.3824 \\ \hline
2 & mean & 0.9085 & -2.7435 & -- & 0.8747 & -2.7212 & -- & 0.8685 & -2.7301 & -- \\ 
 & s.d. & 6.4524 & 6.2657 & -- & 6.3784 & 6.1702 & -- & 6.2699 & 6.0750 & -- \\ \hline
3 & mean & -2.0247 & -- & 2.9749 & -2.0092 & -- & 3.0041 & -2.0528 & -- & 2.9834 \\ 
 & s.d. & 0.5204 & -- & 1.2907 & 0.5638 & -- & 1.2907 & 0.3399 & -- & 1.2396 \\ \hline
$4^{\ast}$ & mean & -- & -2.0000 & 2.9724 & -- & -1.9989 & 3.0037 & -- & -1.9986 & 2.9890 \\ 
 & s.d. & -- & 0.1691 & 0.3137 & -- & 0.1075 & 0.2044 & -- & 0.0712 & 0.1384 \\ \hline
5 & mean & -0.2429 & -- & -- & -0.2132 & -- & -- & -0.2409 & -- & -- \\ 
 & s.d. & 0.4806 & -- & -- & 0.4905 & -- & -- & 0.4746 & -- & -- \\ \hline
6 & mean & -- & -1.8835 & -- & -- & -1.8974 & -- & -- & -1.9033 & -- \\ 
 & s.d. & -- & 0.4777 & -- & -- & 0.4701 & -- & -- & 0.4763 & -- \\ \hline
7 & mean & -- & -- & 3.0079 & -- & -- & 2.9923 & -- & -- & 2.9776 \\ 
 & s.d. & -- & -- & 0.5366 & -- & -- & 0.5401 & -- & -- & 0.5385 \\ \hline\hline 
\multicolumn{1}{|c|}{} & \multicolumn{1}{l ||}{} & \multicolumn{3}{c|}{$n=3000$} & \multicolumn{3}{c|}{$n=5000$} & \multicolumn{3}{c|}{$n=10000$} \\
 & & $\hat{\theta}_{1}$ & $\hat{\theta}_{2}$ & $\hat{\theta}_{3}$ & $\hat{\theta}_{1}$ & $\hat{\theta}_{2}$ & $\hat{\theta}_{3}$ & $\hat{\theta}_{1}$ & $\hat{\theta}_{2}$ & $\hat{\theta}_{3}$ \\ \hline
1 & mean & -0.0098 & -1.9915 & 3.0000 & -0.0432 & -1.9569 & 2.9971 & -0.0040 & -1.9957 & 2.9964 \\ 
 & s.d. & 0.5700 & 0.5601 & 0.2043 & 0.4887 & 0.4802 & 0.1722 & 0.3244 & 0.3193 & 0.1145 \\ \hline
2 & mean & 0.9311 & -2.7877 & -- & 0.7949 & -2.6941 & -- & 0.7585 & -2.5841 & -- \\ 
 & s.d. & 6.2734 & 6.0827 & -- & 6.4740 & 6.2692 & -- & 6.4844 & 6.3014 & -- \\ \hline
3 & mean & -2.0158 & -- & 3.0173 & -2.0390 & -- & 2.9579 & -2.0177 & -- & 3.0653 \\ 
 & s.d. & 0.4537 & -- & 1.2141 & 0.4679 & -- & 1.2394 & 0.4491 & -- & 1.2283 \\ \hline
$4^{\ast}$ & mean & -- & -2.0002 & 2.9953 & -- & -1.9986 & 2.9987 & -- & -1.9999 & 2.9994 \\ 
 & s.d. & -- & 0.0433 & 0.0797 & -- & 0.0335 & 0.0617 & -- & 0.0220 & 0.0450 \\ \hline
5 & mean & -0.2391 & -- & -- & -0.2603 & -- & -- & -0.2159 & -- & -- \\ 
 & s.d. & 0.4794 & -- & -- & 0.4765 & -- & -- & 0.4806 & -- & -- \\ \hline
6 & mean & -- & -1.8916 & -- & -- & -1.9019 & -- & -- & -1.8775 & -- \\ 
 & s.d. & -- & 0.4627 & -- & -- & 0.4842 & -- & -- & 0.4521 & -- \\ \hline
7 & mean & -- & -- & 2.9736 & -- & -- & 3.0103 & -- & -- & 2.9793 \\ 
 & s.d. & -- & -- & 0.5533 & -- & -- & 0.5483 & -- & -- & 0.5326 \\ \hline
\end{tabular}
\label{model_esti3}
\end{center}
\end{table}

\begin{table}[t]
\begin{center}
\caption{\footnotesize The number of models selected by QBIC, BIC and fAIC in Section \ref{Simu2.2} (i--2) over 1000 simulations for various $n$ (1-7 express the models, and the true model is model 4)}
\begin{tabular}{| l || c c c c c c c | c c c c c c c | c c c c c c c |} \hline
\multicolumn{1}{| l ||}{Criterion} & \multicolumn{7}{c|}{$n=200$} & \multicolumn{7}{c|}{$n=500$} & \multicolumn{7}{c|}{$n=1000$} \\
& 1 & 2 & 3 & $4^{\ast}$ & 5 & 6 & 7 & 1 & 2 & 3 & $4^{\ast}$ & 5 & 6 & 7 & 1 & 2 & 3 & $4^{\ast}$ & 5 & 6 & 7 \\ \hline
QBIC & 78 & 1 & 1 & {\bf 920} & 0 & 0 & 0 & 38 & 0 & 7 & {\bf 954} & 1 & 0 & 0 & 27 & 1 & 3 & {\bf 969} & 0 & 0 & 0 \\ 
BIC & 6 & 42 & 245 & {\bf 703} & 4 & 0 & 0 & 7 & 5 & 161 & {\bf 826} & 1 & 0 & 0 & 4 & 1 & 122 & {\bf 873} & 0 & 0 & 0 \\ 
fAIC & 74 & 40 & 236 & {\bf 648} & 2 & 0 & 0 & 94 & 2 & 155 & {\bf 748} & 1 & 0 & 0 & 119 & 1 & 115 & {\bf 765} & 0 & 0 & 0 \\ \hline
\end{tabular}
\label{model_sele4}
\end{center}
\end{table}

\begin{table}[t]
\begin{center}
\caption{\footnotesize The mean and the standard deviation (s.d.) of the estimator $\hat{\theta}_{1}$, $\hat{\theta}_{2}$ and $\hat{\theta}_{3}$ in Section \ref{Simu2.2} (i--2) for various $n$ (1-7 express the models, and the true parameter $\theta_{0}=(0,-2,3)$)}
\smallskip
\begin{tabular}{| l | l || c c c | c c c | c c c |} \hline
\multicolumn{1}{|c|}{} & \multicolumn{1}{l ||}{} & \multicolumn{3}{c|}{$n=200$} & \multicolumn{3}{c|}{$n=500$} & \multicolumn{3}{c|}{$n=1000$} \\
 & & $\hat{\theta}_{1}$ & $\hat{\theta}_{2}$ & $\hat{\theta}_{3}$ & $\hat{\theta}_{1}$ & $\hat{\theta}_{2}$ & $\hat{\theta}_{3}$ & $\hat{\theta}_{1}$ & $\hat{\theta}_{2}$ & $\hat{\theta}_{3}$ \\ \hline
1 & mean & -0.0242 & -1.7665 & 2.9373 & -0.0091 & -1.9145 & 2.9985 & -0.0080 & -1.9306 & 2.9845 \\ 
 & s.d. & 0.2290 & 2.2493 & 0.8445 & 0.1366 & 1.3373 & 0.5130 & 0.1028 & 1.0048 & 0.3855 \\ \hline
2 & mean & 0.0620 & -2.5040 & -- & 0.0870 & -2.7296 & -- & 0.1050 & -2.8925 & -- \\ 
 & s.d. & 0.6529 & 6.3498 & -- & 0.6370 & 6.1847 & -- & 0.6329 & 6.1223 & -- \\ \hline
3 & mean & -0.2047 & -- & 2.9181 & -0.2029 & -- & 2.9897 & -0.2048 & -- & 3.0437 \\ 
 & s.d. & 0.0471 & -- & 1.2826 & 0.0473 & -- & 1.2787 & 0.0503 & -- & 1.2950 \\ \hline
$4^{\ast}$ & mean & -- & -1.9948 & 2.9628 & -- & -1.9976 & 2.9963 & -- & -2.0038 & 2.9883 \\ 
 & s.d. & -- & 0.1712 & 0.3275 & -- & 0.1053 & 0.1918 & -- & 0.0751 & 0.1424 \\ \hline
5 & mean & -0.1029 & -- & -- & -0.0946 & -- & -- & -0.0871 & -- & -- \\ 
 & s.d. & 0.1484 & -- & -- & 0.1527 & -- & -- & 0.1518 & -- & -- \\ \hline
6 & mean & -- & -1.9063 & -- & -- & -1.8867 & -- & -- & -1.8722 & -- \\ 
 & s.d. & -- & 0.4801 & -- & -- & 0.4622 & -- & -- & 0.4716 & -- \\ \hline
7 & mean & -- & -- & 3.0147 & -- & -- & 2.9964 & -- & -- & 2.9648 \\ 
 & s.d. & -- & -- & 0.5319 & -- & -- & 0.5530 & -- & -- & 0.5440 \\ \hline
\end{tabular}
\label{model_esti4}
\end{center}
\end{table}


\medskip

We also conducted similar simulations for the case where $X$ is instead given by
\begin{equation}
X_{t}=\left(1,\, \frac{1}{1+B_{t}^{2}},\, \frac{B_{t}}{1+B_{t}^{2}}\right)^{\prime},
\nonumber
\end{equation}
and quite similar tendencies were observed.

\begin{rem}{\rm 
In each case of Section \ref{Simu2}, we have not paid attention to Assumption \ref{Ass9}(ii), which is not so easy to be verified; we refer to \cite{UchYos13} for several general criterion for the non-degeneracy of the statistical random fields in the present context. 
Let us mention almost sure lower bounds of $\det\big(\int_{0}^{1}X_{t}X_{t}^{\prime}dt\big)$ for the models considered in Sections \ref{Simu2.2} (i) and (ii).
Let $X_{1,0}=a$, then, because of the Schwarz inequality
\begin{align*}
\det\left(\int_{0}^{1}X_{t}X_{t}^{\prime}dt\right) 
&=\det\left\{\left(\begin{array}{ccc}
a^{2} & a\int_{0}^{1}X_{2,t}dt & a\int_{0}^{1}X_{3,t}dt \vspace{1mm}\\ 
a\int_{0}^{1}X_{2,t}dt & \int_{0}^{1}X_{2,t}^{2}dt & \int_{0}^{1}X_{2,t}X_{3,t}dt \vspace{1mm}\\ 
a\int_{0}^{1}X_{3,t}dt & \int_{0}^{1}X_{2,t}X_{3,t}dt & \int_{0}^{1}X_{3,t}^{2}dt
\end{array}\right)\right\} \\
&=a^{2}\left[\left\{\int_{0}^{1}\left(X_{2,t}-\int_{0}^{1}X_{2,t}dt\right)^{2}dt\right\}\left\{\int_{0}^{1}\left(X_{3,t}-\int_{0}^{1}X_{3,t}dt\right)^{2}dt\right\}\right. \\
&\qquad\left.-\left(\int_{0}^{1}X_{2,t}X_{3,t}dt-\int_{0}^{1}X_{2,t}dt\int_{0}^{1}X_{3,t}dt\right)^{2}\right] \\
&\geq a^{2}\left[\left\{\int_{0}^{1}\left(X_{2,t}-\int_{0}^{1}X_{2,t}dt\right)\left(X_{3,t}-\int_{0}^{1}X_{3,t}dt\right)dt\right\}^{2}\right. \\
&\qquad\left.-\left(\int_{0}^{1}X_{2,t}X_{3,t}dt-\int_{0}^{1}X_{2,t}dt\int_{0}^{1}X_{3,t}dt\right)^{2}\right].
\end{align*}
Hence $\det\big(\int_{0}^{1}X_{t}X_{t}^{\prime}dt\big)=0$ holds if and only if at least one of the following conditions is satisfied:
\begin{itemize}
\item[(i)] $X_{2,t}-\int_{0}^{1}X_{2,t}dt=0$ for all $t\in[0,1]$;
\item[(ii)] $X_{3,t}-\int_{0}^{1}X_{3,t}dt=0$ for every $t\in[0,1]$; 
\item[(iii)] There exists a constant $c\neq0$, $X_{2,t}-\int_{0}^{1}X_{2,t}dt=c\big(X_{3,t}-\int_{0}^{1}X_{3,t}dt\big)$ for any $t\in[0,1]$.
\end{itemize}
\label{Simu.rem}
}\qed\end{rem}


\subsection{Non-ergodic diffusion process} \label{Simu3}

Let $\mbX_{n}=(X_{t_{j}})_{j=0}^{n}$ be a data set with $t_{j}=j/n$ and the number of data $n$. We simulate 1000 data sets from
\begin{align*}
dX_{t}=\exp\Big\{\frac{5+2X_{t}}{2(1+X_{t}^{2})}\Big\}dw_{t},\quad t\in[0,1],\quad X_{0}=0,
\end{align*}
where $w$ is an $1$-dimensional standard Wiener process.
We consider the following models:
\begin{align*}
&{\bf Model}\;{\bf 1:}\;dX_{t}=\exp\Big\{\frac{\theta_{1}+\theta_{2}X_{t}+\theta_{3}X_{t}^2}{2(1+X_{t}^{2})}\Big\}dw_{t}; 
\;{\bf Model}\;{\bf 2:}\;dX_{t}=\exp\Big\{\frac{\theta_{1}+\theta_{2}X_{t}}{2(1+X_{t}^{2})}\Big\}dw_{t}; \\
&{\bf Model}\;{\bf 3:}\;dX_{t}=\exp\Big\{\frac{\theta_{1}+\theta_{3}X_{t}^2}{2(1+X_{t}^{2})}\Big\}dw_{t};
\;{\bf Model}\;{\bf 4:}\;dX_{t}=\exp\Big\{\frac{\theta_{2}X_{t}+\theta_{3}X_{t}^2}{2(1+X_{t}^{2})}\Big\}dw_{t}; \\
&{\bf Model}\;{\bf 5:}\;dX_{t}=\exp\Big\{\frac{\theta_{1}}{2(1+X_{t}^{2})}\Big\}dw_{t};
\;{\bf Model}\;{\bf 6:}\;dX_{t}=\exp\Big\{\frac{\theta_{2}X_{t}}{2(1+X_{t}^{2})}\Big\}dw_{t}; \\
&{\bf Model}\;{\bf 7:}\;dX_{t}=\exp\Big\{\frac{\theta_{2}X_{t}^{2}}{2(1+X_{t}^{2})}\Big\}dw_{t}.
\end{align*}
Then the optimal model is Model 2, the true parameter being $\tz=(5,2,0)$.
Table \ref{model_sele6} shows that Model 2 is chosen with high probability as the best model for all criteria.
QBIC tends to take values between BIC and fAIC.
Moreover, the larger the sample size becomes, the higher the frequency that the true model is selected by QBIC and BIC become.
In Table \ref{model_esti6}, the estimators exhibit similar tendency to Tables \ref{model_esti3} and \ref{model_esti4}.

\begin{table}[t]
\begin{center}
\caption{\footnotesize The number of models selected by QBIC, BIC and AIC in Section \ref{Simu3} over 1000 simulations for various $n$ (1-7 express the models, and the true model is Model 2)}
\begin{tabular}{ | l || c c c c c c c | c c c c c c c | c c c c c c c |} \hline
\multicolumn{1}{| l ||}{Criterion} & \multicolumn{7}{c|}{$n=200$} & \multicolumn{7}{c|}{$n=500$} & \multicolumn{7}{c|}{$n=1000$} \\
& 1 & $2^{\ast}$ & 3 & 4 & 5 & 6 & 7 & 1 & $2^{\ast}$ & 3 & 4 & 5 & 6 & 7 & 1 & $2^{\ast}$ & 3 & 4 & 5 & 6 & 7 \\ \hline
QBIC & 291 & {\bf 690} & 19 & 0 & 0 & 0 & 0 & 151 & {\bf 832} & 17 & 0 & 0 & 0 & 0 & 115 & {\bf 874} & 11 & 0 & 0 & 0 & 0 \\ 
BIC & 30 & {\bf 733} & 237 & 0 & 0 & 0 & 0 & 15 & {\bf 842} & 143 & 0 & 0 & 0 & 0 & 20 & {\bf 892} & 88 & 0 & 0 & 0 & 0 \\ 
fAIC & 130 & {\bf 642} & 228 & 0 & 0 & 0 & 0 & 135 & {\bf 728} & 137 & 0 & 0 & 0 & 0 & 151 & {\bf 767} & 82 & 0 & 0 & 0 & 0 \\ \hline
\end{tabular}
\label{model_sele6}
\end{center}
\end{table}

\begin{table}[t]
\begin{center}
\caption{\footnotesize The mean and the standard deviation (s.d.) of the estimator $\hat{\theta}_{1}$, $\hat{\theta}_{2}$ and $\hat{\theta}_{3}$ in Section \ref{Simu3} for various $n$ (1-7 express the models, and the true parameter $\tz=(5,2,0)$)}
\smallskip
\begin{tabular}{| l | l || c c c | c c c | c c c |} \hline
\multicolumn{1}{|l|}{} & \multicolumn{1}{l ||}{} & \multicolumn{3}{c|}{$n=200$} & \multicolumn{3}{c|}{$n=500$} & \multicolumn{3}{c|}{$n=1000$} \\
 & & $\hat{\theta}_{1}$ & $\hat{\theta}_{2}$ & $\hat{\theta}_{3}$ & $\hat{\theta}_{1}$ & $\hat{\theta}_{2}$ & $\hat{\theta}_{3}$ & $\hat{\theta}_{1}$ & $\hat{\theta}_{2}$ & $\hat{\theta}_{3}$ \\ \hline
1 & mean & 4.7972 & 1.6082 & 0.0057 & 4.9333 & 1.7426 & -0.0238 & 4.9679 & 1.8711 & -0.0051 \\ 
 & s.d. & 0.8198 & 1.2508 & 0.5725 & 0.5754 & 0.8357 & 0.3637 & 0.4783 & 0.6970 & 0.2835 \\ \hline 
$2^{\ast}$ & mean & 4.9551 & 1.8910 & -- & 5.0031 & 1.9385 & -- & 5.0135 & 1.9713 & -- \\ 
 & s.d. & 0.7156 & 0.4938 & -- & 0.4686 & 0.3138 & -- & 0.3416 & 0.2260 & -- \\ \hline 
3 & mean & 4.7796 & -- & -0.0972 & 4.8136 & -- & -0.1076 & 4.8542 & -- & -0.1369 \\ 
 & s.d. & 1.0889 & -- & 0.7943 & 0.9700 & -- & 0.7772 & 0.9593 & -- & 0.7626 \\ \hline 
4 & mean & -- & -0.3651 & 1.1604 & -- & -0.2268 & 1.4376 & -- & 0.0324 & 1.2284 \\ 
 & s.d. & -- & 4.1514 & 2.8986 & -- & 3.6764 & 2.7531 & -- & 3.8005 & 2.6084 \\ \hline 
5 & mean & 4.9133 & -- & -- & 4.9294 & -- & -- & 4.9231 & -- & -- \\ 
 & s.d. & 0.5360 & -- & -- & 0.5259 & -- & -- & 0.5477 & -- & -- \\ \hline 
6 & mean & -- & 1.6946 & -- & -- & 1.7188 & -- & -- & 1.7393 & -- \\ 
 & s.d. & -- & 0.4614 & -- & -- & 0.4670 & -- & -- & 0.4742 & -- \\ \hline 
7 & mean & -- & -- & 0.4908 & -- & -- & 0.4826 & -- & -- & 0.4926 \\ 
 & s.d. & -- & -- & 0.2881 & -- & -- & 0.2782 & -- & -- & 0.2872 \\ \hline 
\end{tabular}
\label{model_esti6}
\end{center}
\end{table}

\begin{rem}{\rm
We write $Z_{t}=\big(\frac{1}{1+X_{t}^2},\frac{X_{t}}{1+X_{t}^2},\frac{X_{t}^{2}}{1+X_{t}^2}\big)^{\prime}$. In a similar way to Remark \ref{Simu.rem}, we can see that
\begin{align*}
&\det\left(\int_{0}^{T}Z_{t}Z_{t}^{\prime}dt\right) \\
&=\det\left\{\left(\begin{array}{ccc}
\int_{0}^{1}\frac{1}{(1+X_{t}^2)^{2}}dt & \int_{0}^{1}\frac{X_{t}}{(1+X_{t}^2)^{2}}dt & \int_{0}^{1}\frac{X_{t}^{2}}{(1+X_{t}^2)^{2}}dt \vspace{0.5mm}\\ 
\int_{0}^{1}\frac{X_{t}}{(1+X_{t}^2)^{2}}dt & \int_{0}^{1}\frac{X_{t}^{2}}{(1+X_{t}^2)^{2}}dt & \int_{0}^{1}\frac{X_{t}^{3}}{(1+X_{t}^2)^{2}}dt \vspace{0.5mm}\\ 
\int_{0}^{1}\frac{X_{t}}{(1+X_{t}^2)^{2}}dt & \int_{0}^{1}\frac{X_{t}^{3}}{(1+X_{t}^2)^{2}}dt & \int_{0}^{1}\frac{X_{t}^{4}}{(1+X_{t}^2)^{2}}dt
\end{array}\right)\right\} \\
&=\left(\int_{0}^{1}\frac{1}{(1+X_{t}^2)^{2}}dt\right)\left[\left\{\int_{0}^{1}\frac{1}{(1+X_{t}^2)^{2}}\left(X_{t}-\left(\int_{0}^{1}\frac{1}{(1+X_{t}^2)^{2}}dt\right)^{-1}\int_{0}^{1}\frac{X_{t}}{(1+X_{t}^2)^{2}}dt\right)^{2}dt\right\}\right. \\
&\qquad\times \left\{\int_{0}^{1}\frac{1}{(1+X_{t}^2)^{2}}\left(X_{t}^{2}-\left(\int_{0}^{1}\frac{1}{(1+X_{t}^2)^{2}}dt\right)^{-1}\int_{0}^{1}\frac{X_{t}^{2}}{(1+X_{t}^2)^{2}}dt\right)^{2}dt\right\} \\
&\qquad\left.-\left(\int_{0}^{1}\frac{X_{t}^{3}}{(1+X_{t}^2)^{2}}dt-\left(\int_{0}^{1}\frac{1}{(1+X_{t}^2)^{2}}dt\right)^{-1}\int_{0}^{1}\frac{X_{t}}{(1+X_{t}^2)^{2}}dt\int_{0}^{1}\frac{X_{t}^{2}}{(1+X_{t}^2)^{2}}dt\right)^{2}\right] \\
&\geq0,
\end{align*}
with the last equality holding if and only if at least one of the following holds true for any $t\in[0,1]$:
\begin{itemize}
\item[(i)] $X_{t}-\big(\int_{0}^{1}\frac{1}{(1+X_{t}^2)^{2}}dt\big)^{-1}$ $\int_{0}^{1}\frac{X_{t}^{2}}{(1+X_{t}^2)^{2}}dt=0$;
\item[(ii)] $X_{t}^{2}-\big(\int_{0}^{1}\frac{1}{(1+X_{t}^2)^{2}}dt\big)^{-1}\int_{0}^{1}\frac{X_{t}^{2}}{(1+X_{t}^2)^{2}}dt=0$;
\item[(iii)] There exists a constant $c\neq0$ such that
\begin{equation}
X_{t}-\bigg(\int_{0}^{1}\frac{1}{(1+X_{t}^2)^{2}}dt\bigg)^{-1}\int_{0}^{1}\frac{X_{t}^{2}}{(1+X_{t}^2)^{2}}dt
=c\bigg\{X_{t}^{2}-\bigg(\int_{0}^{1}\frac{1}{(1+X_{t}^2)^{2}}dt\bigg)^{-1}\int_{0}^{1}\frac{X_{t}^{2}}{(1+X_{t}^2)^{2}}dt\bigg\}.
\nonumber
\end{equation}
\end{itemize}
}\qed\end{rem}

\section{Proofs} \label{hm:sec_Proofs}

Recall that $\mbbu_{n}(\tz)=\{u\in\mbbr^{p};\tz+A_{n}(\tz)u\in\Theta\}$. 
In what follows, we deal with the zero-extended version of $\mbbz_{n}$ and use the same notation: 
$\mbbz_{n}$ vanishes outside $\mbbu_{n}(\tz)$, so that
\begin{align*}
\int_{\mbbr^{p}\backslash\mbbu_{n}(\tz)}\mbbz_{n}(u)du=0.
\end{align*}

\subsection{Proof of Theorem \ref{QlfTh3} (ii)}
\label{hm:sec_proof.sc.tpe}

By using the Taylor expansion, we obtain that
\begin{align*}
\mbbz_{n}(u)=\exp\left(\D_{n}[u]-\frac{1}{2}\Gam_{n}(\tet)[u,u]\right)
\end{align*}
for a random point $\tet$ on the segment connecting $\tz$ and $\tz+A_{n}(\tz)u$. 
Then, for any positive $\epsilon$, $\del$ and $M$, we have
\begin{align*}
&\pr\left(\int_{\mbbu_{n}(\tz)\cap\{|u|\geq M\}}\mbbz_{n}(u)du>\epsilon\right) \\
&\leq\pr\left(\int_{|u|\geq M}\mbbz_{n}(u)du>\epsilon\,;
~\inf_{\theta\in\Theta}\lambda_{\min}\big(\Gam_{n}(\theta)\big)<\delta\right)+\pr\left(\int_{|u|\geq M}\mbbz_{n}(u)du>\epsilon\,;
~\inf_{\theta\in\Theta}\lambda_{\min}\big(\Gam_{n}(\theta)\big)\geq\delta\right) \\
&\leq\pr\left(\inf_{\theta\in\Theta}\lambda_{\min}\big(\Gam_{n}(\theta)\big)<\delta\right)+\pr\left\{\int_{|u|\geq M}\exp\left(\D_{n}[u]-\frac{\delta}{2}[u,u]\right)du>\epsilon\right\}.
\end{align*}
Under (\ref{AssTh3}), we can find $\del$ and $N'$ for which 
the first term in the rightmost side can be bounded by $\ep/2$ for $n\ge N'$. 
Given such a $\del$, making use of the tightness of $(\D_{n})$ we can take $N''$ and $M>0$ large enough to ensure
\begin{align*}
\sup_{n\ge N''}\pr\left\{\int_{|u|\geq M}\exp\left(\D_{n}[u]-\frac{\delta}{2}[u,u]\right)du>\epsilon\right\}<\frac{\epsilon}{2}.
\end{align*}
Hence we have $\displaystyle \sup_{n\ge N}\pr\left(\int_{\mbbu_{n}(\tz)\cap\{|u|\geq M\}}\mbbz_{n}(u)du>\epsilon\right)<\epsilon$ 
for $N:=N'\vee N''$, completing the proof.


\subsection{Proof of Theorem \ref{QlfTh2}}

(i) By the change of variable $\theta=\tz+A_{n}(\tz)u$, the marginal quasi-log likelihood function equals
\begin{align*}
\mbbh_{n}(\tz)+\sumk p_{k}\log a_{k,n}(\theta_{k,0})+\log\bigg(\int_{\mbbu_{n}(\tz)}\mbbz_{n}(u)\pi_{n}\big(\tz+A_{n}(\tz)u\big)du\bigg).
\end{align*}
Consequently,
\begin{align*}
\log\bigg(\int_{\Theta}\exp\{\mbbh_{n}(\theta)\}\pi_{n}(\theta)d\theta\bigg)-\bigg(\mbbh_{n}(\tz)+\sumk p_{k}\log a_{k,n}(\theta_{k,0})+\log\bar{Q}_{n}\bigg)=\log\big(\bar{Q}_{n}+\bar{\epsilon}_{n}\big)-\log\bar{Q}_{n},
\end{align*}
where
\begin{align*}
&\bar{Q}_{n}=\pi_{n}(\tz)\int_{\mbbr^{p}}\exp\Big(\D_{n}[u]-\frac{1}{2}\Gam_{0}[u,u]\Big)du,\\
&\bar{\epsilon}_{n}=\int_{\mbbu_{n}(\tz)}\mbbz_{n}(u)\big(\pi_{n}(\tz+A_{n}(\theta_{0})u)-\pi_{n}(\tz)\big)du
+\pi_{n}(\tz)\int_{\mbbr^{p}}\left\{\mbbz_{n}(u)-\exp\left(\D_{n}[u]-\frac{1}{2}\Gam_{0}[u,u]\right)\right\}du.
\end{align*}
We will show that $|\log\big(\bar{Q}_{n}+\bar{\epsilon}_{n}\big)-\log\bar{Q}_{n}|\cip 0$.

First, we note that
\begin{align*}
\bar{Q}_{n}
&=\pi_{n}(\tz)\exp\left(\frac{1}{2}\|\Gam_{0}^{-\frac{1}{2}}\D_{n}\|^{2}\right)\int_{\mbbr^{p}}\exp\left(-\frac{1}{2}\Gam_{0}[u-\Gam_{0}^{-1}\D_{n},u-\Gam_{0}^{-1}\D_{n}]\right)du \\
&=\pi_{n}(\tz)\exp\left(\frac{1}{2}\|\Gam_{0}^{-\frac{1}{2}}\D_{n}\|^{2}\right)(2\pi)^{\frac{p}{2}}|\Gam_{0}|^{-\frac{1}{2}}.
\end{align*}
Because of Assumptions \ref{Ass1} and \ref{Ass2}, $\log\bar{Q}_{n}$ is given by
\begin{align*}
\log\bar{Q}_{n}&=\log\pi_{n}(\tz)+\frac{1}{2}\|\Gam_{0}^{-\frac{1}{2}}\D_{n}\|^{2}+\frac{p}{2}\log(2\pi)-\frac{1}{2}\log|\Gam_{0}|.
\end{align*}
Next we observe that
\begin{align*}
|\bar{\epsilon}_{n}|&\leq\int_{\mbbu_{n}(\tz)}\mbbz_{n}(u)\big|\pi_{n}\big(\tz+A_{n}(\theta_{0})u\big)-\pi_{n}(\tz)\big|du+\pi_{n}(\tz)\int_{\mbbr^{p}}\left|\mbbz_{n}(u)-\exp\left(\D_{n}[u]-\frac{1}{2}\Gam_{0}[u,u]\right)\right|du.
\end{align*}
Fix any $\epsilon>0$. Then, for each $M>0$ we have
\begin{align}
&\pr\left(\int_{\mbbu_{n}(\tz)}\mbbz_{n}(u)\big|\pi_{n}\big(\tz+A_{n}(\tz)u\big)-\pi_{n}(\tz)\big|du>\epsilon\right) \nn\\
&\leq\pr\left(\sup_{|u|<M}\big|\pi_{n}\big(\tz+A_{n}(\tz)u\big)-\pi_{n}(\tz)\big|\sup_{|u|<M}\mbbz_{n}(u)>\frac{\epsilon}{2}\right)
\nn\\
&{}\qquad
+\pr\left(2\sup_{\theta}\pi_{n}(\theta)\int_{|u|\geq M}\mbbz_{n}(u)du>\frac{\epsilon}{2}\right).
\label{hm:add.eq-1}
\end{align}
Let $r_{n}(u):=\frac{1}{2}\left(\Gam_{0}-\Gam_{n}\right)[u,u]+\frac{1}{6}\sum_{i,j,k = 1}^p\big(\p_{\theta^{i}}\p_{\theta^{j}}\p_{\theta^{k}}\mbbh_{n}(\tet) \big)A_{n,ii}(\tz)A_{n,jj}(\tz)A_{n,kk}(\tz)u_{i}u_{j}u_{k}$ for a point $\tet$ between $\tz$ and $\tz+A_{n}(\tz)u$.
Then, under the assumptions we have $\sup_{|u|<K_{0}}|r_{n}(u)|\cip 0$ for every $K_{0}>0$. 
Further, $\sup_{|u|<M}\mbbz_{n}(u)=\sup_{|u|<M}\exp\Big(\D_{n}[u]-\frac{1}{2}\Gam_{0}[u,u]+r_{n}(u)\Big)=O_{p}(1)$, 
so that $\sup_{|u|<M}\big|\pi_{n}\big(\tz+A_{n}(\tz)u\big)-\pi_{n}(\tz)\big|\sup_{|u|<M}\mbbz_{n}(u)=o_{p}(1)$ for each $M>0$. 
Under Assumption \ref{Ass3} we can take a sufficiently large $M$ to conclude that
\begin{align*}
\pr\left(\int_{\mbbu_{n}(\tz)}\mbbz_{n}(u)\big|\pi_{n}\big(\tz+A_{n}(\tz)u\big)-\pi_{n}(\tz)\big|du>\epsilon\right)
<\frac{\ep}{2}+\frac{\ep}{4\sup_{\theta}\pi_{n}(\theta)}.
\end{align*}
It follows that $\int_{\mbbu_{n}(\tz)}\mbbz_{n}(u)\big|\pi_{n}\big(\tz+A_{n}(\tz)u\big)-\pi_{n}(\tz)\big|du\cip0$.
Next, for any $\delta>0$ and $K>0$,
\begin{align*}
&\pr\left\{\int_{\mbbr^{p}}\left|\mbbz_{n}(u)-\exp\left(\D_{n}[u]-\frac{1}{2}\Gam_{0}[u,u]\right)\right|du>\delta\right\} \\
&\le \pr\left\{\int_{|u|<K}\left|\mbbz_{n}(u)-\exp\left(\D_{n}[u]-\frac{1}{2}\Gam_{0}[u,u]\right)\right|du>\frac{\delta}{2}\right\}+\pr\left(\int_{|u|\geq K}\mbbz_{n}(u)du>\frac{\delta}{4}\right) \\
&\qquad+\pr\left\{\int_{|u|\geq K}\exp\left(\D_{n}[u]-\frac{1}{2}\Gam_{0}[u,u]\right)du>\frac{\delta}{4}\right\}.
\end{align*}
We can pick $K>0$ and $N''$ large enough to ensure
\begin{align}
&\sup_{n\ge N''}\pr\left(\int_{|u|\geq K}\mbbz_{n}(u)du>\frac{\delta}{4}\right)<\frac{\delta}{4}, \label{ZnIne}\\
&\sup_{n\ge N''}\pr\left\{\int_{|u|\geq K}\exp\left(\D_{n}[u]-\frac{1}{2}\Gam_{0}[u,u]\right)du>\frac{\delta}{4}\right\}<\frac{\delta}{4}. \label{NorIne}
\end{align}
Since $\D_{n}[u]-\frac{1}{2}\Gam_{0}[u,u]\leq \frac{1}{2}\D_{n}^{\prime}\Gam_{0}\D_{n}$ if and only if $u=\Gam_{0}^{-1}\D_n$, for the same $K>0$ as above we get
\begin{align}
\int_{|u|<K}\left|\mbbz_{n}(u)-\exp\left(\D_{n}[u]-\frac{1}{2}\Gam_{0}[u,u]\right)\right|du&\lesssim\sup_{|u|<K}\left|\exp\left(\D_{n}[u]-\frac{1}{2}\Gam_{0}[u,u]\right)\left\{\exp\big(r_{n}(u)\big)-1\right\}\right| \nonumber\\
&\leq\sup_{|u|<K}\left|\exp\big(r_{n}(u)\big)-1\right|\exp\left(\frac{1}{2}\D_{n}^{\prime}\Gam_{0}\D_{n}\right) \nonumber\\
&\cip0. \label{Conv}
\end{align}
Because of (\ref{ZnIne}) to (\ref{Conv}) we have $\int_{\mbbr^{p}}\left|\mbbz_{n}(u)-\exp\left(\D_{n}[u]-\frac{1}{2}\Gam_{0}[u,u]\right)\right|du \cip 0$, hence $\bar{\epsilon}_{n}\cip0$ and it follows that
\begin{align*}
\log\big(\bar{Q}_{n}+\bar{\epsilon}_{n}\big)-\log\bar{Q}_{n}=\big(\log\bar{Q}_{n}+o_{p}(1)\big)-\log\bar{Q}_{n}=o_{p}(1),
\end{align*}
establishing the claim (i).

\medskip

(ii) By the consistency of $\tes$ we may focus on the event $\{\tes\in\Theta\}\,(\subset\{\p_{\theta}\mbbh_{n}(\tes)=0\})$. Then
\begin{align*}
\D_{n}&=-A_{n}(\tz)\int_{0}^{1}\p_{\theta}^{2}\mbbh_{n}\big(\tes+s(\tz-\tes)\big)dsA_{n}(\tz) [\hat{u}_{n}]=\{\Gam_{0}+o_{p}(1)\} [\hat{u}_{n}],
\end{align*}
so that $\hat{u}_{n}=\Gam_{0}^{-1}\D_{n} + o_{p}(1)$. 
Therefore
\begin{align*}
\mbbh_{n}(\tz)&=\mbbh_{n}(\tes)-\frac{1}{2}\hat{u}_{n}^{\prime}\Gam_{0}\hat{u}_{n}+o_{p}(1) \\
&=\mbbh_{n}(\tes)-\frac{1}{2}(\Gam_{0}^{-1}\D_{n})^{\prime}\Gam_{0}(\Gam_{0}^{-1}\D_{n})+o_{p}(1) \\
&=\mbbh_{n}(\tes)-\frac{1}{2}\|\Gam_{0}^{-\frac{1}{2}}\D_{n}\|^{2}+o_{p}(1),
\end{align*}
which combined with the preceding result (i) and the fact $-A_{n}(\tes)\p_{\theta}^{2}\mbbh_{n}(\tes)A_{n}(\tes)=\Gam_{0}+o_{p}(1)$ 
establishes (ii).

\begin{rem}{\rm 
Recall that Bayes point estimator associated with a loss function 
$\mathfrak{L}:\Theta\times\Theta\to\mbbr$ is defined to be any statistics $\tilde{\theta}_{n}(\mathfrak{L})$ minimizing the random function
\begin{equation}
t\mapsto\int_{\Theta}\mathfrak{L}(t,\theta)\pi_{n}(\theta|\mbX_{n})d\theta,
\nonumber
\end{equation}
where
\begin{equation}
\pi_{n}(\theta|\mbX_{n}):=\frac{\exp\{\mbbh_{n}(\theta)\}\pi_{n}(\theta)}{\int_{\Theta}\exp\{\mbbh_{n}(\theta)\}\pi_{n}(\theta)d\theta}
\nonumber
\end{equation}
denotes the formal posterior density of $\theta$; in particular, the quadratic loss 
$\mathfrak{L}_{2}(t,\theta):=|t-\theta|^{2}$ gives rise to the posterior-mean:
\begin{equation}
\tilde{\theta}_{n}(\mathfrak{L}_{2}):=\frac{\int_{\Theta}\theta\exp\{\mbbh_{n}(\theta)\}\pi_{n}(\theta)d\theta}
{\int_{\Theta}\exp\{\mbbh_{n}(\theta)\}\pi_{n}(\theta)d\theta}
\nonumber
\end{equation}
In the theoretical derivation of the QBIC, we made use of the fact that (at least with sufficiently high probability) 
$\p_{\theta}\mbbh_{n}(\tes)=0$, which does not hold true if we use an integral-type Bayes point estimator.
}\qed\end{rem}

\subsection{Proof of Theorem \ref{hm:thm_mc}}

Let
\begin{align}
F_{n} &:= -2\log\left(\int_{\Theta}\exp\{\mbbh_{n}(\theta)\}\pi_{n}(\theta)d\theta\right), \nn\\
F'_{n} &:= -2\mbbh_{n}(\tz) -2\sumk p_k\log a_{k,n}(\theta_{k,0}) +\log|\Gam_{0}| -p\log2\pi -\|\Gam_{0}^{-\frac{1}{2}}\D_{n}\|^2 -2\log\pi_{n}(\tz).
\nn
\end{align}
We complete the proof by showing that $\E(|F_{n}-F'_{n}|)\to 0$ and $\E(|F'_{n}-\qbic^{\sharp}|)\to 0$ separately. 
Below we will write ``$a_{n}\lesssim b_{n}$'' if the inequality ``$a_{n}\le Cb_{n}$'' for a universal positive constant $C$ holds true.

\medskip

\textit{Proof of  $\E(|F_{n}-F'_{n}|)\to 0$.}  
Obviously Assumption \ref{hm:A_mc1} implies Assumption \ref{Ass1}, hence Theorem \ref{QlfTh2} yields that $F_{n}=F'_{n}+o_{p}(1)$. 

Now pick any $\kappa\in(1, q/p)$ with $q$ being the constant given in Assumption \ref{hm:A_mc3}. To deduce the claim it suffices to show that $\limsup_{n}\E(|F_{n}-F'_{n}|^{\kappa})<\infty$.
Then for some $\delta\in(0,1/\kappa)$,
\begin{align*}
\E(|F_{n}-F'_{n}|^{\kappa}) &\lesssim 1 + \E\bigg(\bigg| \log\bigg(\int_{\mbbu_{n}(\tz)}\mbbz_{n}(u)du\bigg) \bigg|^{\kappa}\bigg) 
+ \E\big(\big|\log|\Gam_{0}|\big|^{\kappa}\big) + \E\left(\left|\Gam_{0}^{-1}[\D_{n}^{\otimes 2}]\right|^{\kappa}\right) \\
&\lesssim 1 + \mbbe\left[\left\{-\log\left(\int_{\mbbu_{n}(\tz)}\mbbz_{n}(u)du\right)\right\}^{\kappa};\int_{\mbbu_{n}(\tz)}\mbbz_{n}(u)du\leq1\right] \\
&\qquad + \mbbe\left[\left\{\log\left(\int_{\mbbu_{n}(\tz)}\mbbz_{n}(u)du\right)\right\}^{\kappa};\int_{\mbbu_{n}(\tz)}\mbbz_{n}(u)du>1\right] \\
&\qquad + \E\big(\lam_{\min}^{-p\kappa}\left(\Gam_{0}\right) + |\Gam_{0}|^{\kappa} \big)
+ \E\left( |\D_{n}|^{2\kappa}\lam_{\min}^{-p\kappa}\left(\Gam_{0}\right) \right) \\
&\lesssim 1 + \mbbe\left\{\left(\int_{\mbbu_{n}(\tz)}\mbbz_{n}(u)du\right)^{-\delta\kappa};\int_{\mbbu_{n}(\tz)}\mbbz_{n}(u)du\leq1\right\} \\
&\qquad + \mbbe\left[\left\{\log\left(\int_{\mbbu_{n}(\tz)}\exp\left(\D_{n}[u]-\frac{1}{2}\inf_{\theta\in\Theta}\lambda_{\min}\big(\Gam_{n}(\theta)\big)[u,u]\right)du\right)\right\}^{\kappa};\int_{\mbbu_{n}(\tz)}\mbbz_{n}(u)du>1\right] \\
&\lesssim 1 + \mbbe\left\{\left(\int_{\mbbu_{n}(\tz)}\mbbz_{n}(u)du\right)^{-1}\right\} \\
&\qquad + \mbbe\left[\left\{\left(\inf_{\theta\in\Theta}\lambda_{\min}\big(\Gam_{n}(\theta)\big)\right)^{-1}[\D_{n}^{\otimes 2}] \right.\right. \\
&\qquad\qquad\left.\left.+\int_{\mbbr^{p}}\exp\left(-\frac{1}{2}\lambda_{\min}\big(\Gam_{n}(\theta)\big)\left[\left(u-\left(\inf_{\theta\in\Theta}\lambda_{\min}\big(\Gam_{n}(\theta)\big)\right)^{-1}\D_{n}\right)^{\otimes2}\right]\right)du\right\}^{\kappa}\right] \\
&\lesssim 1 + \mbbe\left[\left\{\left(\inf_{\theta\in\Theta}\lambda_{\min}\big(\Gam_{n}(\theta)\big)\right)^{-1}[\D_{n}^{\otimes2}]\right\}^{\kappa}+\left|\log\left(\inf_{\theta\in\Theta}\lambda_{\min}\big(\Gam_{n}(\theta)\big)\right)\right|^{\kappa}\right] \\
&\leq 1 + \mbbe\left[\left(\inf_{\theta\in\Theta}\lambda_{\min}\big(\Gam_{n}(\theta)\big)\right)^{-\kappa}|\D_{n}|^{2\kappa}+\left(\inf_{\theta\in\Theta}\lambda_{\min}\big(\Gam_{n}(\theta)\big)\right)^{\kappa}+\left(\inf_{\theta\in\Theta}\lambda_{\min}\big(\Gam_{n}(\theta)\big)\right)^{-\kappa}\right] \\
&\leq 1 + \mbbe\left[\left(\sup_{\theta\in\Theta}\lambda_{\min}^{-\kappa}\big(\Gam_{n}(\theta)\big)\right)(|\D_{n}|^{2\kappa}+1)+\sup_{\theta\in\Theta}|\Gam_{n}(\theta)|^{\kappa}\big)\right] \\
&\lesssim 1,
\end{align*}
where: in the fifth step, we applied \cite[Lemma 2]{Yos11} for the second terms; in the sixth step, we made use of the inequality: for any $s>0$, $|\log x| \lesssim x^{-s} + x^{s}$ for $x>0$.

\medskip

\textit{Proof of $\E(|F'_{n}-\qbic^{\sharp}|)\to 0$.} 
We have $|F'_{n}-\qbic^{\sharp}|\lesssim \overline{R}_{1,n}+\overline{R}_{2,n}+\overline{R}_{3,n}$, where
\begin{align}
\overline{R}_{1,n} &:= \bigg|\log\frac{\pi_{n}(\tes)}{\pi_{n}(\tz)}\bigg| + \sup_{\theta}\big|\p_{\theta}^{3}\mbbh_{n}(\theta)[(\tes-\tz)^{\otimes 3}] \big| 
+\big| \p_{\theta}\mbbh_{n}(\tes) [\tes-\tz] \big|, \nn\\
\overline{R}_{2,n} &:= \left| \log\big| \Gam_{n}(\tes)\big| - \log|\Gam_{0}| \right|, \nn\\
\overline{R}_{3,n} &:= \left| \Gam_{n}(\tes)[\hat{u}_{n}^{\otimes 2}] - \Gam_{0}^{-1}[\D_{n}^{\otimes 2}] \right|.
\nonumber
\end{align}
We will derive $\E(\overline{R}_{i,n})\to 0$ for $i=1,2,3$ separately.

First we look at $\overline{R}_{1,n}$. The convergence  $\pi_{n}(\tes)/\pi_{n}(\tz)\cip 1$ holds under Assumption \ref{hm:A_mc2}: 
indeed, for any $\ep>0$ and $M>0$ we have 
$\pr(|\pi_{n}(\tes)/\pi_{n}(\tz) - 1| > \ep)\le \sup_{n}\pr(|\hat{u}_{n}|>M) + \pr\{(|\hat{u}_{n}|\le M)\cap(\sup_{|u|\le M}|\pi_{n}(\tz+A_{n}(\tz)u) -\pi_{n}(\tz)|\ge C\ep)\}$ 
for some constant $C>0$, from which under Assumption \ref{Ass2} the claim follows on letting $M$ large enough and then $n\to\infty$. 
Then
\begin{equation}
\lim_{n}\E\bigg(\bigg|\log\frac{\pi_{n}(\tes)}{\pi_{n}(\tz)}\bigg|\bigg) \to 0
\nn
\end{equation}
by the bounded convergence theorem. 
We are assuming that $(\hat{u}_{n})_{n}$ is $L^{r}(\pr)$-bounded for some $r>3$ (Assumption \ref{hm:A_mc5}), hence under the assumptions we can apply H\"older's inequality to deduce
\begin{align}
\E\bigg(\sup_{\theta}\big|\p_{\theta}^{3}\mbbh_{n}(\theta)[(\tes-\tz)^{\otimes 3}] \big|\bigg)
&\lesssim \bigg(\max_{i\le p}A_{n,ii}(\tz)\bigg)\E\bigg(\sup_{\theta}\big|A_{n}(\tz)\p_{\theta}^{3}\mbbh_{n}(\theta)A_{n}(\tz)\big| |\hat{u}_{n}|^{3}\bigg)
\nn\\
&\lesssim o(1) \|\hat{u}_{n}\|_{r}\to 0.
\nonumber
\end{align}
Also, for any $s>1$ small enough we have
\begin{align}
\E\left(\big|\p_{\theta}\mbbh_{n}(\tes)[\tes-\tz]\big|^{s}\right) \lesssim \E\left(|\D_{n}|^{s}|\hat{u}_{n}|^{s}\right) 
+ \E\bigg(\sup_{\theta}|\Gam_{n}(\theta)|^{s}|\hat{u}_{n}|^{2s}\bigg) \lesssim \|\hat{u}_{n}\|_{r}\lesssim 1,
\nonumber
\end{align}
thereby $\E(|\p_{\theta}\mbbh_{n}(\tes)[\tes-\tz]|)\to 0$, concluding that $\E(\overline{R}_{1,n})\to 0$.

For handling $\E(\overline{R}_{2,n})$, it suffices to observe that $\overline{R}_{2,n}=|\log(|\Gam_{n}(\tes)|/|\Gam_{0}|)|\cip 0$ and that for any $s'>1$ and $s\in(0,q/p)$,
\begin{align}
\E\left(\overline{R}_{2,n}^{s'}\right)
&\lesssim \E\bigg\{\lam_{\min}^{-sp}(\Gam_{0})\bigg(\sup_{\theta}|\Gam_{n}(\theta)|^{s}\bigg) 
+ |\Gam_{0}|^{s}\bigg(\sup_{\theta}\lam_{\min}^{-sp}(\Gam_{n}(\theta))\bigg)\bigg\} \nn\\
&\lesssim \E\left\{\lam_{\min}^{-q}(\Gam_{0})\right\}+\E\left(\sup_{\theta}\lam_{\min}^{-q}(\Gam_{n}(\theta))\right) \lesssim 1.
\nonumber
\end{align}

To deduce $\E(\overline{R}_{3,n})\to 0$, we note that $\overline{R}_{3,n}=|\{\Gam_{0}+o_{p}(1)\}[\hat{u}_{n}^{\otimes 2}] - \Gam_{0}^{-1}[\{\Gam_{0}[\hat{u}_{n}]\}^{\otimes 2}]| =o_{p}(1)$, 
since $\hat{u}_{n}=\Gam_{0}^{-1}\D_{n} + o_{p}(1)$ as was mentioned in the proof of Theorem \ref{QlfTh2}(ii). The uniform integrability of $(\overline{R}_{3,n})_{n}$ can be verified in a similar manner to the previous case. The proof is complete.

\subsection{Proof of Theorem \ref{ErTh1}}

Under Assumptions \ref{Ass4} to \ref{Ass7}, the argument in \cite[Section 6]{Yos11} ensures the PLDI: for every $L>0$ we can find a constant $C_{L}>0$ such that
\begin{align*}
\pr\left(\sup_{(u_{1},\theta_{2})\in\{r\leq|u_{1}|\}\times\Theta_{2}}\mbbz_{n}^{1}(u_{1};\theta_{1,0},\theta_{2})\geq e^{-r}\right)
+\pr\left(\sup_{u_{2}\in\{r\leq|u_{2}|\}}\mbbz_{n}^{2}(u_{2};\theta_{1,0},\theta_{2,0})\geq e^{-r}\right)\leq\frac{C_{L}}{r^{L}}
\end{align*}
for any $n>0$ and $r>0$. This implies that the inequality (\ref{pldi.2}) holds (see Remark \ref{hm:rem_pldi}).
Assumption \ref{Ass1} readily follows by making use of the lemmas in \cite[Section 6]{Yos11}, and we omit them (see \cite[Section 5.3]{EguMas15} for some details).

\subsection{Proof of Theorem \ref{VolaCor1}}

It is enough to check the conditions [H1] and [H2] of \cite{UchYos13}.

The condition [H1] is a regularity conditions concerning the processes $X$ and $b$, and the non-degeneracy of the diffusion-coefficient function $S(x,\theta)$. 
As a consequence of Assumption \ref{Ass9}(i) and the compactness of $\Theta$, we get
\begin{equation}
\inf_{\omega\in\Omega,t\leq T,\theta\in\Theta}\exp(X_{t}^{\prime}\theta) > 0.
\nonumber
\end{equation}
Based on this inequality and Assumption \ref{Ass8}, it is straightforward to verify [H1].

The condition [H2] is the non-degeneracy of the random field in the limit: 
for every $L>0$, there exists $C_{L}>0$ such that
\begin{equation}
\pr\left(\chi_{0}\leq r^{-1}\right) \leq\frac{C_{L}}{r^{L}},\qquad r>0,
\nonumber
\end{equation}
where
\begin{align*}
\chi_{0}=\inf_{\theta\neq\tz}\frac{1}{2T|\theta-\tz|^{2}}\int_{0}^{T}\Big\{X_{t}^{\prime}(\theta-\tz)+\Big(\exp\big(X_{t}^{\prime}(\tz-\theta)\big)-1\Big)\Big\}dt.
\end{align*}
Since $\exp(x)=1+x+\frac{1}{2}\exp(\xi x)x^{2}$ for some $\xi$ satisfying $0<\xi<1$, letting $x=X_{t}^{\prime}(\theta-\tz)$ we obtain
\begin{align*}
X_{t}^{\prime}(\tz-\theta)+\left\{\exp\big(X_{t}^{\prime}(\tz-\theta)\big)-1\right\}&=\exp\big(X_{t}^{\prime}(\tz-\theta)\big)-1-X_{t}^{\prime}(\tz-\theta) \\
&=\frac{1}{2}\exp\big(\xi X_{t}^{\prime}(\tz-\theta)\big)\big(X_{t}^{\prime}(\tz-\theta)\big)^{2} \\
&=\frac{1}{2}\exp\big(\xi X_{t}^{\prime}(\tz-\theta)\big)(\tz-\theta)^{\prime}X_{t}X_{t}^{\prime}(\tz-\theta) \\
&\geq\frac{1}{2}\exp(-C_{0})(\tz-\theta)^{\prime}X_{t}X_{t}^{\prime}(\tz-\theta)
\end{align*}
for some $C_{0}>0$. Hence
\begin{align*}
\chi_{0}\geq\frac{\exp(-C_{0})}{4T}\inf_{\theta\neq\tz}\frac{\int_{0}^{T}(\tz-\theta)^{\prime}X_{t}X_{t}^{\prime}(\tz-\theta)dt}{|\theta-\tz|^{2}}
\gtrsim\lambda_{\min}\bigg(\int_{0}^{T}X_{t}X_{t}^{\prime}dt\bigg),
\end{align*}
so that $\pr\big(\chi_{0}\leq r^{-1}\big)\leq\pr\big\{\lambda_{\min}(\int_{0}^{T}X_{t}X_{t}^{\prime}dt)\leq r^{-1}\big\}\lesssim C_{L}r^{-L}$. The proof is complete.

\subsection{Proof of Theorem \ref{th.msc}} \label{pr.msc}

We only prove (\ref{msc1.1}) because (\ref{msc1.2}) can be handled analogously. We basically follow Fasen and Kimmig \cite{FasKim15}.

\medskip

(i) Let $\Theta_{m_{0}}$ be nested in $\Theta_{m}$ ($p_{m_{0}}<p_{m}$). 
Define the map $f:\Theta_{m_{0}}\to\Theta_{m}$ by $f(\theta_{m_{0}})=F\theta_{m_{0}}+c$, where $F$ and $c$ satisfy that $\mbbh_{m_{0},n}(\theta_{m_{0}})=\mbbh_{m,n}\big(f(\theta_{m_{0}})\big)$ for any $\theta_{m_{0}}\in\Theta_{m_{0}}$.
If $f(\theta_{m_{0},0})\neq \theta_{m,0}$, $\mbbh_{m_{0},0}(\theta_{m_{0},0})=\mbbh_{m,0}\big(f(\theta_{m_{0},0})\big)<\mbbh_{m,0}(\theta_{m,0})$ and assumption of the optimal model is not satisfied.
Hence we have $f(\theta_{m_{0},0})=\theta_{m,0}$.

By the Taylor expansion of $\mbbh_{m,n}$ around $\hat{\theta}_{m,n}$, we may write
\begin{align*}
\mbbh_{m_{0},n}(\hat{\theta}_{m_{0},n})&=\mbbh_{m,n}\big(f(\hat{\theta}_{m_{0},n})\big) \\
&=\mbbh_{m,n}(\hat{\theta}_{m,n})-\frac{1}{2}\big(\hat{\theta}_{m,n}-f(\hat{\theta}_{m_{0},n})\big)^{\prime}\big(-\p_{\theta_{m}}^{2}\mbbh_{m,n}(\tilde{\theta}_{m,n})\big)\big(\hat{\theta}_{m,n}-f(\hat{\theta}_{m_{0},n})\big),
\end{align*}
where $\tilde{\theta}_{m,n}\cip\theta_{m,0}$ as $n\to\infty$.
Also, $f(\hat{\theta}_{m_{0},n})-\theta_{m,0}=f(\hat{\theta}_{m_{0},n})-f(\theta_{m,0})=F(\hat{\theta}_{m_{0},n}-\theta_{m_{0},0})$.
Since $a_{n}^{-1}\big(\hat{\theta}_{m,n}-f(\hat{\theta}_{m_{0},n})\big)=a_{n}^{-1}(\hat{\theta}_{m,n}-\theta_{m,0})-Fa_{n}^{-1}(\hat{\theta}_{m_{0},n}-\theta_{m_{0},0})=O_{p}(1)$, $\Gam_{m_{0},0}=O_{p}(1)$ and $\Gam_{m,0}=O_{p}(1)$, we have
\begin{align*}
& \pr\left(\qbic^{(m_{0})}-\qbic^{(m)}<0\right) \\
&=\pr\bigg\{\big(\hat{\theta}_{m,n}-f(\hat{\theta}_{m_{0},n})\big)^{\prime}\big(-\p_{\theta_{m}}^{2}\mbbh_{m,n}(\tilde{\theta}_{m,n})\big)\big(\hat{\theta}_{m,n}-f(\hat{\theta}_{m_{0},n})\big) \\
&\qquad\qquad+\log\det\big(-\p_{\theta_{m_{0}}}^{2}\mbbh_{m_{0},n}(\hat{\theta}_{m_{0},n})\big)
-\log\det\big(-\p_{\theta_{m}}^{2}\mbbh_{m,n}(\hat{\theta}_{m,n})\big)<0\bigg\} \\
&=\pr\bigg[\Big\{a_{n}^{-1}\big(\hat{\theta}_{m,n}-f(\hat{\theta}_{m_{0},n})\big)\Big\}^{\prime}\big(-a_{n}^{2}\p_{\theta_{m}}^{2}\mbbh_{m,n}(\tilde{\theta}_{m,n})\big)\Big\{a_{n}^{-1}\big(\hat{\theta}_{m,n}-f(\hat{\theta}_{m_{0},n})\big)\Big\} \\
&\qquad\qquad+\log\det\big(-a_{n}^{2}\p_{\theta_{m_{0}}}^{2}\mbbh_{m_{0},n}(\hat{\theta}_{m_{0},n})\big)
\nn\\
&\qquad\qquad-\log\det\big(-a_{n}^{2}\p_{\theta_{m}}^{2}\mbbh_{m,n}(\hat{\theta}_{m,n})\big)<p_{m}\log a_{n}^{-2}-p_{m_{0}}\log a_{n}^{-2}\bigg] \\
&=\pr\bigg[\Big\{a_{n}^{-1}\big(\hat{\theta}_{m,n}-f(\hat{\theta}_{m_{0},n})\big)\Big\}^{\prime}\Gam_{m,n}(\tilde{\theta}_{m,n})\Big\{a_{n}^{-1}\big(\hat{\theta}_{m,n}-f(\hat{\theta}_{m_{0},n})\big)\Big\} \\
&\qquad\qquad+\log\det\big(\Gam_{m_{0},n}(\hat{\theta}_{m_{0},n})\big)-\log\det\big(\Gam_{m,n}(\hat{\theta}_{m,n})\big)<(p_{m}-p_{m_{0}})\log a_{n}^{-2}\bigg] \\
&\to 1
\end{align*}
as $n\to\infty$.

\medskip

(ii) Let $\mbbh_{m,0}(\theta_{m})\neq\mbbh_{m_{0},0}(\theta_{m_{0},0})$ for every $\theta_{m}\in\Theta_{m}$. 
Because of (\ref{LLN1}) and the consistency of $\hat{\theta}_{m_{0},n}$ and $\hat{\theta}_{m,n}$, we have
\begin{align*}
a_{n}^{2}\mbbh_{m_{0},n}(\hat{\theta}_{m_{0},n})&=a_{n}^{2}\mbbh_{m_{0},n}(\theta_{m_{0},0})+o_{p}(1)=\mbbh_{m_{0},0}(\theta_{m_{0},0})+o_{p}(1), \\
a_{n}^{2}\mbbh_{m,n}(\hat{\theta}_{m,n})&=a_{n}^{2}\mbbh_{m,n}(\theta_{m,0})+o_{p}(1)=\mbbh_{m,0}(\theta_{m,0})+o_{p}(1).
\end{align*}
Since $\mbbh_{m_{0},0}(\theta_{m_{0},0}) > \mbbh_{m,0}(\theta_{m,0})$ a.s. and $a_{n}^{2}\log a_{n}^{-2}\to 0$,
\begin{align*}
& \pr\left(\qbic^{(m_{0})}-\qbic^{(m)}<0\right) \\
&=\pr\Big\{-2\mbbh_{m_{0},n}(\hat{\theta}_{m_{0},n})+2\mbbh_{m,n}(\hat{\theta}_{m,n})+\log\det\big(-a_{n}^{2}\p_{\theta_{m_{0}}}^{2}\mbbh_{m_{0},n}(\hat{\theta}_{m_{0},n})\big) \\
&\qquad\qquad-\log\det\big(-a_{n}^{2}\p_{\theta_{m}}^{2}\mbbh_{m,n}(\hat{\theta}_{m,n})\big)<(p_{m}-p_{m_{0}})\log a_{n}^{-2}\Big\} \\
&=\pr\left\{a_{n}^{2}\big(\mbbh_{m_{0},n}(\hat{\theta}_{m_{0},n})-\mbbh_{m,n}(\hat{\theta}_{m,n})\big) > o_{p}(1)\right\} \nn\\
&= \pr\left\{ \mbbh_{m_{0},0}(\theta_{m_{0},0})-\mbbh_{m,0}(\theta_{m,0})\big) > 0 \right\} + o(1)\nn\\
&\to 1
\end{align*}
as $n\to\infty$.


\subsection{Proof of Theorem \ref{th.msc2}}

We have
\begin{align}
\pr\left(\qbic^{(m_{1,0},m_{2,0})}-\qbic^{(m_{1},m_{2})}\geq0\right)&\leq\pr\left(\qbic^{(m_{1,0},m_{2,0})}-\qbic^{(m_{1,0},m_{2})}\geq0\right) \notag\\
&\quad+\pr\left(\qbic^{(m_{1,0},m_{2})}-\qbic^{(m_{1},m_{2})}\geq0\right).
\end{align}
Applying the proof of Theorem \ref{th.msc} (i) we see that the both terms in the right-hand side tends to zero, hence the claim.


\subsection{Proof of Theorem \ref{th.msc3}}

As with Theorem \ref{th.msc} (i), under assumptions of Theorem \ref{th.msc3} we can deduce that
\begin{align}
&\pr\left(\qbic^{(m_{1,0})}-\qbic^{(m_{1})}<0\right)\to 1,
\end{align}
which means that $\pr(m_{1,n}^{\ast}=m_{1,0})\to 1$.
This together with Theorem \ref{th.msc}(i) then gives
\begin{align*}
& \pr\left(\qbic^{(m_{2,0}|m_{1,n}^{\ast})}-\qbic^{(m_{2}|m_{1,n}^{\ast})}\geq0\right) \nn\\
&=\pr\left(\qbic^{(m_{2,0}|m_{1,n}^{\ast})}-\qbic^{(m_{2}|m_{1,n}^{\ast})}\geq0,~m_{1,n}^{\ast}=m_{1,0}\right) \\
&\qquad+\pr\left(\qbic^{(m_{2,0}|m_{1,n}^{\ast})}-\qbic^{(m_{2}|m_{1,n}^{\ast})}\geq0,~m_{1,n}^{\ast}\neq m_{1,0}\right) \\
&\leq\pr\left(\qbic^{(m_{2,0}|m_{1,0})}-\qbic^{(m_{2}|m_{1,0})}\geq0\right)+\pr\left(m_{1,n}^{\ast}\neq m_{1,0}\right) \to 0,
\end{align*}
completing the proof.


\bigskip

\subsection*{Acknowledgements}
The authors are grateful to Prof. Yoshinori Kawasaki for drawing authors' attention to some relevant literature in econometrics. 
They also thank Prof. Masayuki Uchida and Prof. Nakahiro Yoshida for their valuable comments. 
This work was partly supported by CREST, JST.

\bigskip 

\bibliographystyle{abbrv} 

\end{document}